\documentclass{article}

\usepackage[margin=1in, top=1.5in, bottom=1.5in, left=1.2in, right=1.2in]{geometry} 
\usepackage{graphicx,amsmath,amsfonts,amsthm,amssymb,float,bm,epsfig,color,enumitem}
\usepackage[colorlinks=true,linkcolor=blue,citecolor=red]{hyperref}
\usepackage{subfigure}
\usepackage{tikz}
\usepackage{authblk} 
\usepackage{sectsty} 
\usepackage[utf8]{inputenc}   

\DeclareRobustCommand{\orcidicon}{%
    \begin{tikzpicture}
    \draw[lime, fill=lime] (0,0) 
    circle [radius=0.16] 
    node[white] {{\fontfamily{qag}\selectfont \tiny ID}};
    \draw[white, fill=white] (-0.0625,0.095) 
    circle [radius=0.007];
    \end{tikzpicture}
    \hspace{-2mm}
}

\newcommand{\orcidauthorA}{0009-0007-6656-430X}

\newcommand{\orcidA}{\href{https://orcid.org/\orcidauthorA}{\orcidicon}}
\newcommand{\orcidB}{\href{https://orcid.org/\orcidauthorA}{\orcidicon}}
\newcommand{\orcidC}{\href{https://orcid.org/\orcidauthorA}{\orcidicon}}

\newcommand{\email}[1]{\texttt{\small \href{mailto:#1}{#1}}} 
\newcommand{\affilnorm}[1]{\normalsize #1} 

\sectionfont{\fontfamily{phv}\selectfont} 

\newtheorem{theorem}{Theorem} 
\newtheorem{remark}[theorem]{Remark}

\begin{document}

\title{Turing pattern induced by cross-diffusion in cancer immunotherapy model}

\author[1]{S. Bonfiglio\orcidA\footnote{\email{simone.bonfiglio1@studenti.unime.it}}}
\author[1]{C. F. Munaf\'o\orcidB\footnote{\email{carmelofilippo.munafo@unime.it}}}
\author[1]{P. Rogolino\orcidC\footnote{Corresponding author: \email{progolino@unime.it}} }


\affil[1]{\affilnorm{Department of Mathematical and Computer Sciences, Physical Sciences and Earth Sciences, University of Messina, Viale F. Stagno d'Alcontres 31, 98166, Messina, Italy}}

\date{}

\maketitle

\begin{abstract}
	In this paper, we investigate a mathematical model describing the interactions between effector cells ($E$), cancer cells ($T$), and the IL-2 compound ($I_L$). The model considered here is a generalization, taking into account some cross-diffusion effects, of a spatial cancer immunotherapy model proposed by S. Suddin et al in 2021. These modifications allow us to describe two biologically relevant scenarios: a patient treated with Adoptive Cell Immunotherapy (ACI) and a patient not receiving any treatment/therapy. 
	
	Cross-diffusion effects are particularly relevant in the interactions between tumor cells and the immune system, in fact they play a key role in immune response dynamics and cannot be neglected. We analyze the equilibrium points of the homogeneous system, along with their stability and bifurcation mechanisms. Furthermore, adopting the Turing approach for reaction–diffusion systems, we investigate the diffusion-driven instability and the emergence of spatial regular structures (stationary in time), \emph{i.e.} the patterns.
	
	Finally, numerical simulations based on the Finite Difference Method (FDM) are presented for the two previously mentioned scenarios.
	
\end{abstract}

\section{Introduction}
The importance of immune system in limiting cancer progression has recently been highlighted through several direct molecular interactions between cancers and immune effector cells. In these interactions, the interleukin-2 (IL-2) plays a fundamental role in the activation of immune system and has been shown to be effective in mediating tumor regression. As a monotherapy, IL-2 has received approval for the treatment of metastatic renal cell carcinoma and metastatic melanoma.

Cervical cancer is a malignant tumor whose pathogenesis is caused by the human papilloma virus (HPV), a retrovirus of DNA family, which consists of nine open reading frames (ORFs): from E1 to E7, L1 and L2. Each protein has a specific role related to viral infection and replication processes. HPV infection causes lesions in cervical tissue without forming clots. Several treatment options are currently available for patients with cervical cancer, including chemotherapy, radiotherapy, and surgery. Although these therapies remain important in clinical context, none of these methods have proven to be completely effective. 

Hence, in addition to preventive strategies, there is a growing need to develop more successful and effective therapeutic strategies. One promising avenue is the immunotherapy, which has been extensively investigated in recent years \cite{curti1996influence,gause1996phase,hara1996rejection,kaempfer1996prediction,keilholz1994immunotherapy,valle2016pleiotropic}. 

In the last decades, among emerging immune-therapeutic treatment for cervical cancer, Adoptive Cellular Immunotherapy (ACI) has been introduced as a novel and promising strategy. ACI involves the injection of cultured immune cells with anti-tumor properties into the patient, typically in combination with high doses of IL-2 \cite{curti1996influence}. In fact, it has been shown that cervical cancer cells are detected as antigens by the immune system. A component of the immune system that plays a central role in attacking antigens is the Cytotoxic T lymphocytes (CD8+) that can be activated by CD4+ T cells. These responses are enhanced by some IL-2 compound, which is produced in higher quantities by CD4+ T cells and to a lesser quantities by CD8+ T cells. 

The ACI treatment it is usually carried out in conjunction with large amounts of IL-2. This can take two approaches, lymphokine-activated killer (LAK) cell therapy or tumor-infiltrating lymphocyte (TIL) therapy, see \cite{kirschner1998modeling}. 

Understanding the mechanisms behind these immune responses and their impact on tumor dynamics has motivated the development, in recent decades, of various mathematical models. These models have been developed to describe the interactions between tumor cells and the corresponding immune responses \cite{adam1997survey,bagheri2021,wang2021turing,wangberg2023}. It can be useful for clinicians in optimizing both single and combination therapies (\emph{i.e. in designing optimal treatment protocols}), and can provide researchers understand the mechanisms underlying the effectiveness or ineffectiveness of therapeutic combinations. Therefore, in an important field such as cancer immunotherapy, mathematical modeling can play a significant role in leading the future direction of research.

Among the various approaches proposed in the literature to describe how cell density evolves in space and time due to local interactions and diffusion, reaction-diffusion models are particularly prominent \cite{turing1990chemical,Murray2,LamLou}. These systems, taking into account both self and cross-diffusions, are widely applied across various scientific fields, such as ecology \cite{kan1998singular,sherratt2013pattern,
carfora2025turing,ali2025turing,ali2026turing,tulumello2014cross}, biology and medicine \cite{giunta2021pattern,lombardo2017demyelination}, epidemiology \cite{wang2012complex,avila2022dynamics,della2022mathematical,della2025modelling}, and the social sciences \cite{inferrera2024reaction,petrovskii2020modelling,wen2009global}.

The emergence of spatially heterogeneous structures via Turing-type instabilities in cancer–immune interaction models has been explored by several authors. 
The studies conducted have highlighted how spatial diffusion and instability mechanisms can significantly influence cancer dynamics and treatment outcomes. 
Ko and Ahn \cite{ko2011stationary} investigated a reaction-diffusion model incorporating effector cells, tumor cells, and IL-2 under immunotherapy in a spatially inhomogeneous in vivo environment, demonstrating the influence of spatial heterogeneity on therapeutic efficacy. Similarly, Suddin et al. \cite{suddin2021reaction} analyzed a related model and found that high diffusivity of effector cells facilitates the emergence of spatially heterogeneous steady states, while excluding the existence of stable spatially homogeneous positive (tumor-coexistence) equilibria. Brennan in \cite{brennan2025pattern} further explored the role of subcritical Turing instabilities, showing that even when the cancer-free state is linearly stable, stable spatially inhomogeneous solutions can persist, allowing localized tumor survival. 
This emphasizes how spatial structure at the tissue level can significantly affect tumor resistance to immunotherapy, particularly in models that incorporate IL-2 and effector cell dynamics as treatment parameters. Moreover, Oluwatosin in \cite{oluwatosin2022investigating} investigated Turing pattern formation in a spatially distributed cancer–immune cell interaction model, deriving the conditions under which such patterns emerge due to diffusion-driven instability. Similarly, Wang in \cite{wang2021turing,wangberg2023} studied a comparable reaction–diffusion system and also derived conditions for Turing instability, incorporating diffusion terms to model the spatial spread of tumor cells (specifically myelogenous cells) and immune cells within the circulating blood.

The role of spatial structure in the emergence of Turing instabilities has also been investigated in the context of network-organized systems \cite{othmer1969interactions,othmer1971instability,othmer1974non,nakao2010turing,hata2014dispersal,sarker2025spatial}. Classical models focus on chemical mixtures in rods or sheets with constrained diffusion, offering a basis for analyzing spatial pattern formation. In parallel, compartmental network models, initiated by Othmer and Scriven \cite{	othmer1969interactions,othmer1971instability,othmer1974non}, have shown that diffusion-driven instabilities can also emerge in discrete structures, where network topology plays a key role in determining system stability.

Our aim is to develop a mathematical model describing the immunotherapeutic treatment of cancer cells, based on a reaction-diffusion system, which can exhibit diffusion-driven Turing instability. This mechanism, where small perturbations destabilize a spatially homogeneous state, leading to the emergence of spatial patterns, may play  a crucial role in identifying the conditions under which treatment resistance can arise. This aspect is particularly relevant when simulating therapeutic scenarios that vary in both time and space.
 
To this end, we start from a model proposed by Suddin in 2021 \cite{suddin2021reaction}, which describes the interactions between cancer cells, effector cells and IL-2 compounds. In this model, the population is divided into three subgroups, whose spatial distributions evolve over time. This predicts various kinds of dynamical outcomes, with potential implications for the spread of the cancer cells within the tissue.

The main novelty with respect to the Suddin model, \cite{suddin2021reaction}, relies on the introduction of a cross-diffusive contribution describing the spread of the IL-2 compounds towards the  cancer cells. In particular, cross-diffusion plays a fundamental role in the dynamics of the immune response and cannot be neglected.
Based on this extended model, we examine two scenarios, one describing tumor growth behavior in untreated patients, and the other in those receiving ACI therapy untreated patients. In both cases, the system admits asymptotically stable homogeneous equilibria, which lose their stability due to self- and cross-diffusive terms, leading to the emergence of stable spatial Turing pattern, as in \cite{wang2021turing,ko2011stationary,brennan2025pattern,oluwatosin2022investigating}.

This paper is organized as follows. 
In Section \ref{sec_2}, we propose our mathematical model and analyze the stability of the homogeneous equilibrium configurations. Furthermore, Section \ref{sec_3} is devoted to the study of Turing instability; including  restrictions on some parameters under which spatial patterns emerge. Section \ref{sec_4} presents numerical simulations illustrating  the formation of these patterns. Finally, Section \ref{sec_conclu} provides a discussion and  concluding remarks.
 
\section{A nonlinear model of cancer immunotherapy}\label{sec_2}
In this section, a formulation of the fundamental assumptions underlying  our model is presented.

Let us consider a section of cervical tissue denoted by  $\Omega$. We assume that the cell population  consists of three compartments, namely  cancer cells,  effector cells and  IL-2 compounds. Let $E(\mathbf{x},t)$, $T(\mathbf{x},t)$ and $I_L(\mathbf{x},t)$ denote their densities at the spatial position $\mathbf{x} \in \Omega$ and time $t>0$, respectively. The interactions and diffusion of these quantities are described in  \cite{suddin2021reaction} by the following set of equations:
\begin{subequations}\label{sys}
    \begin{align}
    	& \frac{\partial E}{\partial t} = c \,T - \mu _1 E + \frac{p_1 E \,I_L}{g_1 + I_L} + s_1 + D_{11} \Delta E, & \forall \, \mathbf{x} \in \Omega, \, \forall t>0, \\
	& \frac{\partial T}{\partial t} = r_2 T (1 - b T) - \frac{p_2 E \,T}{g_2 + T} + D_{22} \Delta T, & \forall \, \mathbf{x} \in \Omega, \, \forall t>0, \\
	& \frac{\partial I_L}{\partial t} = \frac{p_3 E \,T}{g_3 + T} - \mu _3 I_L + s_3 + D_{33} \Delta I_L, & \forall \, \mathbf{x} \in \Omega, \, \forall t>0,
    \end{align}
\end{subequations}
where all the parameters involved are assumed to be positive, except for the parameters $s_1$ and $s_3$ which are assumed to be non-negative. 

Specifically, the coefficient $c$ represents the antigenicity of the cancer cells and measures the ability of the immune system to recognize cancer through non-self protein antigens. The parameter $\mu_1$ denotes the average of the natural lifespan of  effector cells; while $\mu_3$ represents the death rate of IL-2 compounds. The parameter $s_1$ accounts for Adoptive Cellular Immunotherapy (ACI) to the effector cells; and $s_3$ is cytokine therapy to increase IL-2 compounds. The parameters $p_1$, $p_2$ and $p_3$ denote, respectively, the maximum proliferation rate of effector cells induced by IL-2 compounds, the maximum degradation rate of cancer cells that have an interaction with the effector cells, and the maximum production rates of the IL-2 compounds by the effector cells. The parameters $g_1$, $g_2$ and $g_3$ represent, respectively, the proliferation kinetics of the effector cells, the degradation kinetics of the cancer cells, and the growth of the IL-2 compounds caused by the interaction between  cancer and  effector cells. Finally, the constants $D_{ii}$ with $i \in \{1,2,3\}$ represent the self diffusion coefficients of the three cell groups.
 
We observe that  cancer growth is supposed to follow a logistic model, in which $r_2$ represents the constant birth rate and $1/b$ denotes the carrying capacity.

In \cite{suddin2021reaction}, an analytical study has been provided for the equilibria of the reaction-diffusion system \eqref{sys} describing the interactions mentioned above. Additionally,  after analyzing the initial behavior of the tissue, numerical simulations have been conducted to illustrate the spatial structures emerging within the tissue,  representing the effects of immunotherapy, as well as the steady-state conditions of the system which exhibit the long-term behavior of these interactions.

Here, we propose a variant of the model \eqref{sys} by introducing some cross-diffusion effects, as  follows:
\begin{subequations}\label{sys_mod}
	\begin{align}
		& \frac{\partial E}{\partial t} = c \,T - \mu _1 E + \frac{p_1 E \,I_L}{g_1 + I_L} + s_1 + D_{11} \Delta E, & \forall \, \mathbf{x} \in \Omega, \, \forall t>0, \label{sys_mod_a}\\
	   	& \frac{\partial T}{\partial t} = r_2 T (1 - b T) - \frac{p_2 E\, T}{g_2 + T} + D_{22} \Delta T, & \forall \, \mathbf{x} \in \Omega, \, \forall t>0, \label{sys_mod_b} \\
	   	& \frac{\partial I_L}{\partial t} = \frac{p_3 E \,T}{g_3 + T} - \mu _3 I_L + s_3 + D_{33} \Delta I_L + D_{32} \Delta T, & \forall \, \mathbf{x} \in \Omega, \, \forall t>0. \label{sys_mod_c}
    	\end{align}
\end{subequations}
wherein $\Delta = \frac{\partial^2}{\partial^2 x}+\frac{\partial^2}{\partial^2 y}$ is the Laplacian operator in  two-dimensional space, which
describes the random movement (diffusion) of molecules within the cell. In particular, in \eqref{sys_mod}, when $i=j$, the positive coefficients $D_{ij}$ represent self-diffusion, whereas when $i \neq j$, $D_{ij}$ represent cross-diffusion effects.
In general, the positive/negative cross-diffusion terms imply that a species of cells tends to diffuse toward regions with lower/greater concentration of the other cell type. In this context, the term $D_{32}\Delta T$ with a positive value of the coefficient $D_{32}$ indicates that IL-2 cells move toward regions far from the tumor cells, whereas a negative value indicates that IL-2 cells migrate toward zones with a high concentration of tumor cells.

These effects are particularly relevant in the interactions between tumor cells and the immune system and therefore cannot be neglected \cite{valle2016pleiotropic}.

The model will be studied in a two-dimensional rectangular domain $\Omega=[0,L_x]\times [0,L_y] \subset \mathbb{R}^2$, where $L_x$ and $L_y$ denote its dimensions. We assign the following initial conditions
\begin{align*}
	E (\mathbf{x},0) &= E_0 \, f_E(\mathbf{x}), & \forall \, \mathbf{x} \in \Omega, \\ 
	T (\mathbf{x},0) &= T_0 \, f_T(\mathbf{x}), & \forall \, \mathbf{x} \in \Omega, \\
	I_L(\mathbf{x},0) &= I_L^0 \, f_I(\mathbf{x}), & \forall \, \mathbf{x} \in \Omega,
\end{align*}
and boundary conditions
\begin{align*}
	\nabla E(\mathbf{x},t) \cdot \mathbf{\hat{n}} &= 0, & \forall \, \mathbf{x} \in \partial \Omega, \, \forall t>0, \\
	\nabla T(\mathbf{x},t) \cdot \mathbf{\hat{n}} &= 0, & \forall \, \mathbf{x} \in \partial \Omega, \, \forall t>0, \\
	\nabla I_L(\mathbf{x},t) \cdot \mathbf{\hat{n}} &= 0, & \forall \, \mathbf{x} \in \partial \Omega, \, \forall t>0,
\end{align*}
where $f_E$, $f_T$, $f_I$ are suitable functions, $\partial\Omega$ and $\mathbf{\hat{n}}$ being the boundary of $\Omega$ and the outward unit normal to $\partial\Omega$, respectively. The zero-flux boundary conditions mean that no effector cells, cancer cells or IL-2 compounds cross the tissue boundary. 

For the sake of convenience, by applying the following scalings
\begin{equation*}
	\hat{t} = \frac{t}{\tau} \qquad \hat{x} = \frac{x}{L_x}, \qquad \hat{y} = \frac{y}{L_y}, \qquad u = \frac{E}{E_0}, \qquad v = \frac{T}{T_0}, \qquad w = \frac{I_L}{I_L^0},
\end{equation*}
suggested by the problem, the system \eqref{sys_mod} is simplified into its dimensionless form as follows
\begin{subequations}\label{eq:3}
	\begin{align}
		& \frac{\partial u}{\partial t} = \hat{c} v - \hat{\mu _1} u + \frac{\hat{p_1} u w}{\hat{g_1} + w} + \hat{s_1} + d_{11} \left( \frac{\partial^2 u}{\partial \hat{x}^2} + \tau_L \frac{\partial^2 u}{\partial \hat{y}^2} \right), \\
        		& \frac{\partial v}{\partial t} = \hat{r_2} v (1 - \hat{b} v) - \frac{\hat{p_2} u v}{\hat{g_2} + v} + d_{22} \left( \frac{\partial^2 v}{\partial \hat{x}^2} + \tau_L \frac{\partial^2 v}{\partial \hat{y}^2} \right), \\
        		& \frac{\partial w}{\partial t} = \frac{\hat{p_3} u v}{\hat{g_3} +v}- \hat{\mu _3} w + \hat{s_3} + d_{33} \left( \frac{\partial^2 w}{\partial \hat{x}^2} + \tau_L \frac{\partial^2 w}{\partial \hat{y}^2} \right) \\
		& \qquad \quad+ d_{32} \left( \frac{\partial^2 v}{\partial \hat{x}^2} + \tau_L \frac{\partial^2 v}{\partial \hat{y}^2} \right) \nonumber
	\end{align}
\end{subequations}
The initial and the boundary conditions of system \eqref{eq:3} become
\begin{align*}
	u(\mathbf{x},0) &= \hat{f}_E(\mathbf{x}), & \forall \, \mathbf{x} \in \hat\Omega, \\ 
	v(\mathbf{x},0) &= \hat{f}_T(\mathbf{x}), & \forall \, \mathbf{x} \in \hat\Omega, \\
	w(\mathbf{x},0) &= \hat{f}_I(\mathbf{x}), & \forall \, \mathbf{x} \in \hat\Omega,
\end{align*}
and
\begin{align*}
    \nabla u(\mathbf{\hat{x}},t) \cdot \mathbf{\hat{n}} &= 0, & \forall \, \mathbf{x} \in \partial \hat\Omega, \, \forall t>0, \\
    \nabla v(\mathbf{x},t) \cdot \mathbf{\hat{n}} &= 0, & \forall \, \mathbf{x} \in \partial \hat\Omega, \, \forall t>0, \\
    \nabla w(\mathbf{x},t) \cdot \mathbf{\hat{n}} &= 0, & \forall \, \mathbf{x} \in \partial \hat\Omega, \, \forall t>0,
\end{align*}
respectively, wherein $\hat\Omega = [0,1]\times[0,1]$ is the dimensionless domain. The dimensionless parameters assume the following expressions
\begin{align*}
	& \hat{p}_1 = \tau p_1, \qquad \hat{p}_2 = \frac{\tau p_2 E_0}{T_0}, \qquad \hat{p}_3 = \frac{\tau p_3 E_0}{I_L^0}, \qquad  \hat{g}_1 = \frac{g_1}{I_L^0}, \qquad \hat{g}_2 = \frac{g_2}{T_0}, \\
	& \hat{g}_3 = \frac{g_3}{T_0}, \qquad \hat{c} = \frac{\tau c T_0}{E_0}, \qquad \hat{\mu}_1 = \tau \mu_1, \qquad \hat{\mu}_3 = \tau \mu_3, \qquad  \hat{s}_1 = \frac{\tau s_1}{E_0}, \\
	& \hat{s}_3 = \frac{\tau s_3}{I_L^0}, \qquad  \hat{r}_2 = \tau r_2, \qquad  \hat{b} = b T_0, \qquad  \hat{d}_{11} = \frac{\tau D_{11}}{L_x^2} , \qquad \hat{d}_{22} = \frac{\tau D_{22}}{L_x^2} , \\ 
	& \hat{d}_{33} = \frac{\tau D_{33}}{L_x^2}, \qquad \hat{d}_{32} = \frac{\tau D_{32} T_0}{L_x^2 I_L^0}
\end{align*}
and $\tau_L = \left(\dfrac{L_x}{L_y}\right)^2 $ is a scaling factor.

In order to simplify the notation let us drop the hats and we set $L_x=L_y$ (thus $\tau_L=1$), whence, introducing the vectors $\mathbf{U}$ and the matrix $\mathcal{D}$, say
\begin{equation*}
    \mathbf{U} = 
    \begin{pmatrix}
        u \\
        v \\
        w
    \end{pmatrix},
    \qquad
        \mathcal{D} =
    \begin{pmatrix}
        d_{11} & 0      & 0\\
        0      & d_{22} & 0\\
        0      & d_{32} & d_{33}
    \end{pmatrix},
\end{equation*}
we may write dimensionless model \eqref{eq:3} in compact form as
\begin{equation}\label{syst_comp}
	\frac{\partial \mathbf{U}}{\partial t} = \mathbf{F}(\mathbf{U}) + \mathcal{D} \,\Delta \mathbf{U}.
\end{equation}
where
\begin{equation*}
\mathbf{F}(\mathbf{U})=
	\begin{pmatrix}
		c v - \mu_1 u + \dfrac{p_1 u\, w}{g_1 + w} + s_1 \\
        		r_2 v (1 - b v) - \dfrac{p_2 u\, v}{g_2 + v} \\
        		\dfrac{p_3 u \,v}{g_3 + v} - \mu_3 w + s_3
   	\end{pmatrix}
\end{equation*}

\subsection{Linear stability analysis without diffusion}
In this section, we establish the existence of steady-state points of the system and analyze their local stability. By considering the system \eqref{syst_comp}, where the spatial terms are neglected
\begin{equation}\label{sys_adiff}
	\frac{ \text{d} \mathbf{U}}{\text{d} t} = \mathbf{F}(\mathbf{U}).
\end{equation}
It can be easily shown that there are at least six admissible equilibria. However, the focus of our analysis will be on the existence and stability of two particular equilibria:
\begin{itemize}
	\item $\mathbf{U}_1^*=(u_1^*,0,w_1^*)$ commonly known as \textbf{Cancer Free Equilibrium} (CFE), since cancer cells are absent,
    	\item $\mathbf{U}_2^*=(u_2^*,v_2^*,w_2^*)$ referred to as \textbf{Cancer Coexistence Equilibrium} (CCE), corresponds to a state in which cancer cells, immune cells 
 and IL-2  coexist.
\end{itemize} 
where the quantities $u_1^*$, $w_1^*$, $u_2^*$, $v_2^*$, $w_2^*$ must take positive values.

The focus on only two of the six equilibria that the system mathematically admits is due to their the biological meaning. In particular, the CFE, which corresponds to the absence of tumor cells (\emph{i.e.} $v=0$), represents a scenario of healing or control of the tumor under the action of the immune system. For this reason it is necessary to verify its stability, since a stable CFE indicates that the system tends to evolve towards a tumor free state. Similarly, the CCE, where both  tumor and immune cells persist, corresponds to a condition in which the tumor is not eradicated but remains in a controlled, non-aggressive state. A stable CCE suggests a persistent not contained tumor that does not progress toward  more severe scenarios requiring additional treatment or surgery, when possible.
\begin{theorem}
	The existence of a positive constant steady state, corresponding to the CFE, is guaranteed if and only if 
	\begin{equation*}
    		u_1^* = - \frac{s_1(s_3 + g_1 \mu _3)}{p_1 s_3 - s_3 \mu _1 - g_1 \mu _1 \mu _3 }, \qquad w_1^* = \frac{s_3}{\mu _3},
	\end{equation*}
	under the following constraint:
	\begin{equation*}
    		p_1 s_3 - s_3 \mu _1 - g_1 \mu _1 \mu _3 < 0.
	\end{equation*}

	\begin{proof}
 		To ensure that the CFE exists we set $v^*=0$ in the condition $\mathbf{F}(\mathbf{U}^*)=\mathbf{0}$ obtaining
		\begin{align*}
			- \mu_1 u^* + \frac{p_1 u^*\, w^*}{g_1 + w^*} + s_1 &= 0  \\
        			- \mu_3 w^* + s_3 &= 0
		\end{align*}
		from which
 		\begin{equation*}
    			u_1^* = - \frac{s_1(s_3 + g_1 \mu _3)}{p_1 s_3 - s_3 \mu _1 - g_1 \mu _1 \mu _3 }, \qquad w_1^* = \frac{s_3}{\mu _3}.
		\end{equation*}
		For this equilibrium to be biologically meaningful, $u_1^*$ must be positive. Since the numerator is always positive, this holds if and only if the denominator satisfies the following constraint:
		\begin{equation*}
    			p_1 s_3 - s_3 \mu _1 - g_1 \mu _1 \mu _3 < 0.
		\end{equation*}
	\end{proof}
\end{theorem}

\begin{theorem}
	The CFE equilibrium point of the system \eqref{syst_comp} is locally asymptotically stable (LAS) if and only if the condition $u_1^*>\dfrac{r_2 g_2}{p_2}$ holds.
    	\begin{proof}
		The Jacobian matrix concerned with the linearization for model \eqref{sys_adiff} at CFE is written as 
		\begin{equation*}
         		\mathcal{J} =
        			\begin{bmatrix}
            			\dfrac{p_1 w}{g_1+w}-\mu _1 &  c &  -\dfrac{p_1 u w}{(g_1 +w)^2} \\ \\
            			- \dfrac{p_2 u}{g_2+v} &  -b r_2 v +r_2 (1-b v)+ \dfrac{p_2 u}{g_2+v}-\dfrac{p_2 u v}{(g_2+v)^2} &  0 \\ \\
           			 \dfrac{p_3 v}{g_3+v} &  -\dfrac{p_3 u v}{(g_3+v)^2}+\dfrac{p_3 u}{g_3+v} &  \mu_3
       			 \end{bmatrix}.
  		\end{equation*} 
		Evaluating the Jacobian at $\mathbf{U}_1^*$, \emph{i.e.} $\mathcal{J}\vert_{\mathbf{U}_1^*} = \nabla_{\mathbf{U}} \mathbf{F}\vert_{\mathbf{U}_1^*} $, we find that its eigenvalues are
    		\begin{equation*}
            		\lambda _1  =-\mu_3, \qquad
            		\lambda _2  =\frac{p_1 s_3 - s_3 \mu _1 - g_1 \mu _1 \mu _3}{s_3 + g_1 \mu _3}, \qquad
            		\lambda _3 =r_2-\frac{p_2}{g_2} u_1^*,
       	 	\end{equation*}
  		given that all treatment parameters applied to the system are non-negative, we deduce that $\mathbf{U}_1^*$ is LAS if 
       		\begin{equation*}
            		u_1 ^* > \frac{r_2 g_2}{p_2}.
        		\end{equation*}
	\end{proof}
\end{theorem}

\begin{remark}
	It easy to deduce that, in the absence of therapy, \emph{i.e.} $s_1=s_3=0$, the CFE, $\mathbf{U}_1^*=(u_1^*,0,w_1^*)$, reduces to $(0,0,0)$. Moreover, this equilibrium point is unstable.
\end{remark}

\begin{theorem}
	The existence of a positive constant steady state corresponding to the CCE is ensured if and only if  
	\begin{equation*}
		u_2^* = \frac{r_2 (g_2 + v_2^*) (1 - b v_2^*)}{p_2}, \qquad w_2^* = \frac{s_3 p_2 (g_3 + v_2^*) + p_3 r_2 v_2^* (g_2 + v_2^*) (1 - b v_2^*)}{p_2 (g_3 + v_2^*) \mu _3},
	\end{equation*}
	under the following conditions: 
	\begin{itemize}
		\item[1)]  $v_2^*<1/b $
		\item[2)]  $a_5 v^5 + a_4 v^4 + a_3 v^3 + a_2 v^2 + a_1 v + a_0 = 0$
	\end{itemize}
	where $1/b$ is the carrying capacity of the tumor cells and the coefficients are given by
	\begin{align*}
   		a_5 & = b^2 p_3 r_2^2 (p_1 - \mu _1), \\
    		a_4 & = -b c p_2 p_3 r_2 - 2 b p_3 r_2^2 (p_1 - \mu_ 1) + 2 b^2 g_2 p_3 r_2^2 (p_1 - \mu _1), \\
    		a_3 &= c p_2 p_3 r_2 - b c g_2 p_2 p_3 r_2 - b  p_2  p_3  r_2  s_1 - b p_1  p_2  r_2  s_3 + p_3 r_2^2 (p_1 - \mu _1)  \\
    		& \quad  - 4 b g_2 p_3 r_2^2 (p_1 - \mu _1) + b^2 g_2^2 p_3 r_2^2 (p_1 - \mu _1) + b  p_2 r_2 \mu _1(s_3 + g_1 \mu _3), \\
    		a_2 & = c g_2 p_2 p_3 r_2 + p_2  p_3  r_2  s_1 - b g_2 p_2 p_3 r_2 s_1 + p_1 p_2 r_2 s_3 - b g_2 p_1  p_2  r_2  s_3  \\
    		& \quad + 2 g_2 p_3 r_2^2 (p_1 - \mu _1) - 2 b g_2^2 p_3 r_2^2 (p_1 - \mu _1) + c  p_2^2 (s_3+g_1 \mu _3) \\
    		& \quad  - p_2 r_2 \mu _1 (s_3 + g_1 \mu _3)+ b g_2 p_2 r_2 \mu _1 (s_3 + g_1 \mu _3) + b g_3 p_2 r_2 (p_1 s_3 - \mu _1 (s_3 + g_1 \mu _3)), \\
    		a_1 & = g_2^2 p_3 r_2^2 (p_1 - \mu _1) + c  g_3 p_2^2(s_3 + g_1 \mu _3) - g_2 p_2 r_2 \mu _1 (s_3 + g_1 \mu _3) + g_3 p_2 r_2 (p_1 s_3 \\
    		& \quad - \mu _1 (s_3 + g_1 \mu _3) ) + b  g_2 g_3 p_2 r_2(-p_1 s_3 + \mu _1 (s_3 + g_1 \mu _3)), \\
    		a_0 &= g_3 p_2^2 s_1 (s_3 + g_1 \mu _3) - g_2 g_3 p_2 r_2 (-p_1 s_3 + \mu _1 (s_3 + g_1  \mu _3)).
	\end{align*}
	
	\begin{proof}
		By solving the algebraic system $\mathbf{F}(\mathbf{U}^*)=\mathbf{0}$ we obtain the expressions
    		\begin{equation*}
			u_2^* = \frac{r_2 (g_2 + v_2^*) (1 - b v_2^*)}{p_2}, \qquad w_2^* = \frac{s_3 p_2 (g_3 + v_2^*) + p_3 r_2 v_2^* (g_2 + v_2^*) (1 - b v_2^*)}{p_2 (g_3 + v_2^*) \mu _3}, 
    		\end{equation*}
		Since all parameters are assumed to be positive, the positivity of $u_2^*$ and $w_2^*$ requires that the factor $(1 - b v_2^*)$ be positive. Hence, the condition 1), \emph{i.e.} $v_2^*<1/b$ must hold. Moreover, substituting the expressions for $u_2^*$ and $w_2^*$, the system reduces to a single equation in the unknown $v^*_2$ (which we denote simply by $v$) given by
		\begin{equation}\label{eq_v}
			\frac{a_5 v^5 + a_4 v^4 + a_3 v^3 + a_2 v^2 + a_1 v + a_0}{p_2 \left[ p_3 r_2 v (g_2 + v) (1 - b v) + g_3 p_2 (s_3 + g_1 \mu _3) + p_2 v (s_3 + g_1  \mu _3) \right]} = 0
		\end{equation}
		Since the denominator of \eqref{eq_v} is strictly positive and never vanishes, the equation is satisfied if and only if
		\begin{equation}\label{eq_v_new}
    			a_5 v^5 + a_4 v^4 + a_3 v^3 + a_2 v^2 + a_1 v + a_0 = 0.
		\end{equation}
	\end{proof}
\end{theorem}

Since finding positive solutions to the fifth-degree equation \eqref{eq_v_new} is complex, we limit our analysis to two significant cases in order to examine tumor progression and evaluate the effectiveness of ACI therapy in both treated and untreated patients.

In the following subsections, we therefore investigate the existence of a real positive root, $v_2^*$, of  equation \eqref{eq_v_new}, considering both for the untreated patient (\emph{i.e.} $s_1=s_3=0$) and for the patient undergoing ACI therapy (\emph{i.e.} $s_1,s_3 \neq 0$). Let us remark that, under the assumption $p_1=\mu_1$ equation \eqref{eq_v_new} reduces to a fourth-degree equation. For a detailed discussion on the existence of equilibria in this case, see \cite{KusSan}.

\subsubsection{Stability analysis for untreated patients}
In this case, assuming $s_1=s_3=0$, we consider $p_2$ and $c$ as variable parameters, while adopting the remaining parameters values from \cite{kirschner1998modeling}, as summarized in Table \ref{tab:par}. We then analyze the existence and stability of CCE as $p_2$ and $c$ vary.

\begin{table}[ht]
	\centering
    	\begin{tabular}{|c|c|c|c|}
        		\hline
        		Parameter & Value & Parameter & Value \\ \hline
        		$\mu_1$ & 0.167 & $p_3$ & 27.778 \\ \hline
        		$p_1$  & 0.69167 & $g_3$ & 0.001 \\ \hline
        		$g_1$ & 20 & $\mu_3$ & 55.55556 \\ \hline
        		$r_2$ & 1 & $b$ & 1 \\ \hline
        		$g_2$ & 0.1 & & \\ \hline
    	\end{tabular}
	\caption{Parameter values}
	\label{tab:par}
\end{table}
The coefficients of equation \eqref{eq_v_new} become as follows
\begin{align*}
	& a_5 = 14.57, \\
    	& a_4 = (-26.23 - 27.78 c \,p_2), \\
    	& a_3 = (8.89 + 185.56 p_2 + 25 c \,p_2), \\
    	& a_2 = (2.62 - 166.81 p_2 + 2.78 c \,p_2 + 1111.11 c \,p_2^2), \\
	& a_1 = (0.15 - 18.72 p_2 + 1.11 c \,p_2^2), \\
    	& a_0 = -0.02 p_2.
\end{align*}
Using the Routh-Hurwitz criterion, we investigate the minimal conditions required to ensure that the real root of equation \eqref{eq_v_new}, which exists because the polynomial has odd degree, is positive. Since $a_5>0$, $a_4<0$, $a_3>0$ and $a_0<0$ for every value of $p_2$, our goal is achieved if one of the following sets of conditions fulfilled
\begin{itemize}
	\item $R_4<0$, $R_5<0$, 
	\item $R_4>0$, $R_5>0$, 
\end{itemize}
where $R_4$, $R_5$ are the fourth and the fifth element of the first column of the Routh matrix, respectively.

\subsubsection{Stability analysis for patients treated  with ACI therapy}
In this case, assuming $s_1=0.0035$, $s_3=0.2$, we consider $p_2$ and $c$ as variable parameters and the remaining parameters are fixed as in Table \ref{tab:par}. We investigate the existence and stability of the CCE as $p_2$ and $c$ vary. The coefficients of equation \eqref{eq_v_new} are then given as follows 
\begin{align*}
	& a_5 = 14.57, \\ 
    	& a_4 = (-26.23 - 27.78 c p_2), \\
    	& a_3 = (8.89 + 185.35 p_2 + 25 c p_2), \\
    	& a_2 = (2.62 - 166.63 p_2 + 2.78 c p_2 + 1111.31 c p_2^2), \\
    	& a_1 = (0.15 - 18.7 p_2 +3.89 p_2^2 + 1.11 c p_2^2), \\
    	& a_0 = -0.02 p_2 + 0.004 p_2^2.
\end{align*}
Using Routh-Hurwitz criterion, we investigate the minimal conditions required to ensure that the real root of equation \eqref{eq_v_new}, which exists because the polynomial has odd degree, is positive. Given that $a_5>0$, $a_4<0$ and $a_3>0$ for every value of $p_2$, achieving this requires that one of the following sets of conditions is satisfied 
\begin{itemize}
	\item $a_0<0$, $R_4<0$, $R_5<0$,
	\item $a_0<0$, $R_4>0$, $R_5<0$, 
	\item $a_0<0$, $R_4>0$, $R_5>0$,
\end{itemize}
where $R_4$ and $R_5$ denote the fourth and the fifth element of the first column of the Routh matrix, respectively.

\begin{figure}[ht!]
	\centering
    	\subfigure[]{\includegraphics[scale=0.34]{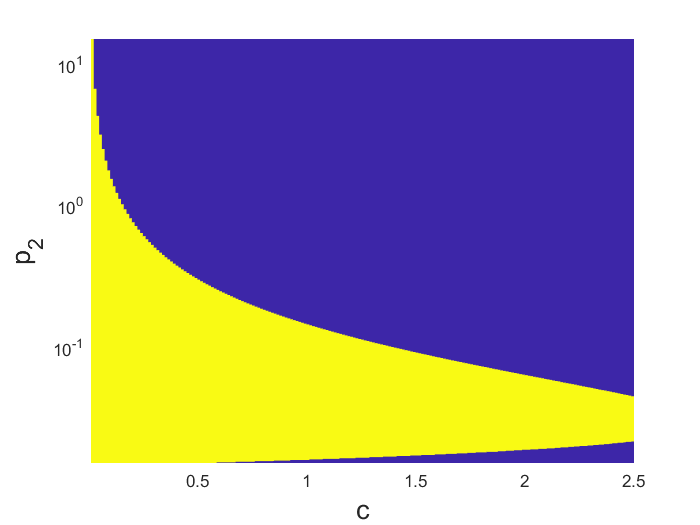}\label{cond1}}
    	\subfigure[]{\includegraphics[scale=0.34]{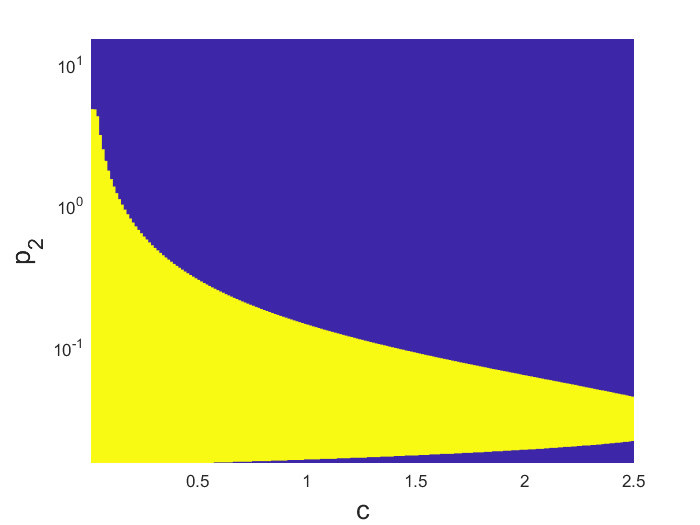}\label{cond2}}
    	\caption{Region of the  $p_2\,c$-plane where the Routh conditions are satisfied: (a) untreated case and (b) treated case. Note that a logarithmic scale is applied to the $p_2$-axis. The remaining parameters are listed in Table \ref{tab:par}.}
    		\label{cond}
\end{figure}

The region of the (p2 c)-plane where the Routh conditions are satisfied, both untreated and treated case are shown in Figure \ref{cond}, only in the first set of constrains. 

\begin{figure}[ht!]
	\centering
    	\subfigure[]{\includegraphics[scale=0.34]{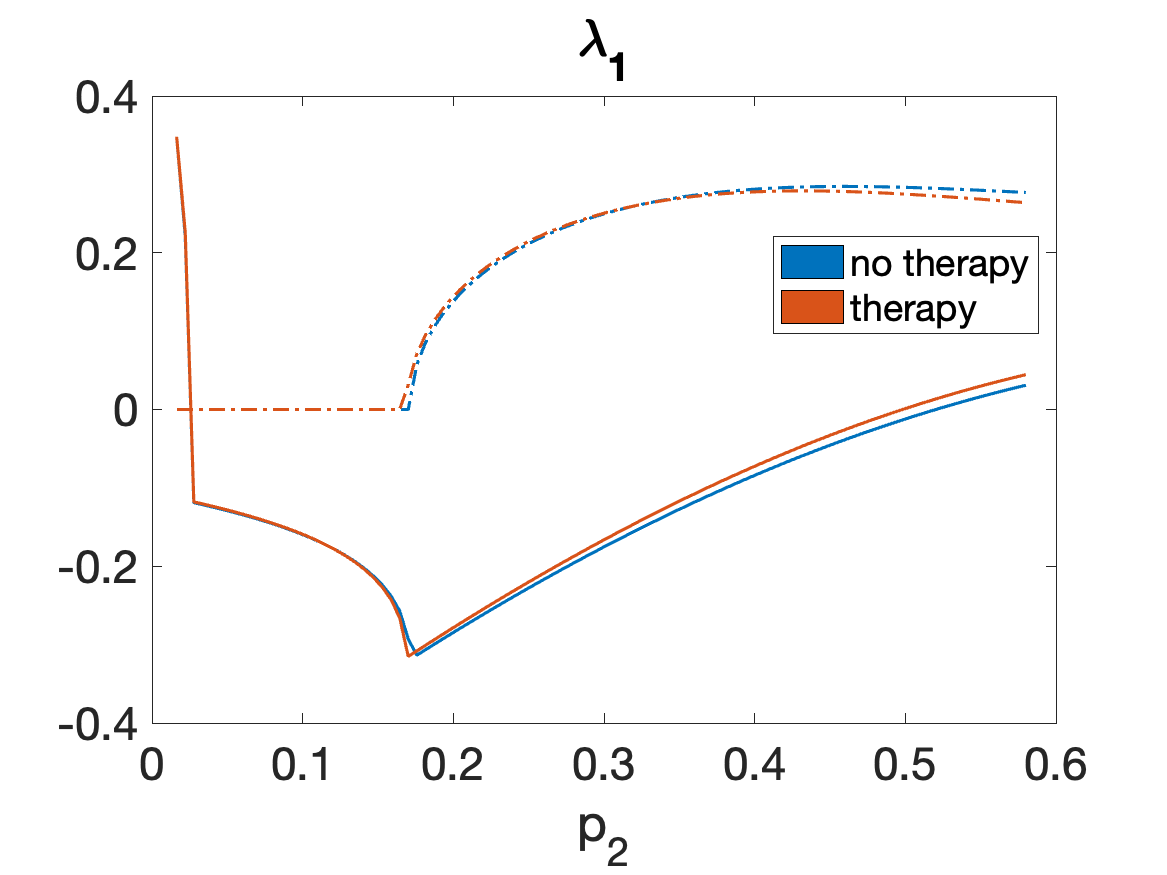}}
    	\subfigure[]{\includegraphics[scale=0.34]{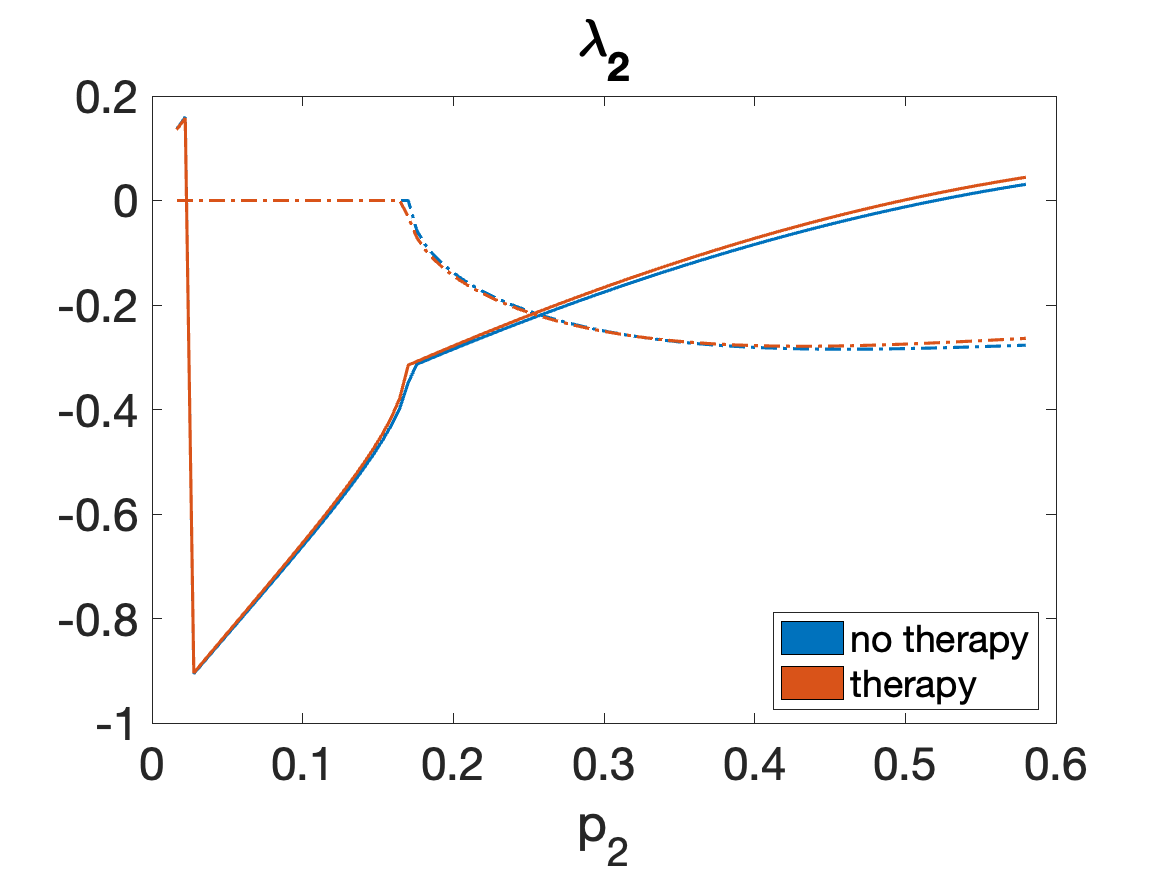}}
    	\subfigure[]{\includegraphics[scale=0.34]{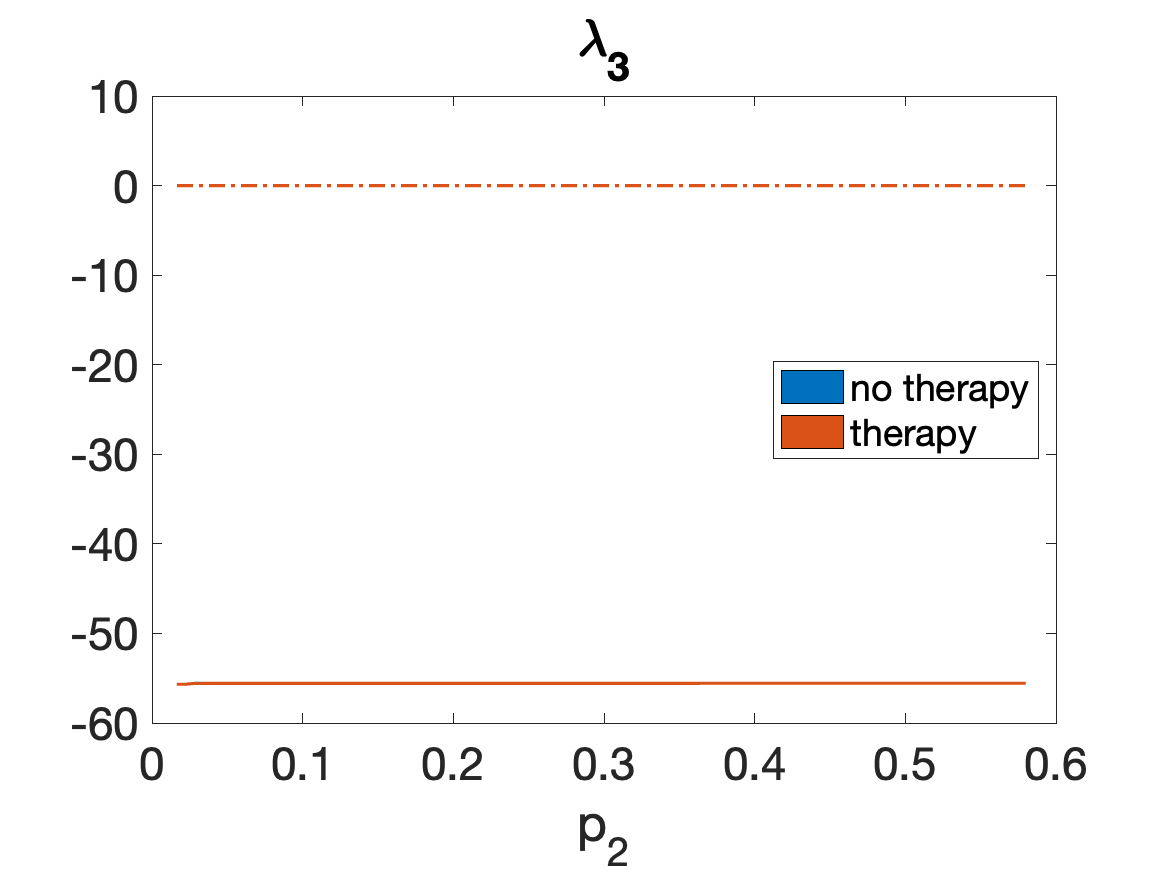}}
    	\caption{The solid and dotted lines represent the real part $\text{Re}(\lambda)$ and the imaginary part $\text{Im}(\lambda)$ of the eigenvalues $\lambda$, respectively, both with and without therapy. The parameters values are given  in Table \ref{tab:par}, with $c=0.25$.}
       \label{fig:eigenvalues}
\end{figure}

Moreover, we fix $c=0.25$ and study the stability as $p_2$ varies within the interval $[0.017,0.58]$, where the existence of the CCE is guaranteed. In this framework, the eigenvalues $\lambda$ of the Jacobian matrix, (both with and without therapy)
evaluated at $\mathbf{U}_2^*$, for $p_2 \in [0.017,0.58]$, are shown in Figure \ref{fig:eigenvalues}. Hence, we can conclude that there exists a subinterval $(\alpha, \beta) \subset [0.017,0.58]$ in which the CCE is stable for all $p_2 \in (\alpha, \beta)$, where $\alpha$ and $\beta$ are appropriate values. Their values may differ under two distinct conditions, depending on the presence or absence of therapy.

\subsubsection{Hopf bifurcation: periodic solutions}
\begin{figure}[ht!]
	\centering
    	\subfigure[]{\includegraphics[scale=0.24]{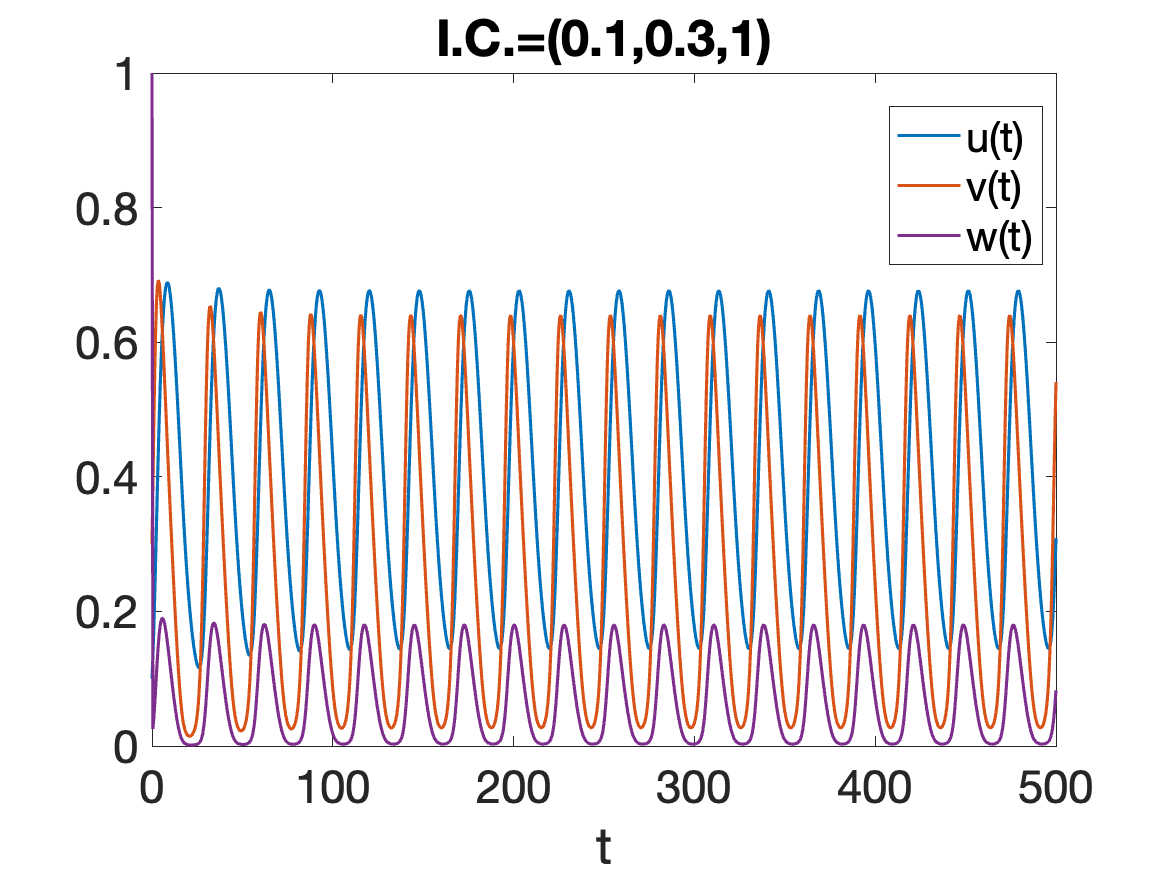}}
    	\subfigure[]{\includegraphics[scale=0.24]{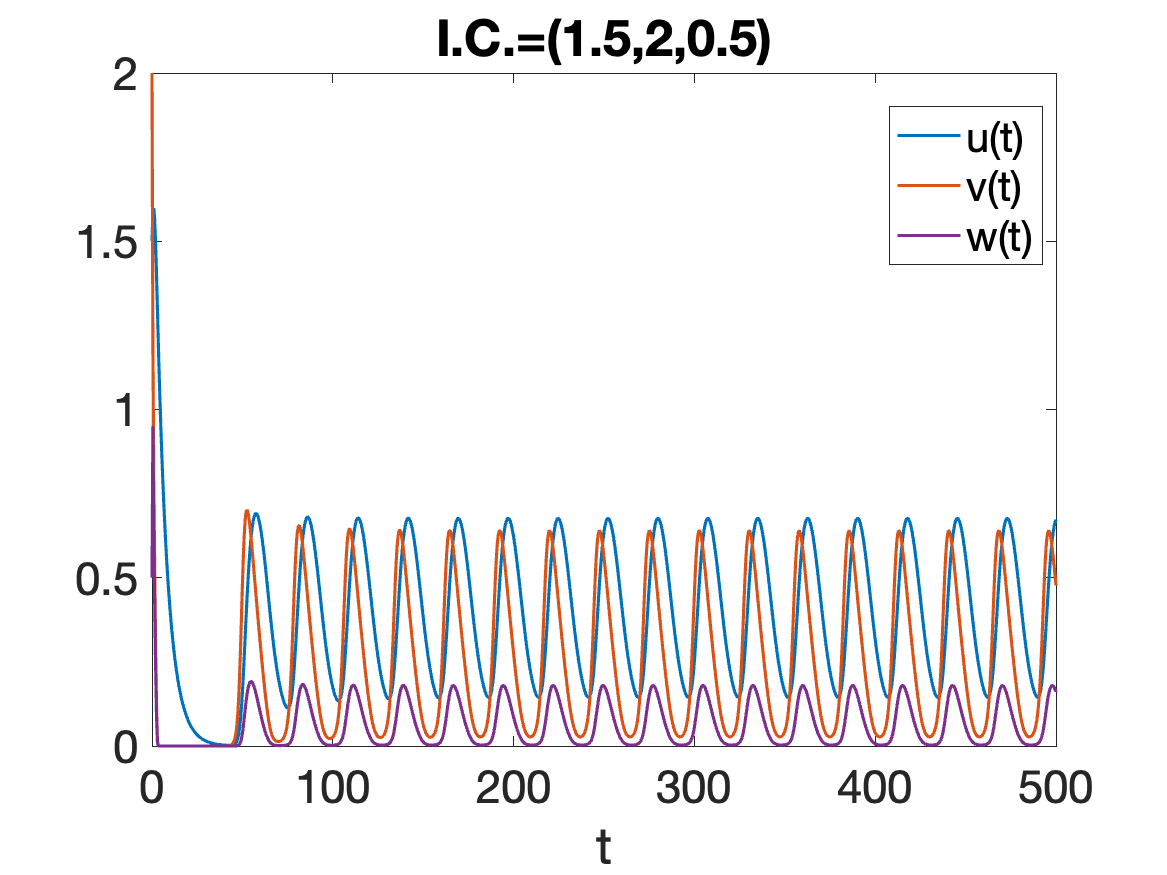}}
    	\subfigure[]{\includegraphics[scale=0.2]{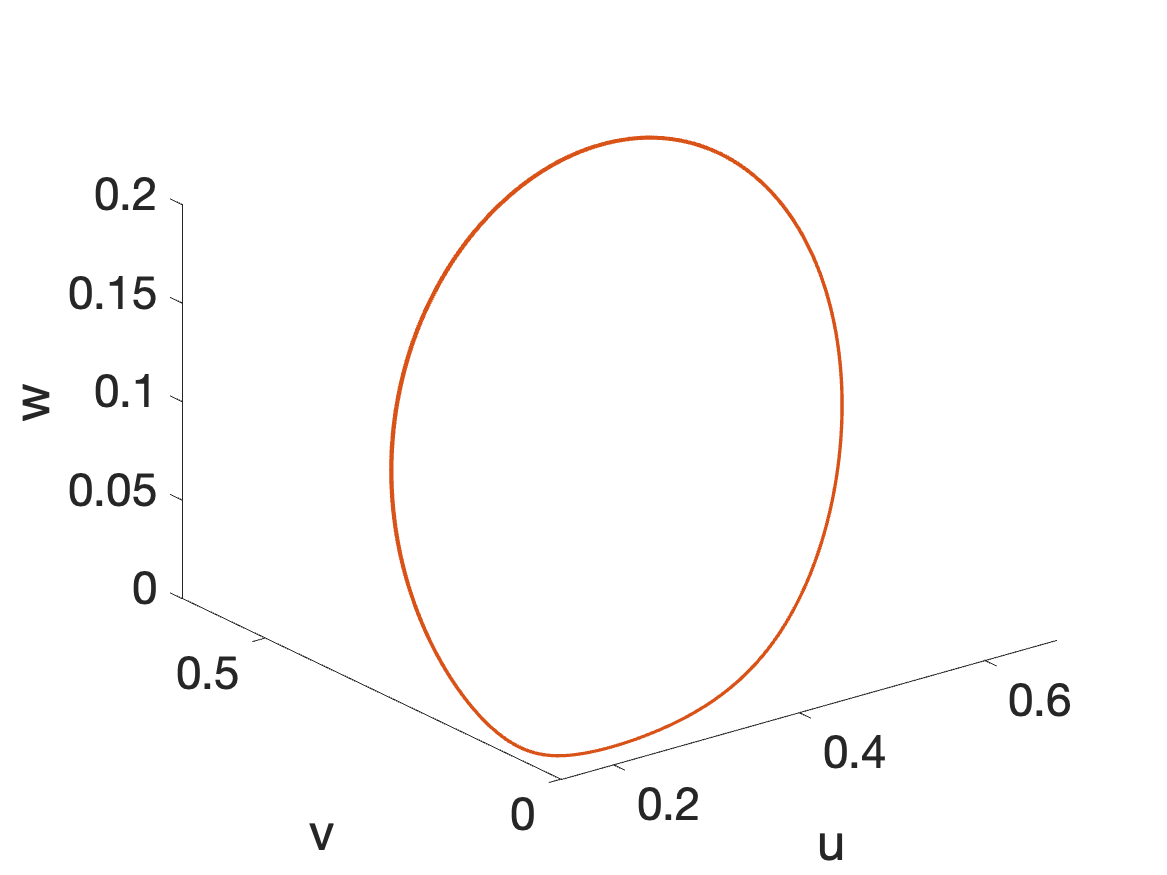}} \\
	\subfigure[]{\includegraphics[scale=0.24]{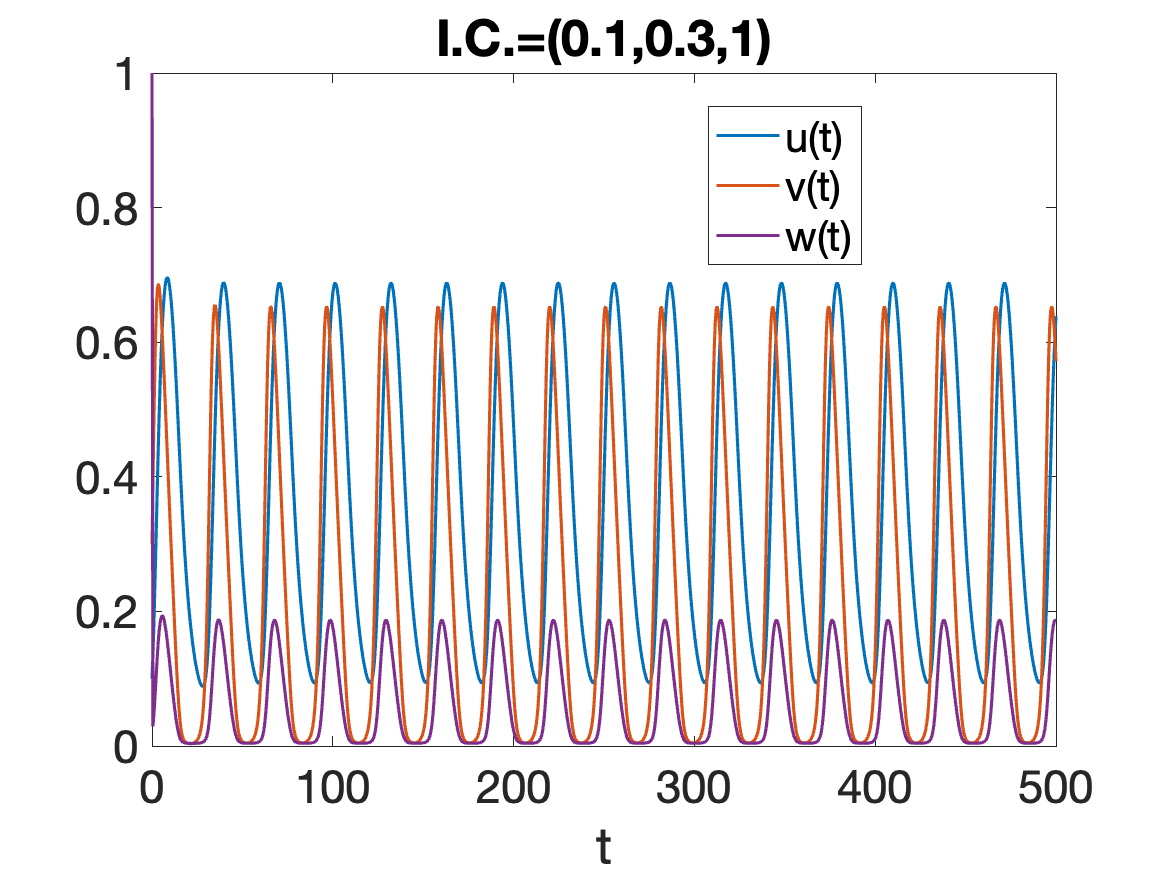}}
    	\subfigure[]{\includegraphics[scale=0.24]{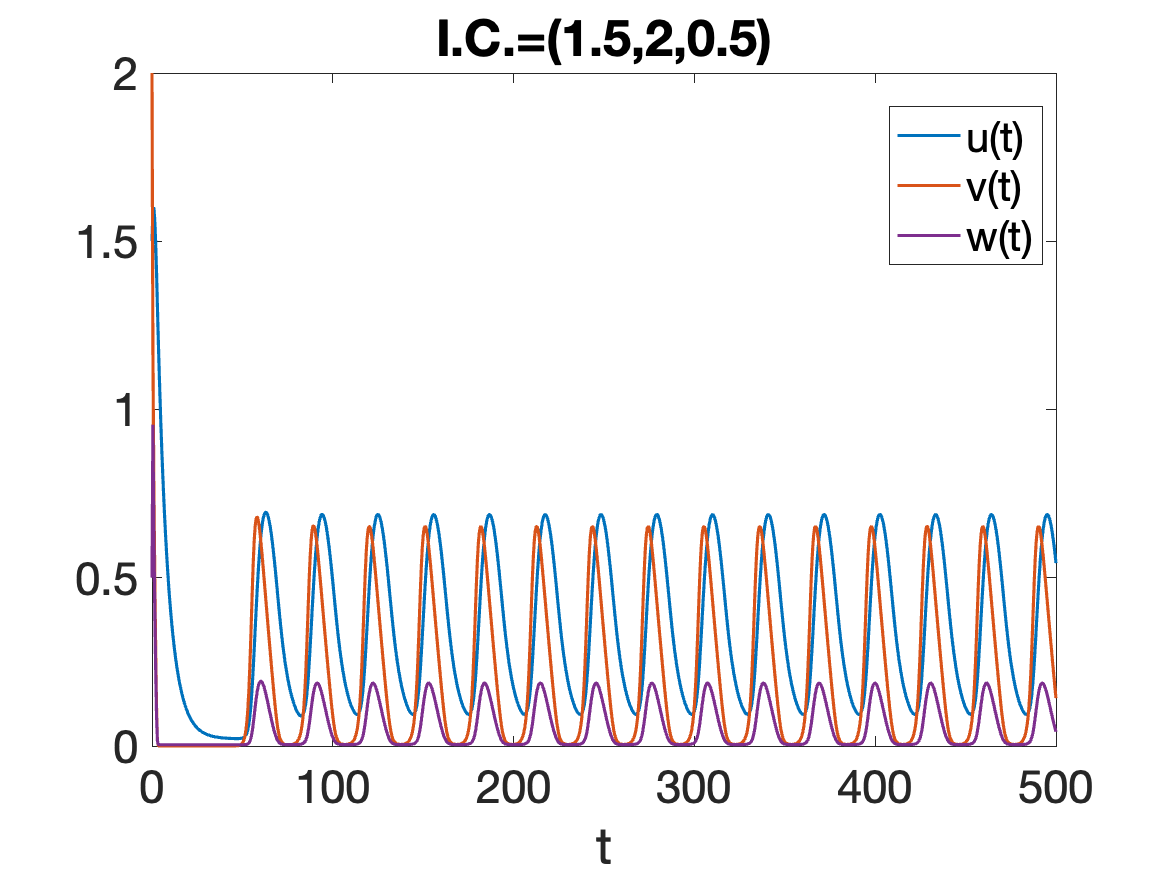}}
    	\subfigure[]{\includegraphics[scale=0.2]{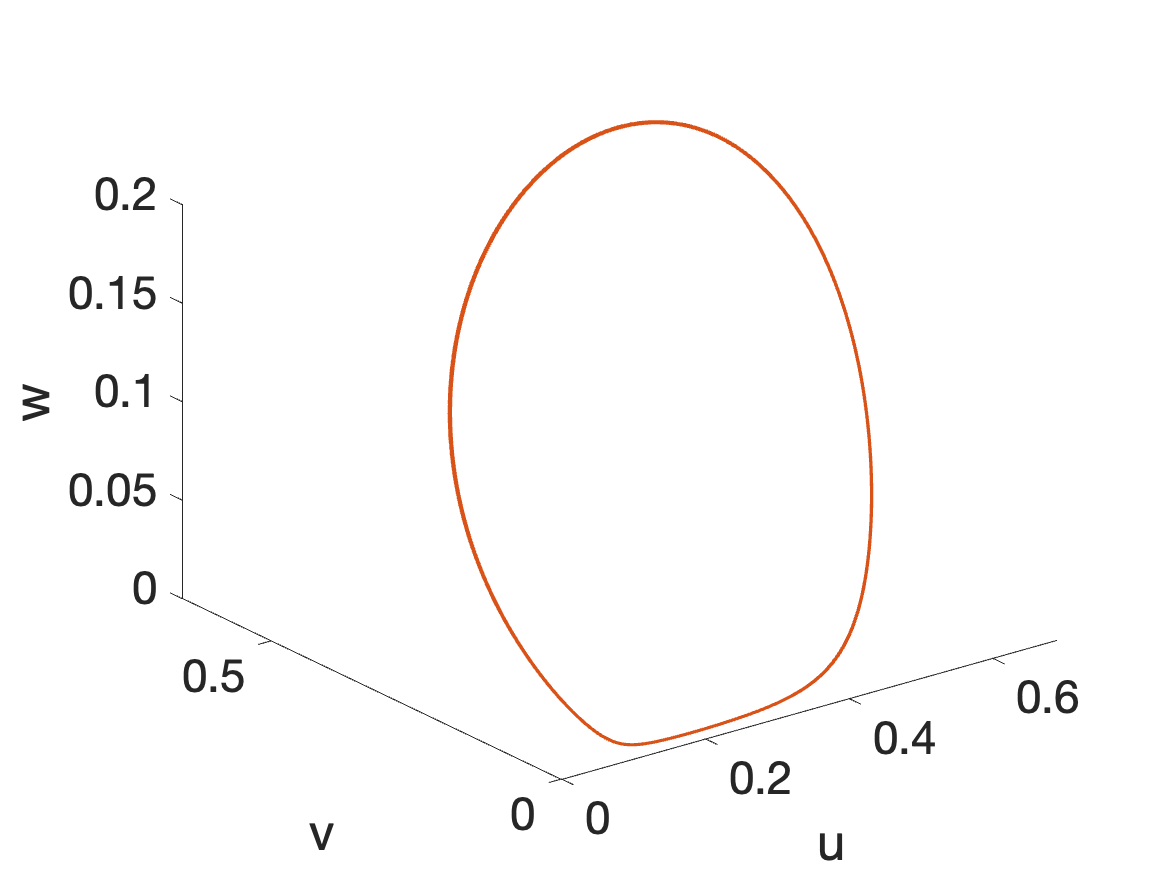}} 
    	\caption{First row: $s_1=s_3=0$ and Second row: $s_1=0.0035$, $s_3=0.2$. The remaining  parameters are listed in Table \ref{tab:par}, with $c=0.25$. The time evolutions for two different initial conditions are shown in the first and second column, \emph{i.e.} $u_0=0.1$, $v_0=0.3$, $w_0=1$ and $u_0=1.5$, $v_0=2$, $w_0=0.5$, respectively. The third column shows the limit cycle corresponding to the first set of initial conditions. } 
  \label{fig:hopf_time}
\end{figure}

It is possible to investigate when the 	Hopf bifurcation occurs. In particular, in the absence of therapy (\emph{i.e.}, $s_1= s_3=0$),  a stable limit cycle appears when the bifurcation parameter   $p_2$ reaches the critical value $p_2^c = 0.520$.  
However, even in the presence of therapy, such a bifurcation can still be observed, corresponding to the critical value $p_2^c = 0.498$. Some numerical simulations illustrating the presence of a stable limit cycle for $p_2 \geq p_2^c$ are presented in Figure \ref{fig:hopf_time}.

\section{Diffusion-driven instability}\label{sec_3}
The effect of self- and cross-diffusion terms in a reaction-diffusion system can result in  the loss of stability of an equilibrium point and lead to  the emergence of distint structures or patterns. This widely recognized phenomenon was first described in 1952 by Alan Turing in his pioneering paper  \cite{turing1990chemical} and studied in various contexts by a large number of authors 
\cite{carfora2025turing,tulumello2014cross,giunta2021pattern,della2022mathematical,inferrera2024reaction,vanag2009cross,xie2012cross,hao2020spatial,zincenko2021turing,chakraborty2021diffusion,aymard2022pattern}

Let us begin the analysis by considering the general system \eqref{syst_comp}. We will briefly outline the necessary computational steps and then determine the numerical conditions on the cross-diffusion coefficients introduced in order to obtain  Turing instability for the CCE.

 This  is studied in  both scenarios, with and without therapy, by setting $c=0.25$, $p_2=0.5$, $d_{11}=0.001$, $d_{22}=1.99\cdot 10^{-5}$, $d_{33}=0.01$, while keeping the other parameters as shown in Table \ref{tab:par}. Denoting by $\mathbf{U^*}=(u^*,v^*,w^*)$ an asymptotically stable equilibrium of system \eqref{syst_comp} without the spatial terms, and taking $\mathbf{U}_0=(u_0,v_0,w_0)$ as a constant vector, let us consider the following perturbation:
 
\begin{equation}\label{tur}
	\mathbf{U}=  \mathbf{U}^* + \mathbf{U}_0 \exp{(\lambda t + i \,\mathbf{k}\cdot \mathbf{x})}. 
\end{equation}
By substituting relation \eqref{tur} and linearizing system \eqref{syst_comp}, we get
\begin{equation*}
	[\lambda \,\mathcal{I} - (\nabla_{\mathbf{U}}\mathbf{F}(\mathbf{U})\vert _{\mathbf{U} = \mathbf{U^*}} - \mathcal{D} \vert \mathbf{k}\vert^2) ]\,\mathbf{U}_0  = \mathbf{0}.
\end{equation*}
It is a homogeneous linear system that admits non-zero solutions for $\mathbf{U}_0$ provided that the following condition is satisfied 
\begin{equation*}
	\det (\lambda \,\mathcal{I} - (\mathcal{J}\vert _{\mathbf{U} = \mathbf{U^*}} - \mathcal{D} \vert \mathbf{k}\vert^2))=0,
\end{equation*}
where $\mathcal{J}\vert _{\mathbf{U} = \mathbf{U^*}} = \nabla_{\mathbf{U}} \mathbf{F}(\mathbf{U})\vert _{\mathbf{U}=\mathbf{U^*}}$. The equilibrium state loses its stability, if at least one of the eigenvalues of the following matrix
\begin{equation}\label{mat_tur}
	\mathcal{A} = \mathcal{J}\vert _{\mathbf{U} = \mathbf{U^*}} - \mathcal{D} \vert \mathbf{k}\vert^2
\end{equation}
has positive real part.

\subsubsection{Untreated patients}
As previously stated, in this case we have $s_1=s_3=0$. Let us consider the stable equilibrium $\mathbf{U}_2=( 0.592878, 0.372148, 0.295647)$, and write the characteristic polynomial of the matrix $\mathcal{A}$ given in \eqref{mat_tur},
\begin{equation*}
	- \lambda^3 + a_2(k) \lambda^2 + a_1(k) \lambda + a_0(k,d_{32})
\end{equation*}
wherein
\begin{align*}
	a_2(k) & = -55.59- 0.011 k^2,\\
     	a_1(k) &= -1.427 - 0.057 k^2 , \\
     	a_0(k,d_{32}) &= -4.471 + 0.006 k^2  + 0.008 (-0.118 + d_{32} k^2).
\end{align*}

It is easily recognized that the coefficients $a_2(k)$ and $a_1(k)$ are negative. Therefore,  Turing instability may occur if 
\begin{equation} 
	a_0(k,d_{32})>0 \qquad \text{for some $k$}. \label{condition_1}
\end{equation}
One can observe that $a_0(k,d_{32})$ describes a parabola in $k$, so condition \eqref{condition_1} is satisfied if $d_{32}>-1.0668$.

\subsection{Treated patients}
As mentioned earlier, in the case of patients undergoing treatment, $s_1=0.0035$, $s_3=0.2$. Let us consider the stable equilibrium $\mathbf{U}_2=(0.586645,0.3542,0.296099)$, and write the characteristic polynomial of the matrix $\mathcal{A}$ given in \eqref{mat_tur}, 
\begin{equation*}
	- \lambda^3 + a_2(k) \lambda^2 + a_1(k) \lambda + a_0(k,d_{32})
\end{equation*}
where,
\begin{align*}
	a_2(k) & = -55.5898 - 0.0110199 k^2,\\
	a_1(k) &= -1.427 - 0.057 k^2, \\
	a_0(k,d_{32}) &= -4.194 + 0.007 k^2  + 0.008 (-0.129161 + d_{32} k^2).
\end{align*}
Similarly to the previous case, Turing instability can arise if \eqref{condition_1} is fulfilled. This occurs when $d_{32}>-1.45136 $.

\section{Numerical results}\label{sec_4}
In this Section, we present some numerical results. It is important to investigate the occurrence of Turing instability  arising in the neighborhood of a a equilibrium point in a reaction-diffusion system used to model cancer, both with and without treatment, where a stable and homogeneous state becomes unstable and leads to the formation of spatial structures. By analyzing the stability properties of the system in the neighborhood of the point where the transition from stable to unstable behavior occurs, researchers can derive useful information to predict tumor development under different conditions and, consequently, adopt therapeutic strategies to slow its progression.

Numerical integrations are performed by using an explicit finite difference method \cite{MitGri} with steps $\Delta x = \Delta y = 0.01$, and $\Delta t = 10^{-3}$. We also considered smaller steps and the accuracy of the results has not been affected.

The solutions clearly exhibit the formation of well recognized patterns. For both scenarios (without and with therapy) described by \eqref{syst_comp}, we consider initial perturbations that take the following gaussian forms
\begin{equation}\label{IC}
	\begin{aligned}
		& u(x,y,0)= \sum_{k=1}^4 b_j(k)\exp\left(-\dfrac{(x-x_0^j(k))^2+(y-y_0^j(k))^2}{\sigma_1^2}\right), \\
		& v(x,y,0)= \exp\left(-\dfrac{(x-x_1^j)^2+(y-y_1^j)^2}{\sigma_2^2}\right), &\qquad j=1,2,3 \\
		& w(x,y,0)= \sum_{k=1}^4 b_j(k)\exp\left(-\dfrac{(x-x_0^j(k))^2+(y-y_0^j(k))^2}{\sigma_1^2}\right), 
	\end{aligned}
\end{equation}
where
\begin{align*}
	& b_1=(1,1,1,1),\quad x_0^1=(0,0,1,1), \quad y_0^1=(0,1,0,1),\\
	& b_2=(1,0,0,0),\quad x_0^2=(1,0,0,0), \quad y_0^2=(1,0,0,0),\\
	& b_3=(1,0,0,0),\quad x_0^3=(0,0,0,0), \quad y_0^3=(0,0,0,0),
\end{align*}
and
\begin{align*}
	& x_1^1=0.5, \qquad y_1^1=0.5, \qquad x_1^2=1,\qquad y_1^2=1,\qquad  x_1^3=0, \qquad y_1^3=0, \\
	& \sigma_1^2=0.02, \qquad \sigma_2^2=0.06.
\end{align*}

\subsection{Turing pattern: Untreated patients}
\begin{figure}[ht!]
		\centering
    		\subfigure[]{\includegraphics[scale=0.3]{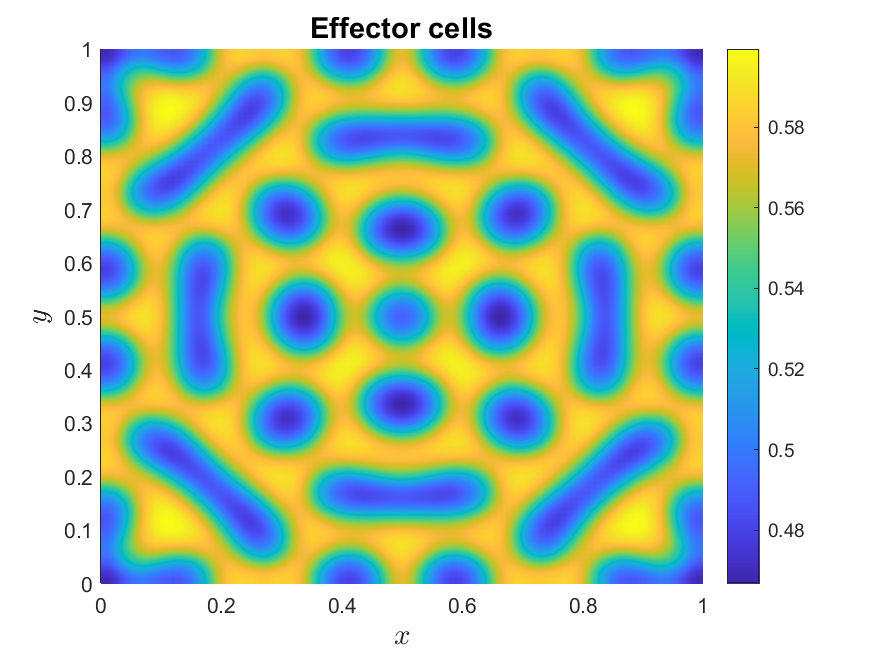}}
    		\subfigure[]{\includegraphics[scale=0.3]{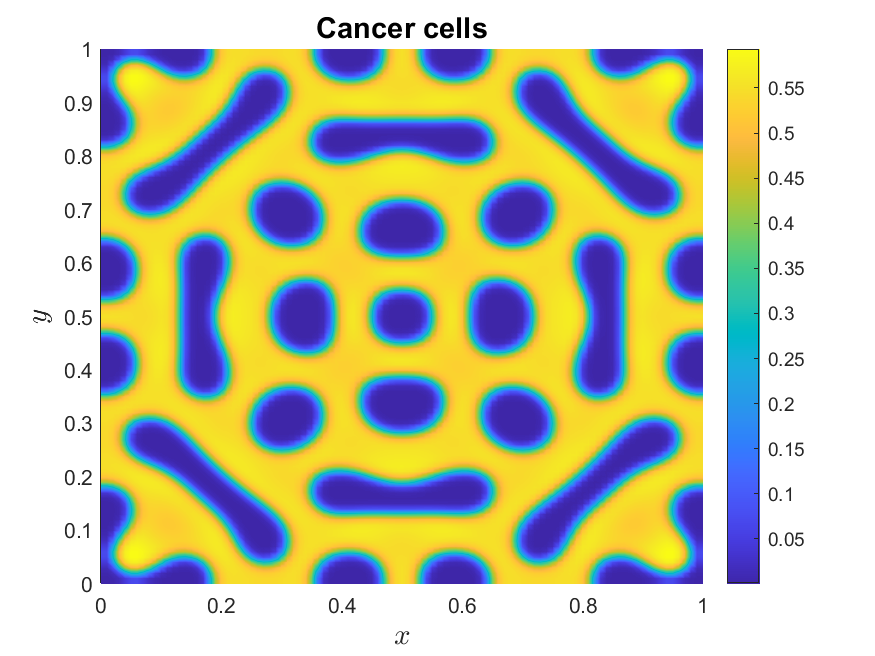}}
    		\subfigure[]{\includegraphics[scale=0.3]{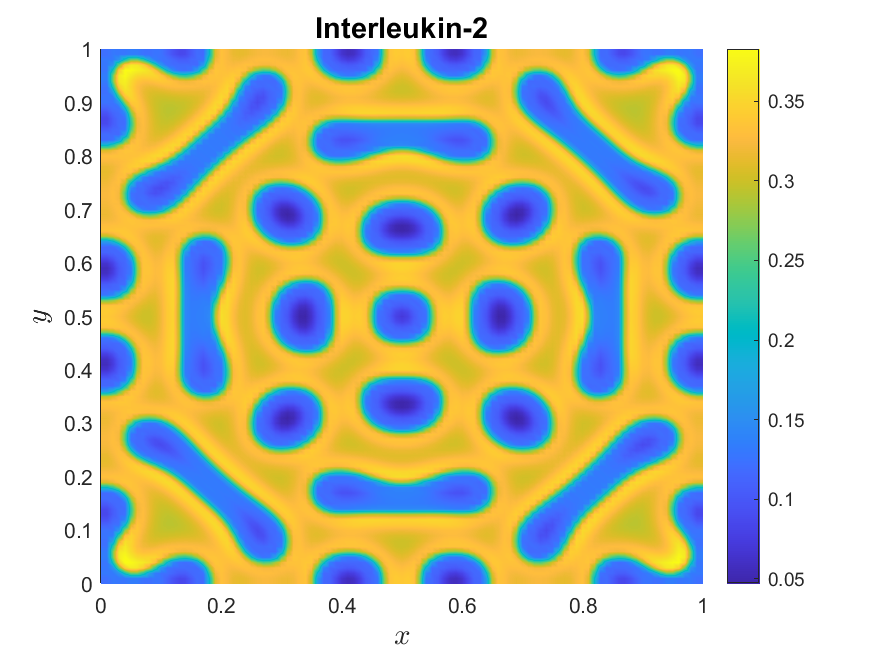}}\\
    		\subfigure[]{\includegraphics[scale=0.3]{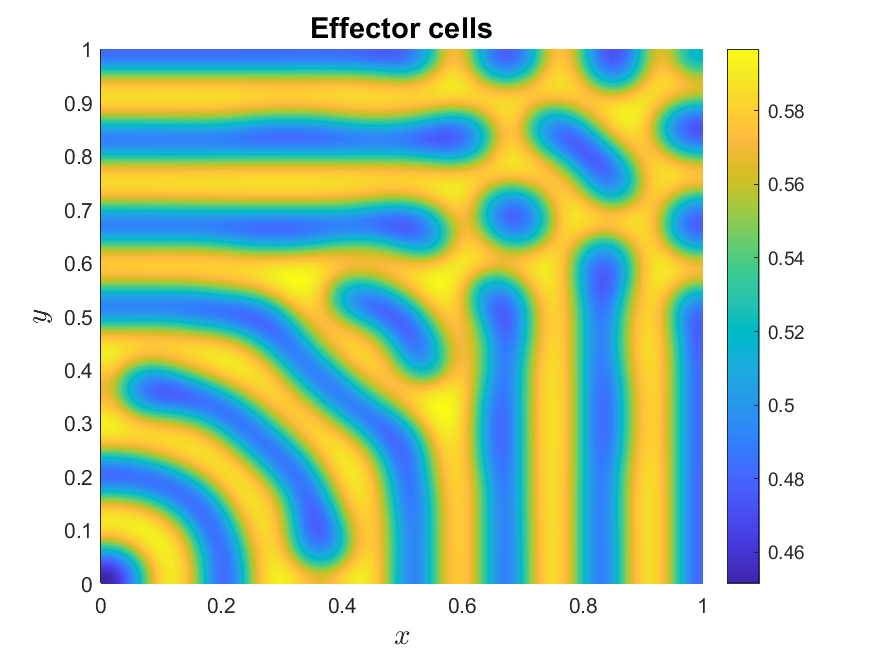}}
   		\subfigure[]{\includegraphics[scale=0.3]{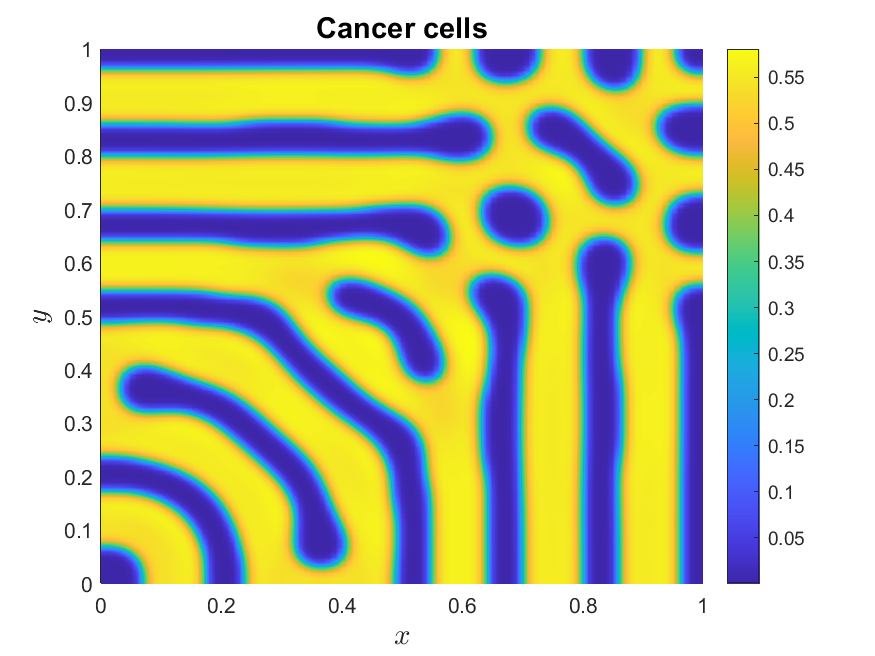}}
    		\subfigure[]{\includegraphics[scale=0.3]{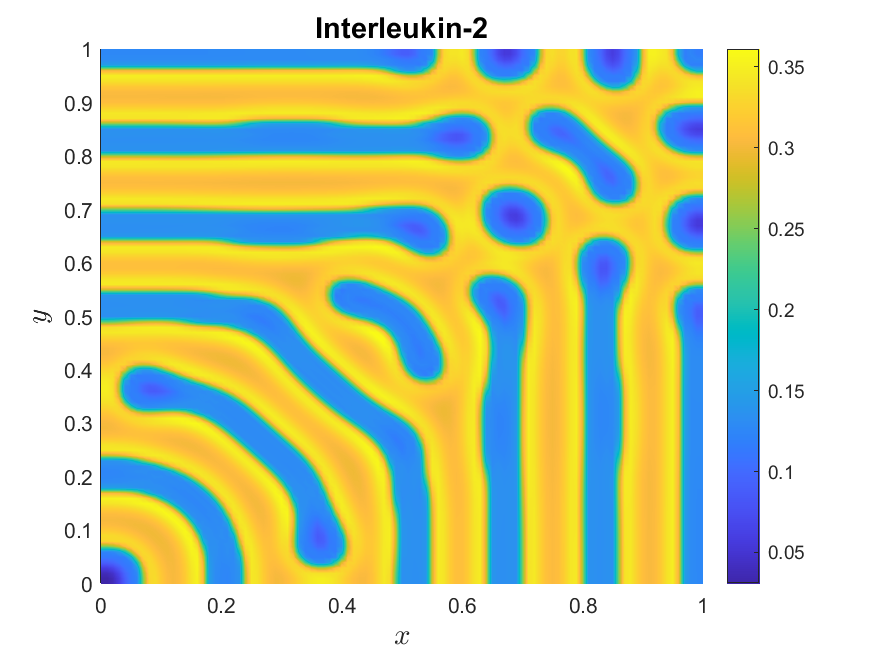}}\\
    		\subfigure[]{\includegraphics[scale=0.3]{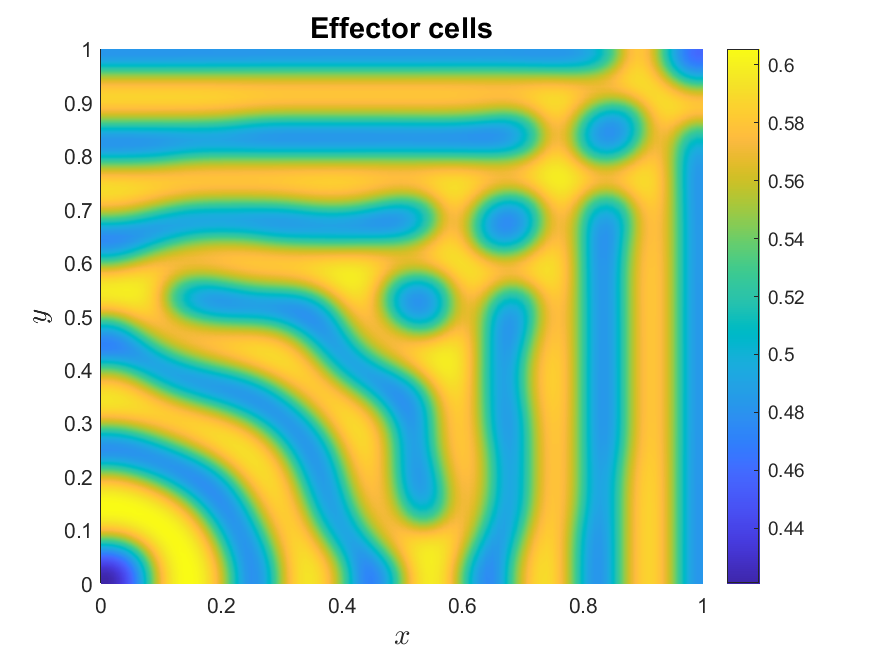}}
    		\subfigure[]{\includegraphics[scale=0.3]{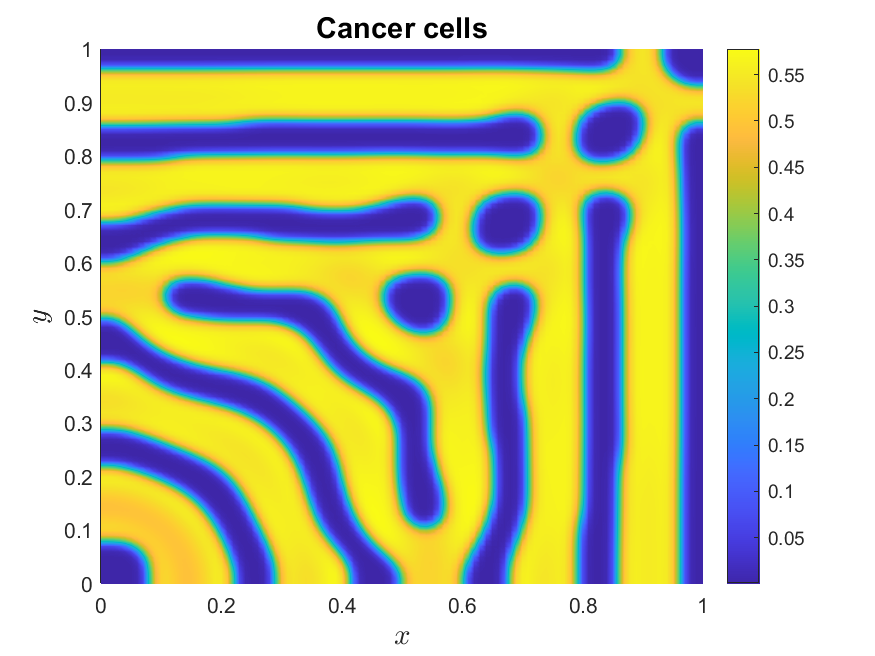}}
    		\subfigure[]{\includegraphics[scale=0.3]{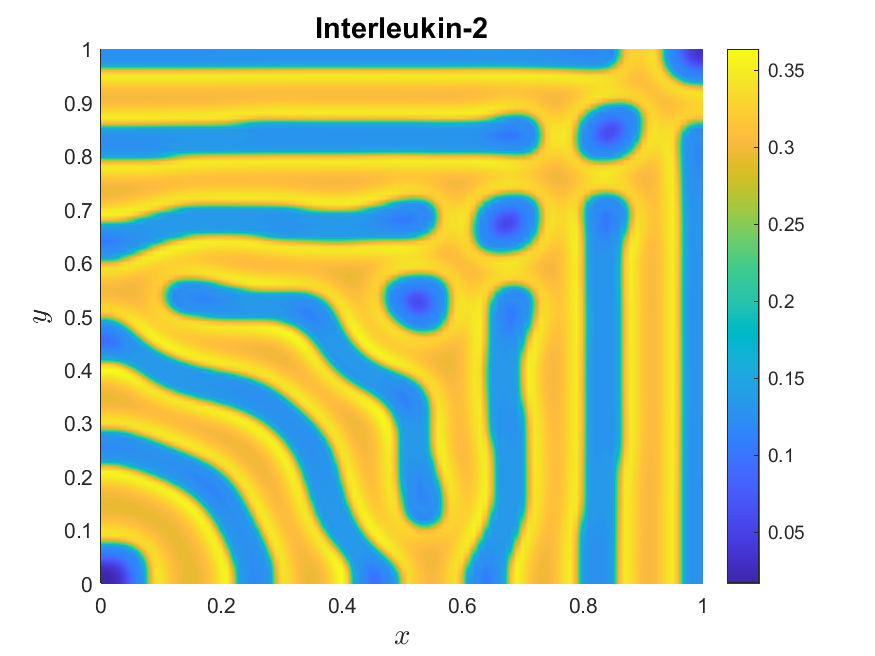}}
    		\caption{Spatial distributions of $u,\,v,\, w$ with $ d_{32}=-0.01$ in case of patient non-treated. The parameters are $c=0.25$, $p_2=0.5$, $d_{11}=0.001$, $d_{22}=1.99\cdot 10^{-5}$, $d_{33}=0.01$ and the remaining ones are provided in Table \ref{tab:par}}. 
		\label{sen_ter}
          \end{figure}

\begin{figure}[ht!]
		\centering
    		\subfigure[]{\includegraphics[scale=0.3]{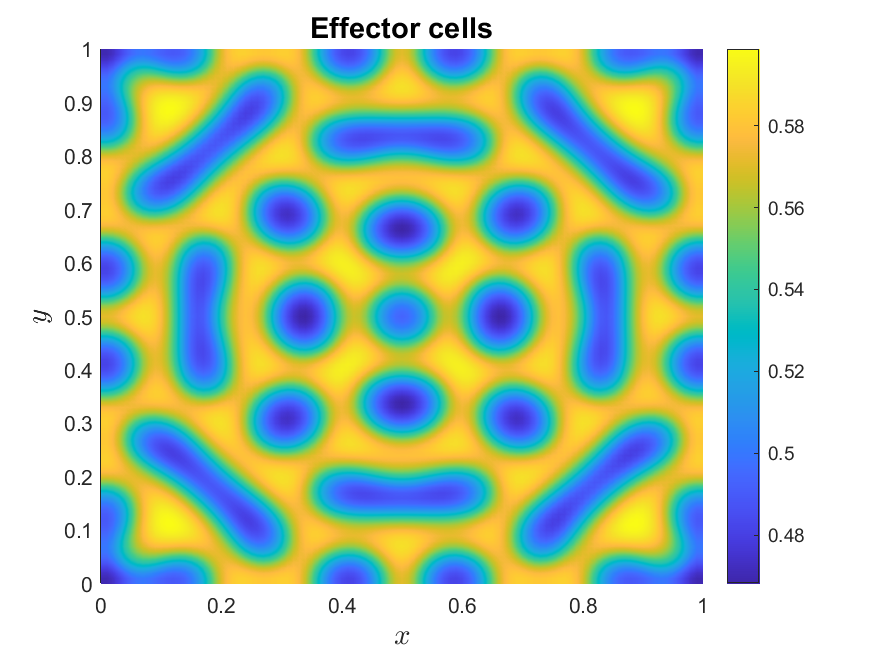}}
    		\subfigure[]{\includegraphics[scale=0.3]{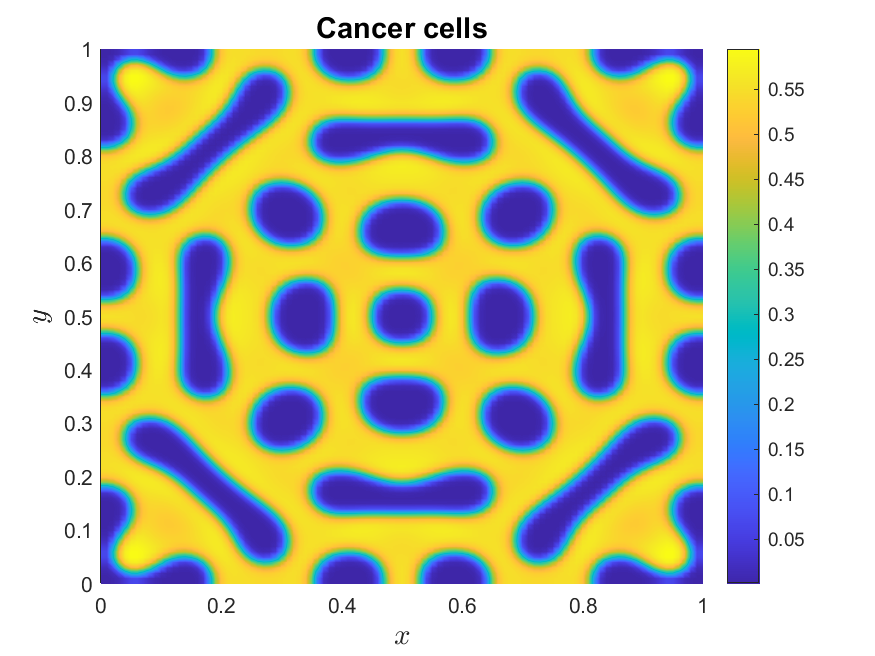}}
    		\subfigure[]{\includegraphics[scale=0.3]{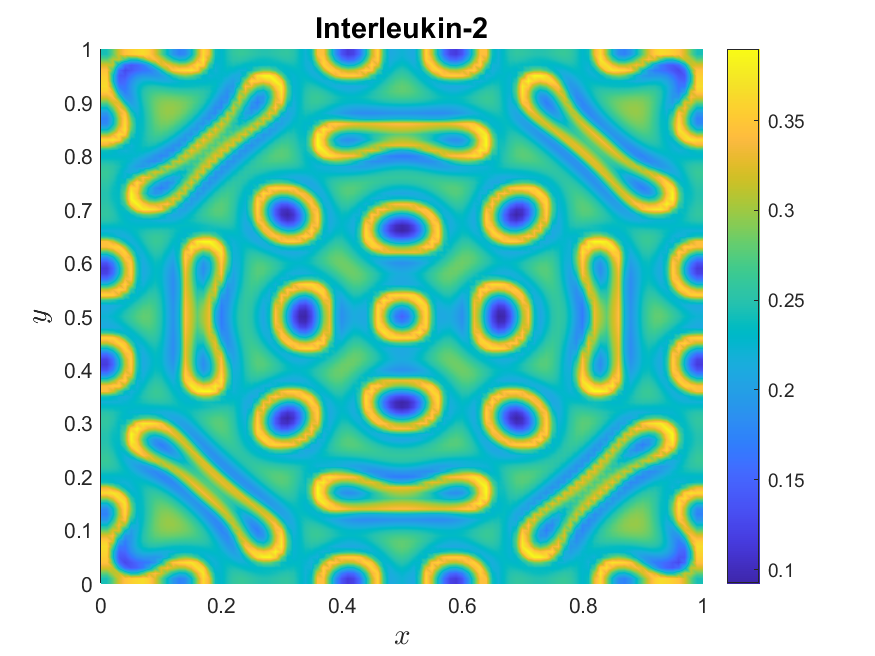}}\\
    		\subfigure[]{\includegraphics[scale=0.3]{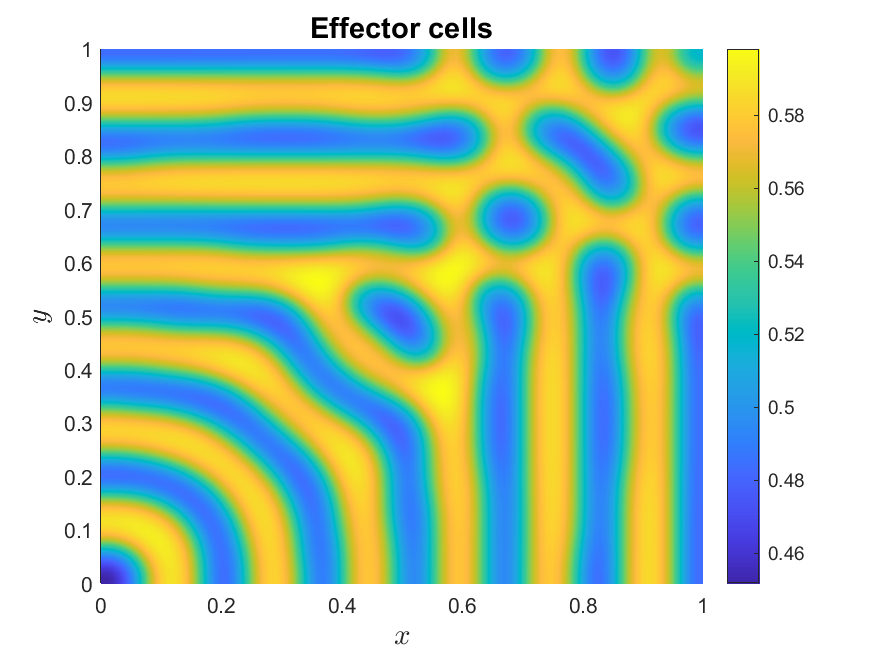}}
   		\subfigure[]{\includegraphics[scale=0.3]{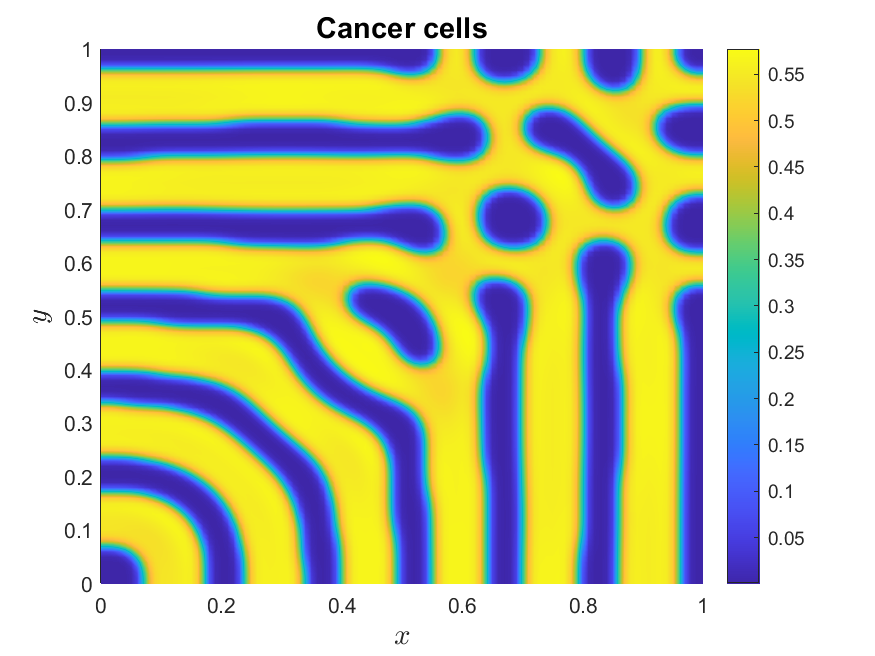}}
    		\subfigure[]{\includegraphics[scale=0.3]{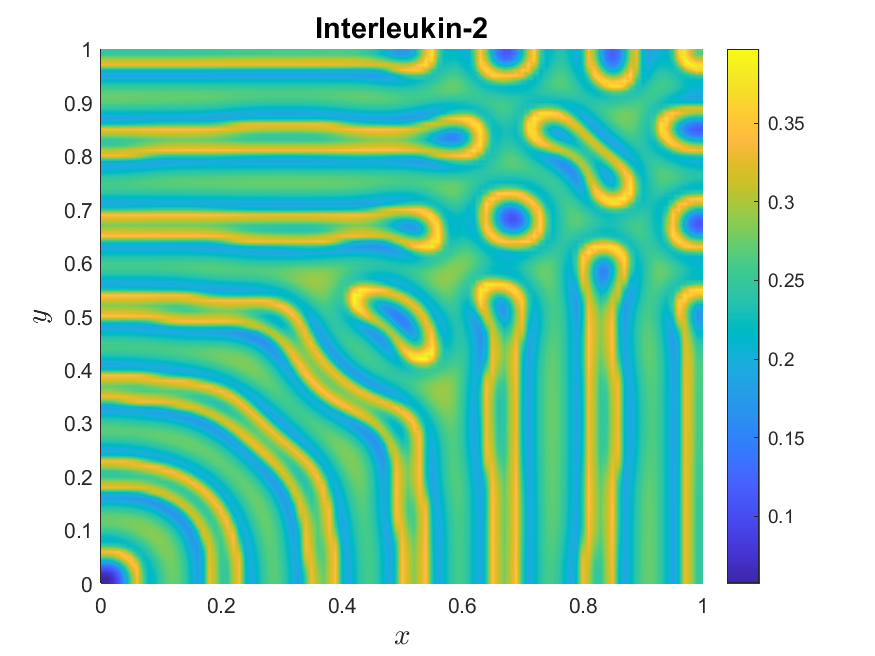}}\\
    		\subfigure[]{\includegraphics[scale=0.3]{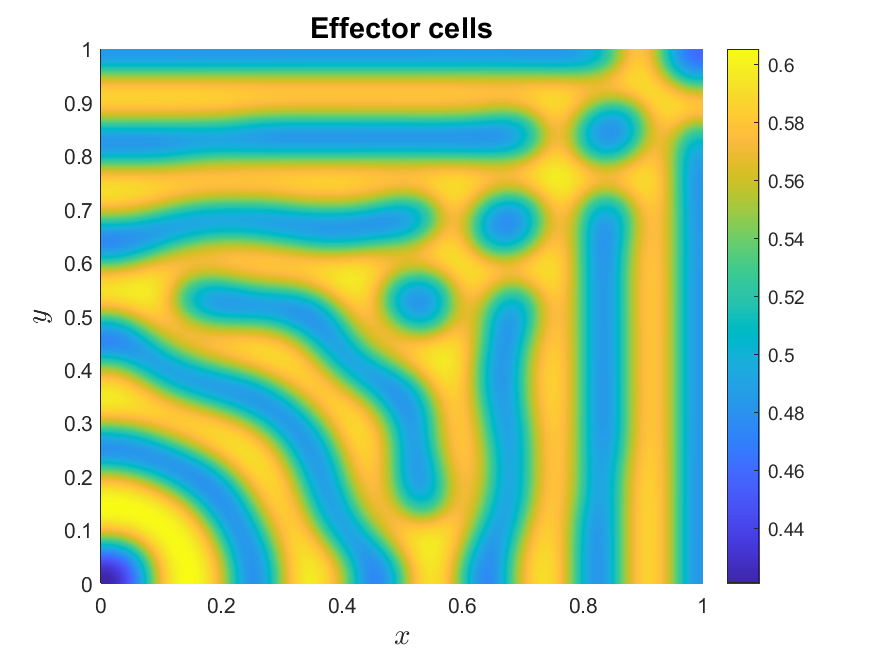}}
    		\subfigure[]{\includegraphics[scale=0.3]{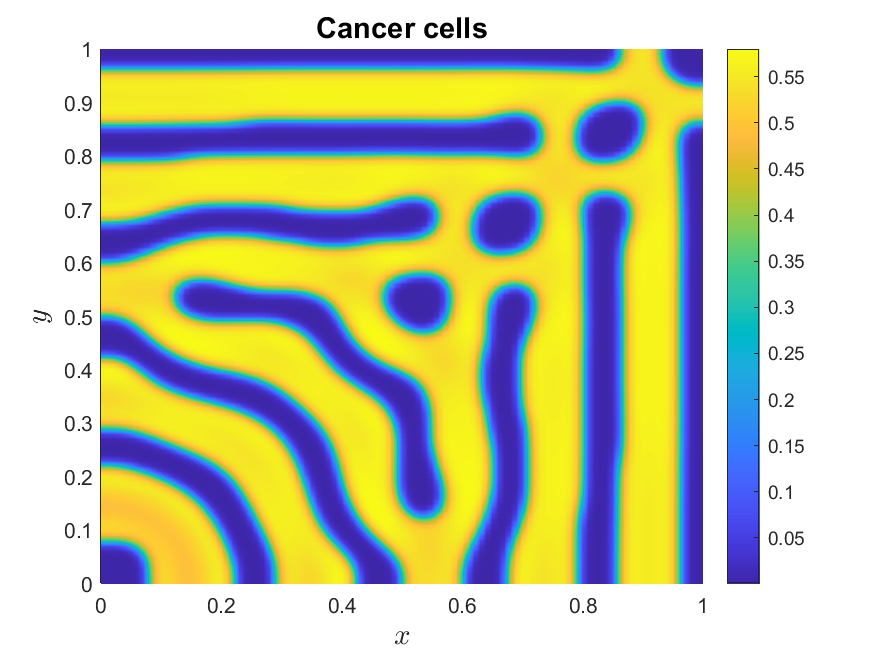}}
    		\subfigure[]{\includegraphics[scale=0.3]{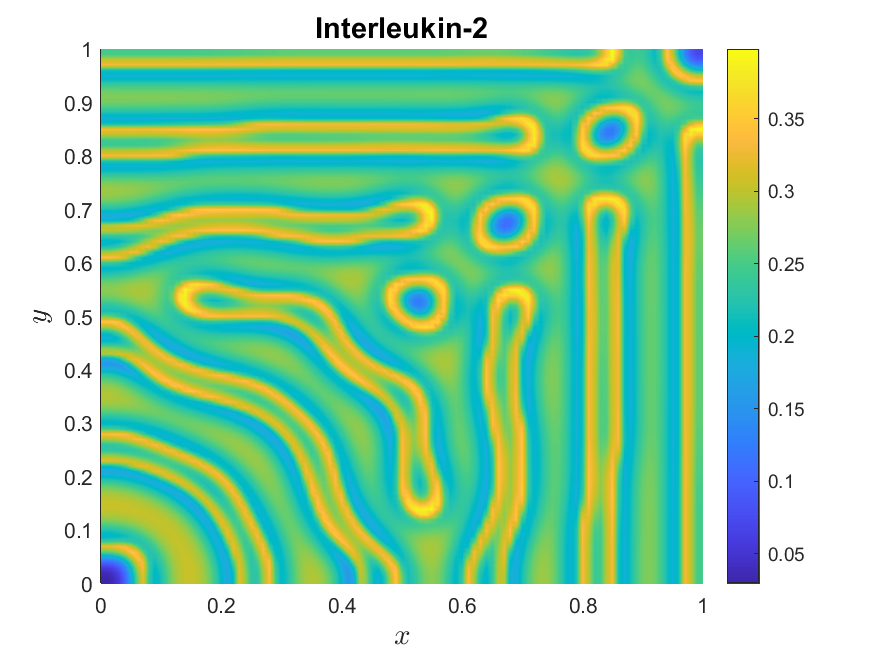}}
   		\caption{Spatial distributions of $u,\,v,\, w$ with $d_{32}=0.01$ in the case of an untreated patient. The parameters are $c=0.25$, $p_2=0.5$, $d_{11}=0.001$, $d_{22}=1.99\cdot 10^{-5}$, $d_{33}=0.01$; the remaining ones are listed in Table \ref{tab:par}}.
    		\label{sen_ter_1}
	\end{figure}

Figs. \ref{sen_ter} and \ref{sen_ter_1} show the contour plots of the solution representing the spatial distribution of effector cells, cancer cells and IL-2 over the 2D spatial domain obtained at $t=1000$. These patterns correspond to the fixed parameters and two different values of $d_{32}$ greater than the critical value, with three different initial conditions defined by equation \eqref{IC}, corresponding to $j=1,2,3$. From this point onward, in all figures, the color shading denotes cell concentrations, with yellow indicating high values while blue indicating low values. 

It can be observed that, as the cross-diffusion coefficient varies, $d_{32}=-0.01$ and $0.01$, a transition occurs in the stationary Turing patterns from a configuration consisting of cold spots mixed with stripe patterns to a structure predominantly covered by stripe patterns.

For higher values of the cross-diffusion coefficient, the cytokine subpopulation displays less pronounced stripe patterns, as a result of greater diffusion toward the tumor cells.

You can observe that, as time increases, the hot spots gradually merge until, at $t=1000$, the hot spots of the subpopulations have mostly coalesced.

\subsection{Turing pattern: Treated patients}  
\begin{figure}[ht!]
		\centering
    		\subfigure[]{\includegraphics[scale=0.3]{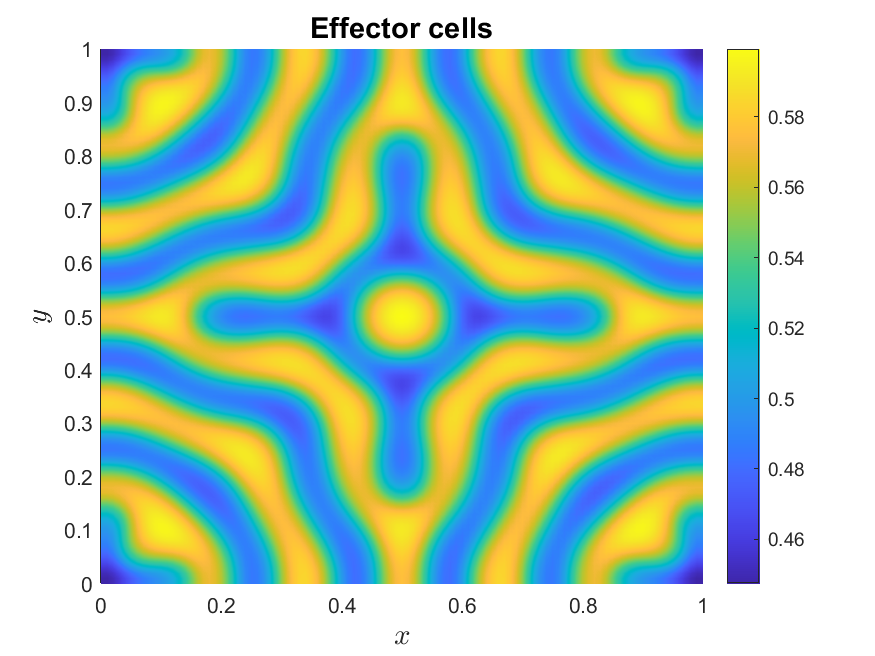}}
    		\subfigure[]{\includegraphics[scale=0.3]{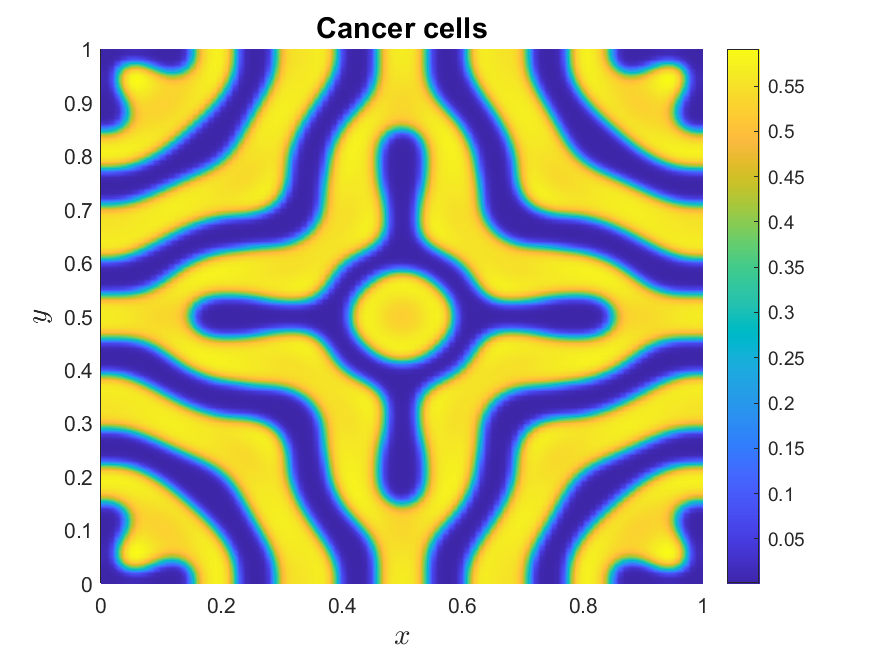}}
    		\subfigure[]{\includegraphics[scale=0.3]{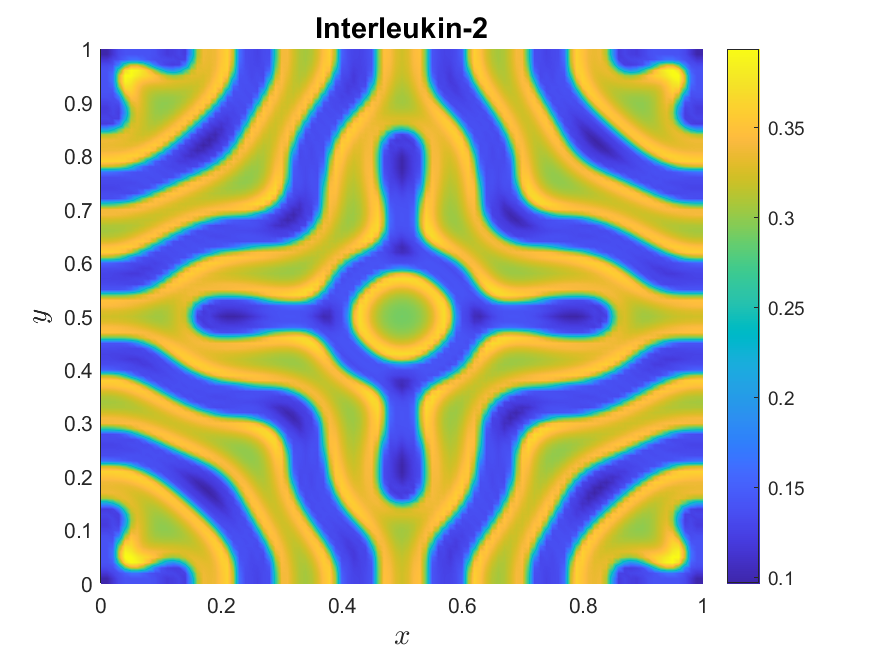}}\\
    		\subfigure[]{\includegraphics[scale=0.3]{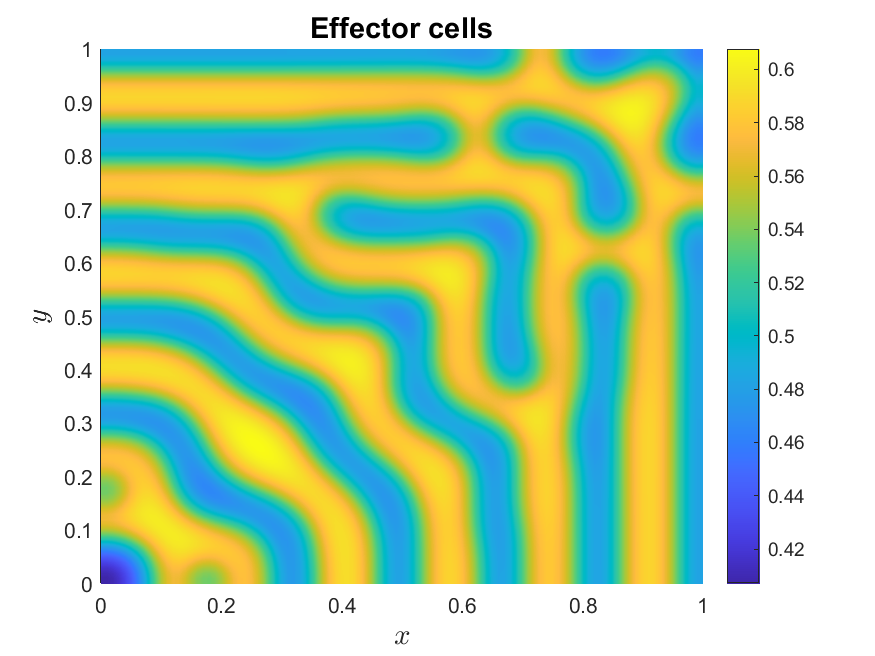}}
   		\subfigure[]{\includegraphics[scale=0.3]{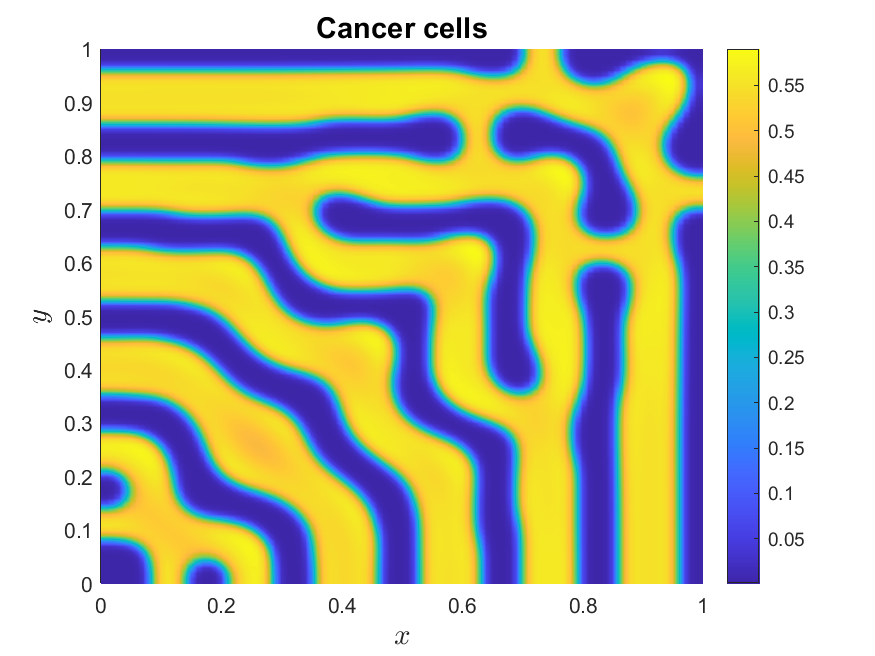}}
    		\subfigure[]{\includegraphics[scale=0.3]{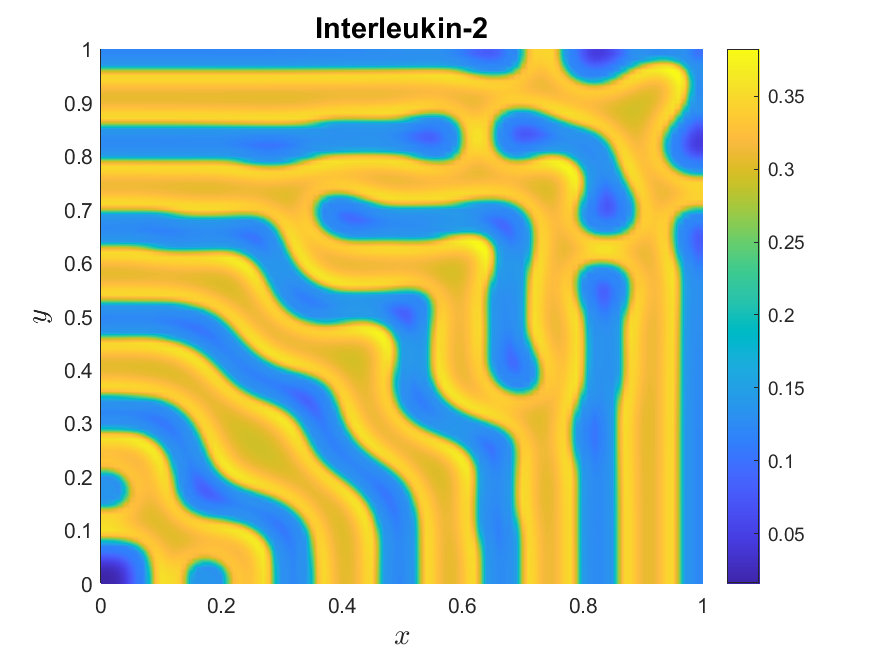}}\\
    		\subfigure[]{\includegraphics[scale=0.3]{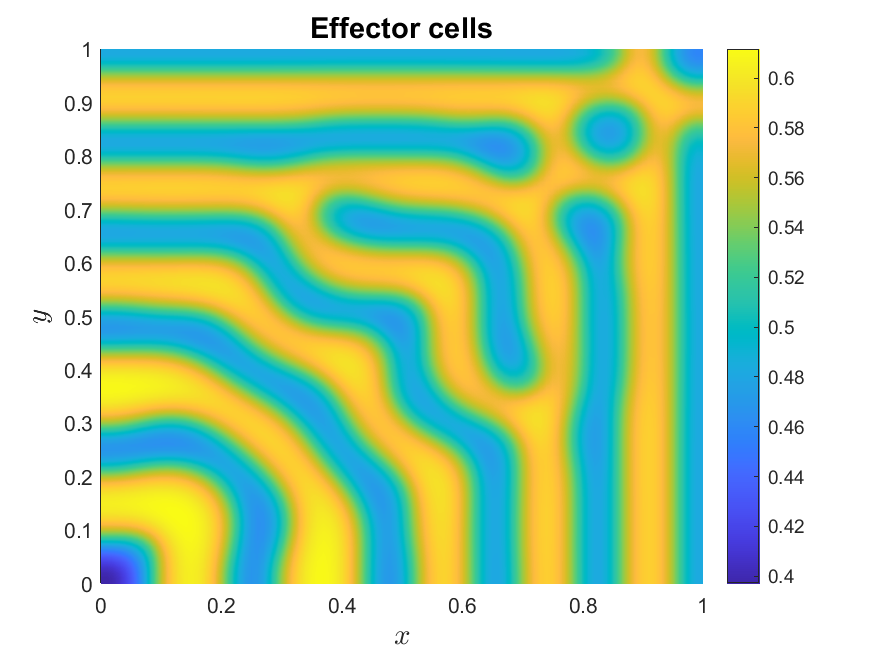}}
    		\subfigure[]{\includegraphics[scale=0.3]{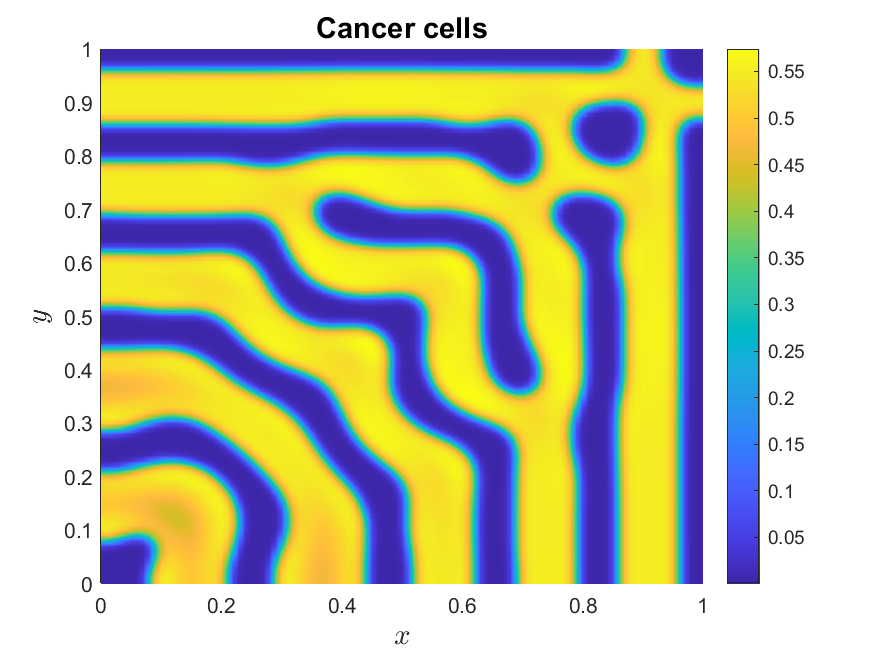}}
    		\subfigure[]{\includegraphics[scale=0.3]{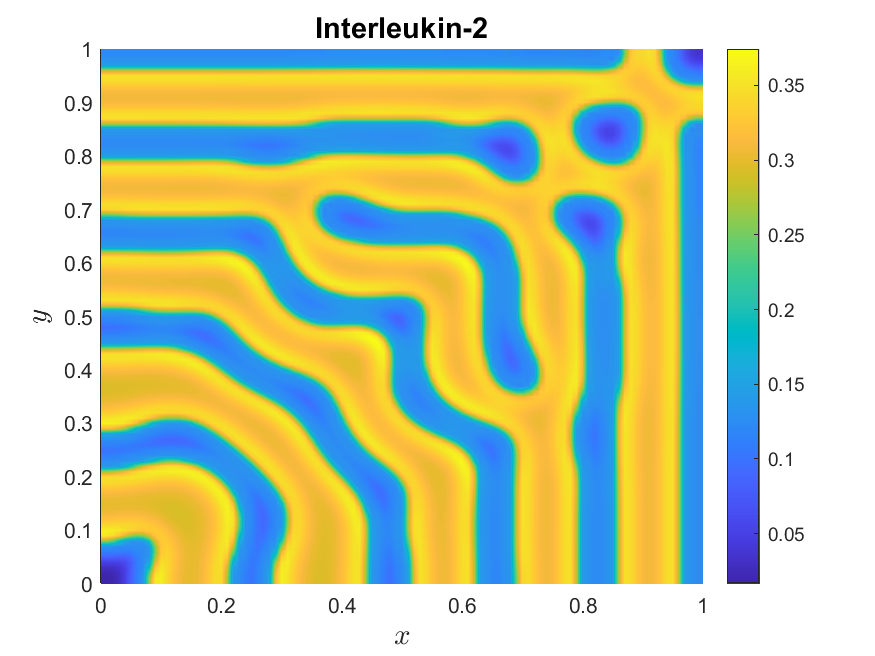}}
    		\caption{Spatial distributions of $u,\,v\,and w$ with $ d_{32}=-0.01$ in the case of a treated patient. The parameters are $c=0.25$, $p_2=0.5$, $d_{11}=0.001$, $d_{22}=1.99\cdot 10^{-5}$, $d_{33}=0.01$; the remaining ones are provided  in Table \ref{tab:par}}.
    		\label{con_ter}
	\end{figure}

\begin{figure}[ht!]
		\centering
    		\subfigure[]{\includegraphics[scale=0.3]{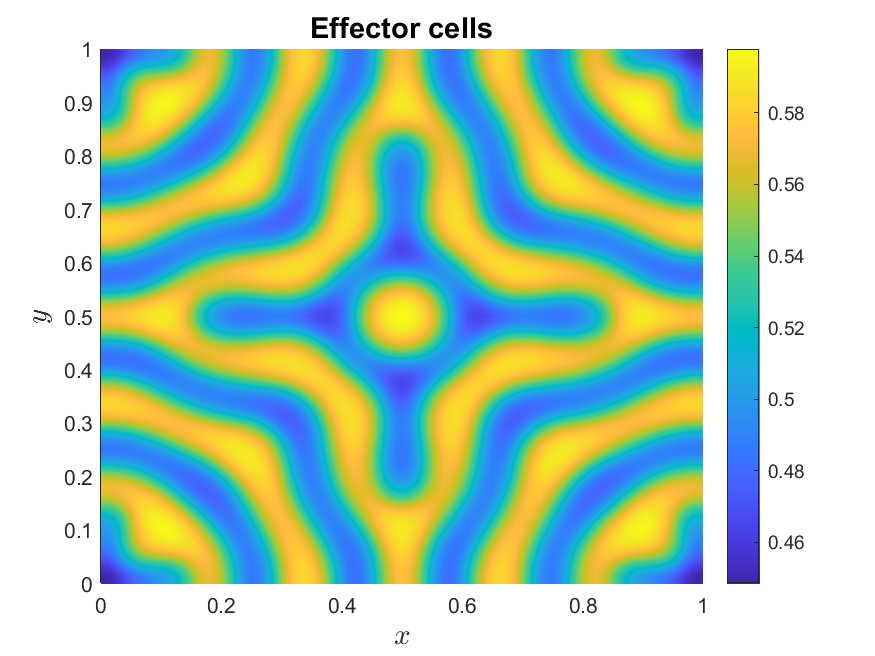}}
    		\subfigure[]{\includegraphics[scale=0.3]{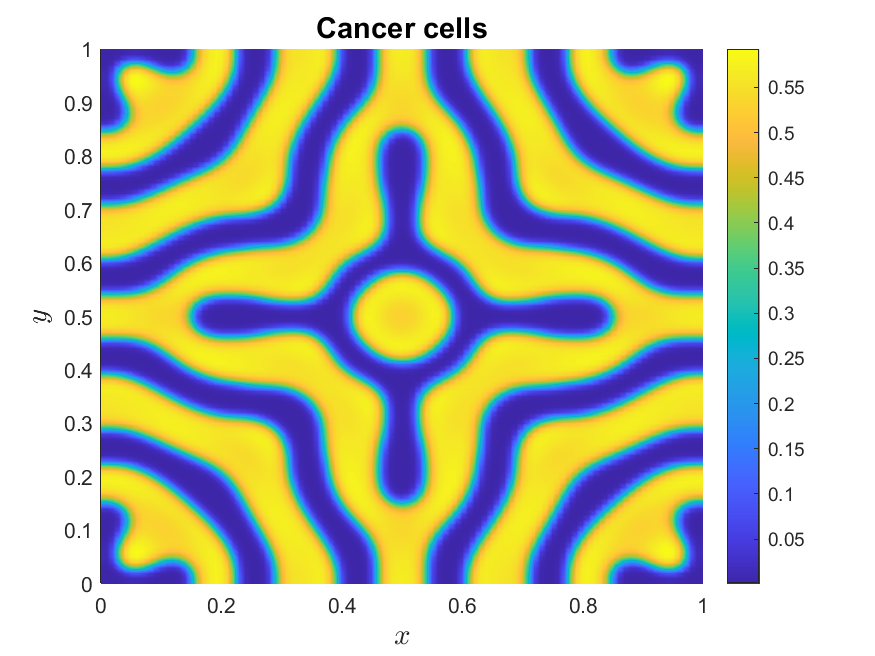}}
    		\subfigure[]{\includegraphics[scale=0.3]{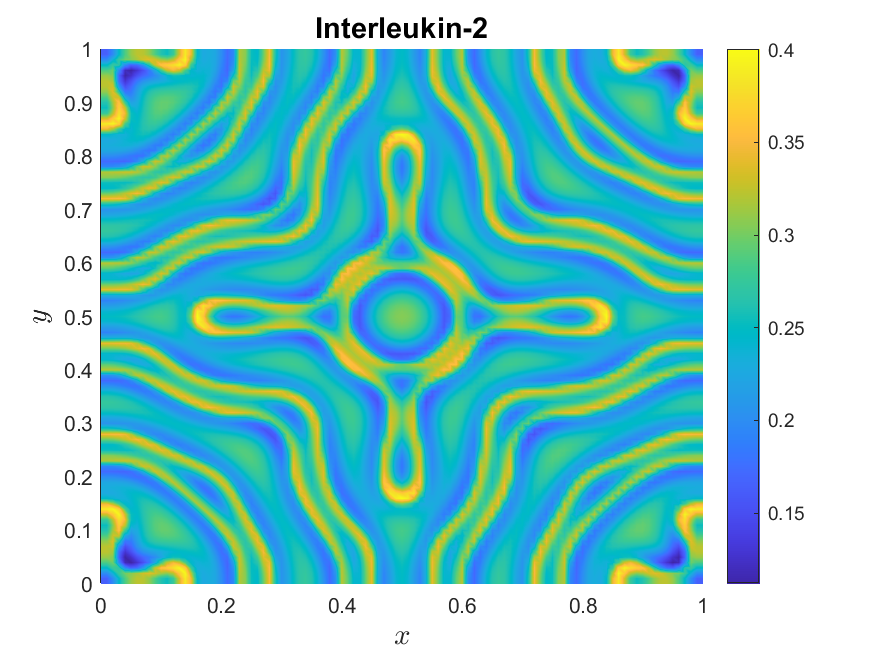}}\\
    		\subfigure[]{\includegraphics[scale=0.3]{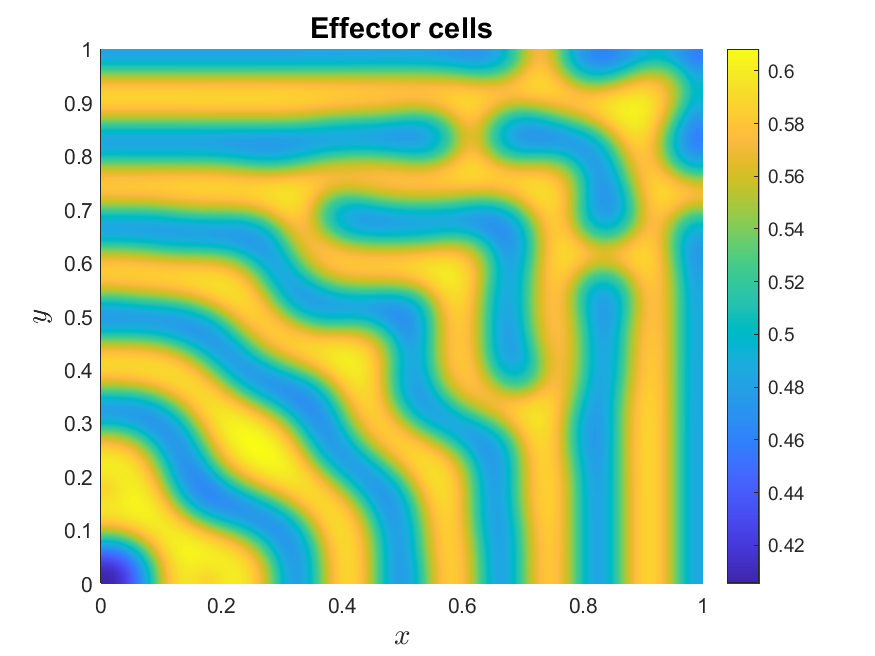}}
   		\subfigure[]{\includegraphics[scale=0.3]{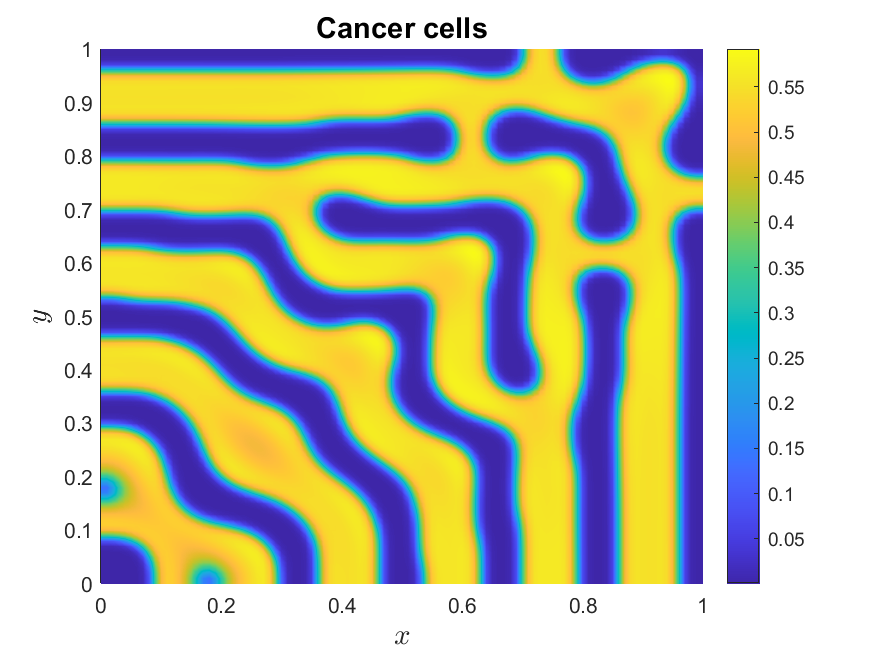}}
    		\subfigure[]{\includegraphics[scale=0.3]{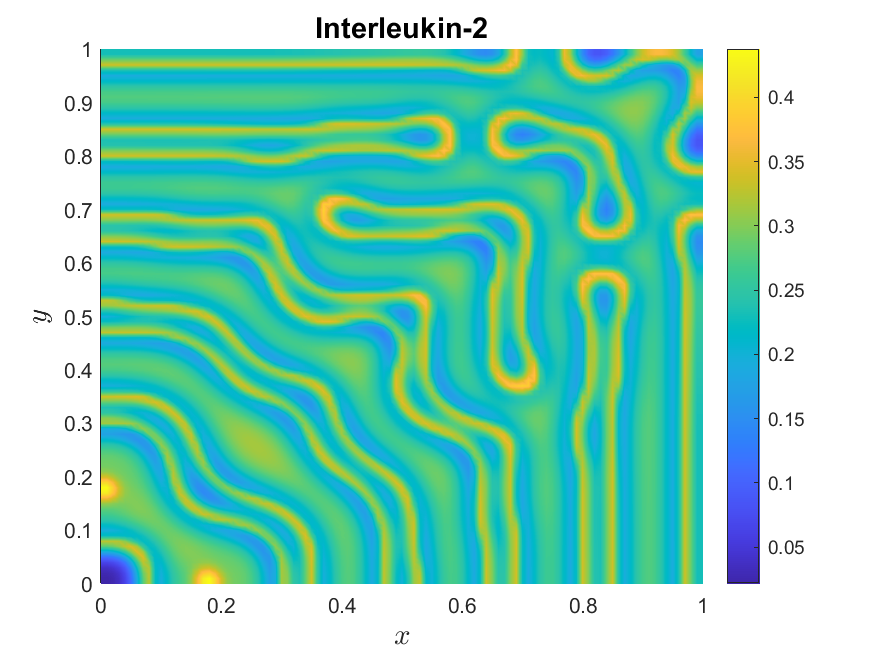}}\\
    		\subfigure[]{\includegraphics[scale=0.3]{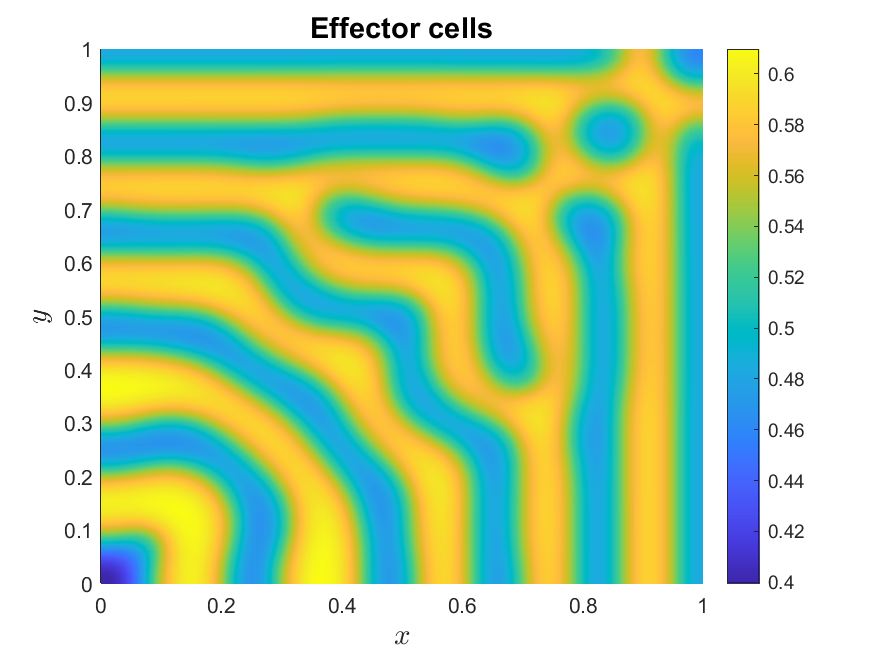}}
    		\subfigure[]{\includegraphics[scale=0.3]{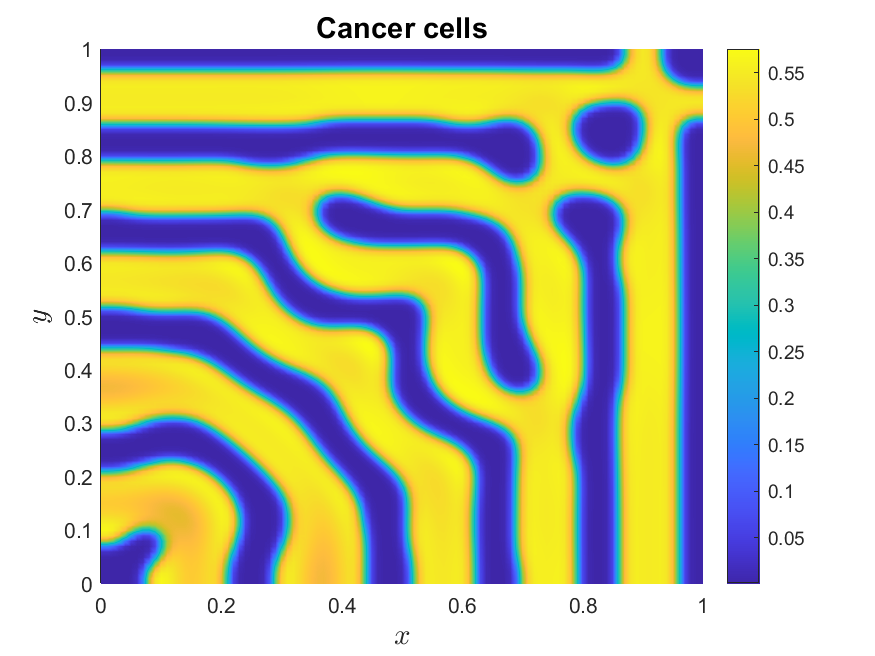}}
    		\subfigure[]{\includegraphics[scale=0.3]{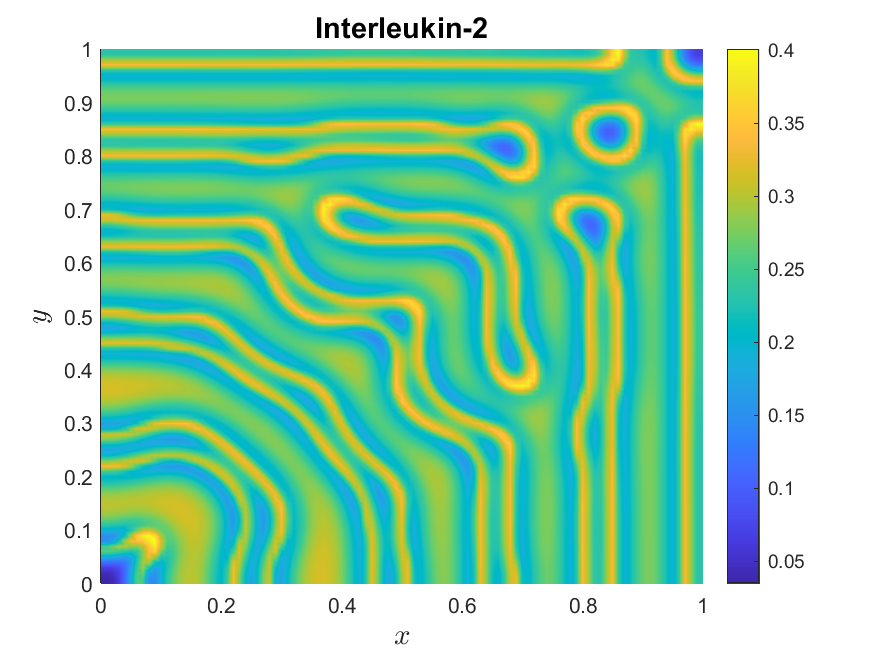}}
    		\caption{Spatial distributions of $u,\,v and\, w$ with $ d_{32}=0.01$ in the case of a treated patient.  The parameters are $c=0.25$, $p_2=0.5$, $d_{11}=0.001$, $d_{22}=1.99\cdot 10^{-5}$, $d_{33}=0.01$; the remaining ones are listed  in Table \ref{tab:par}}.
    		\label{con_ter_1}
	\end{figure}

Figures \ref{con_ter} and \ref{con_ter_1} show the Turing patterns in the three cell subpopulations over the 2D spatial domain, obtained at $t=1000$, for two different values of the cross-diffusion coefficient $d_{32}$. In agreement with the role played by effector cells and cytokines in the presence of tumor cells, it can be observed that the mixed spot–stripe patterns appears similarly across all three subpopulations. By comparing the figures related to treated and untreated patients (see Figs. \ref{sen_ter}, \ref{sen_ter_1} and Figs. \ref{con_ter},\ref{con_ter_1}), a clear difference emerges. Cold spots are transformed into intermediate stripe-like structures mixed with a few isolated cold spots when patients undergo therapy. Therefore, immunotherapy based on ACI and cytokine administration results in an improvement.
Such a transition is nearly reversed when the cross-diffusion coefficient takes a positive value (recalling that the coefficient $d_{32}<0$ induces the diffusion of cytokines toward regions with a high density of tumor cells). It can be clearly observed that hot region (yellow) containing only a few isolated cold spots, corresponding to high concentrations of tumor cells, evolve into stripes (blue), indicating lower tumor cell concentrations (compare Figs. \ref{sen_ter_1} and \ref{con_ter_1}).

\subsection{Hopf bifurcations: oscillatory solutions in untreated and treated patients}
\begin{figure}[ht!]
		\centering
    		\subfigure[]{\includegraphics[scale=0.22]{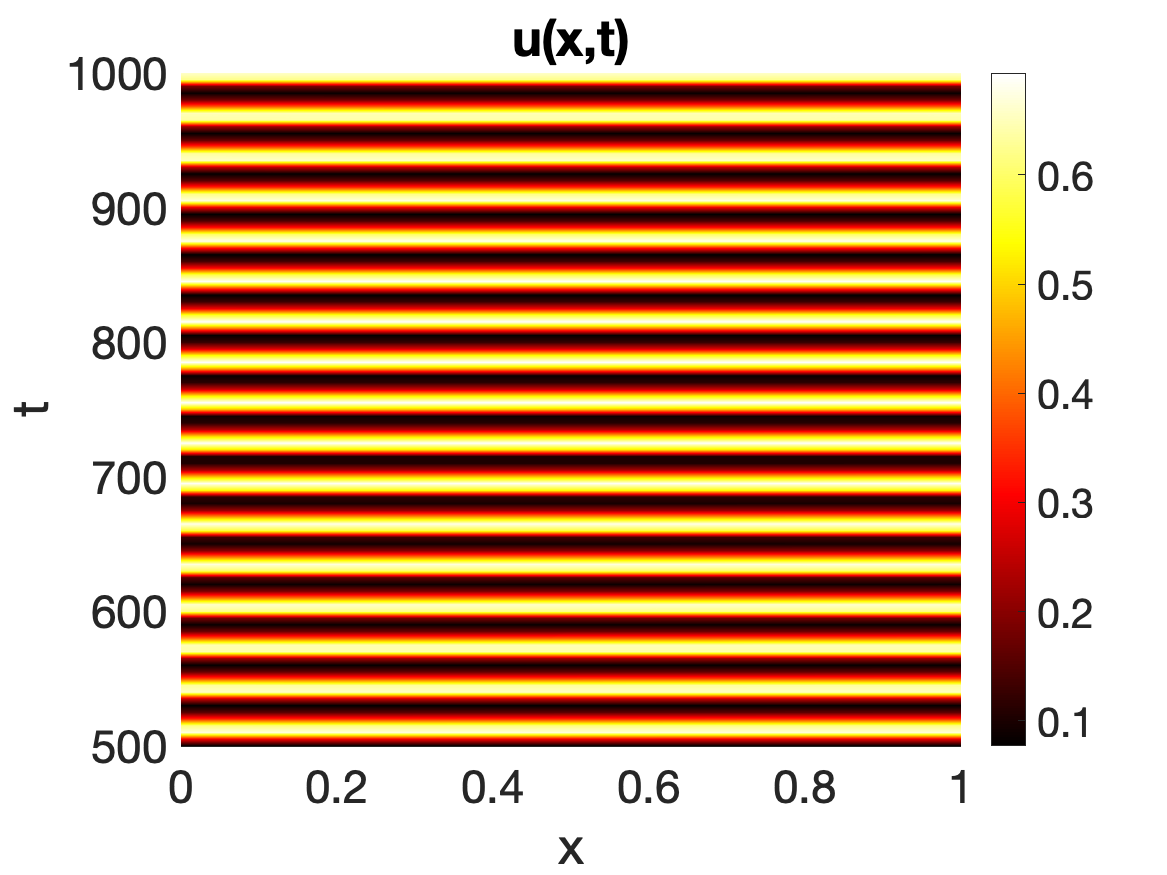}}
    		\subfigure[]{\includegraphics[scale=0.22]{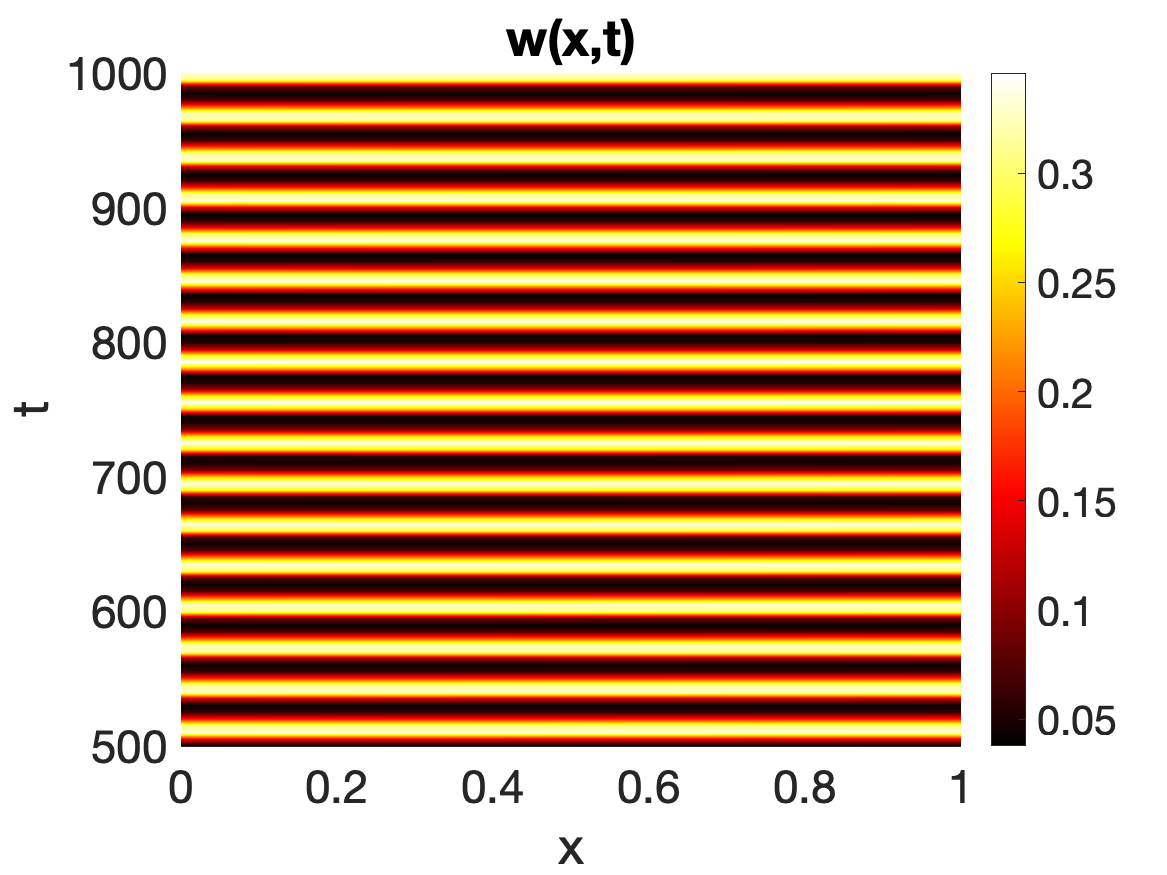}}
    		\subfigure[]{\includegraphics[scale=0.22]{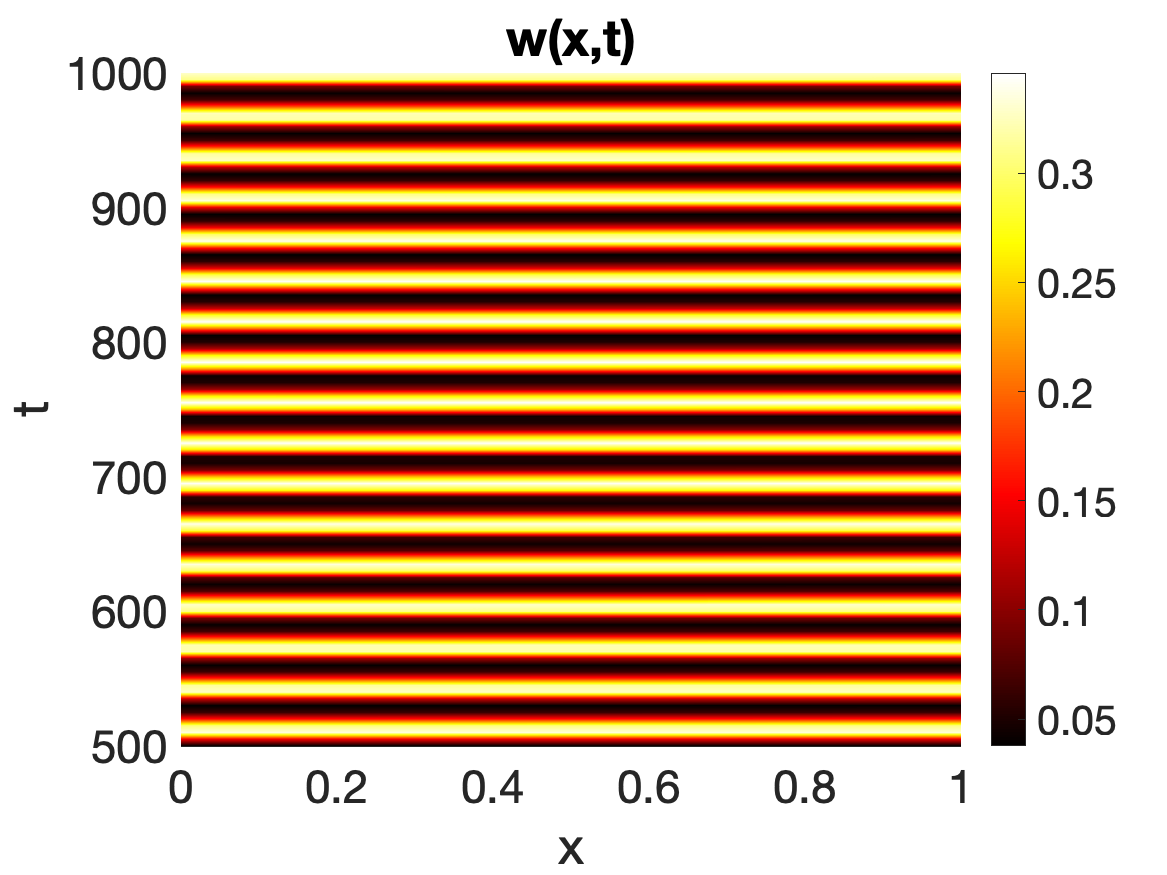}} \\
    		\subfigure[]{\includegraphics[scale=0.22]{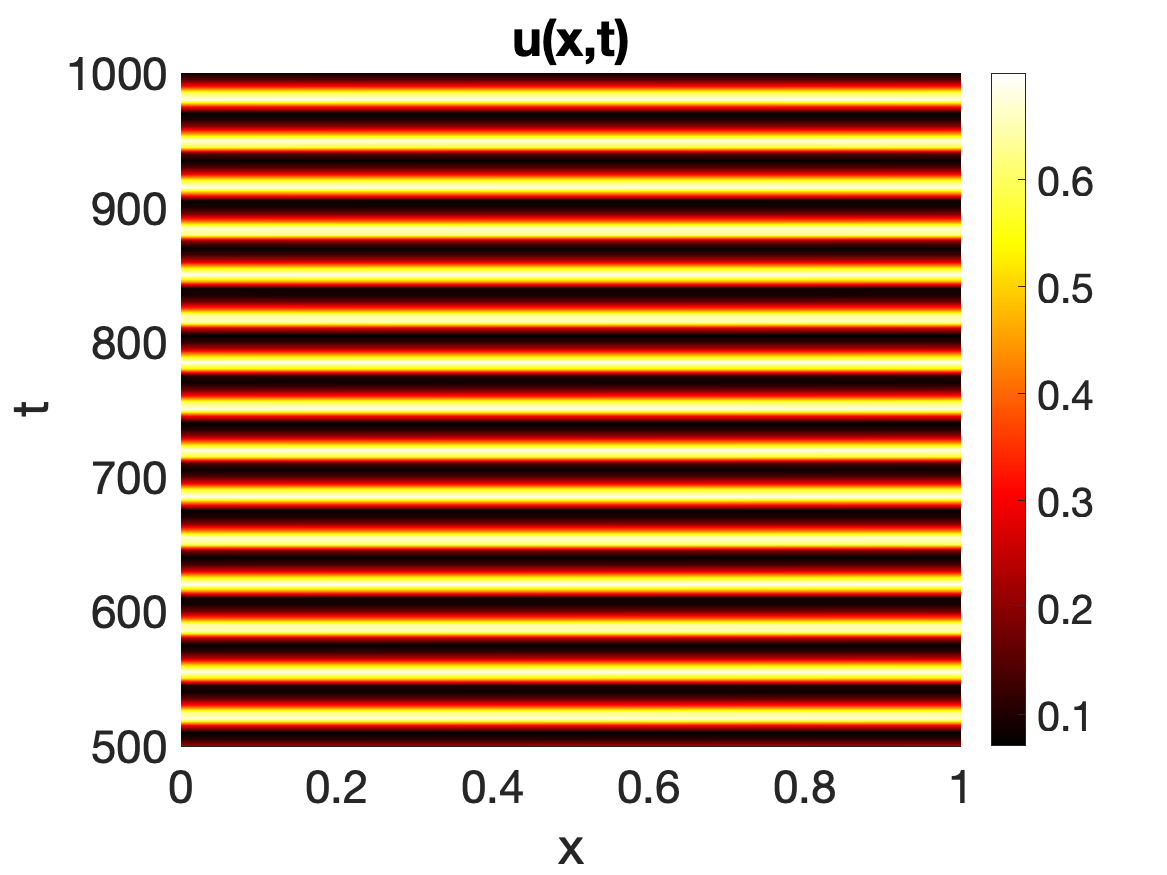}}
    		\subfigure[]{\includegraphics[scale=0.22]{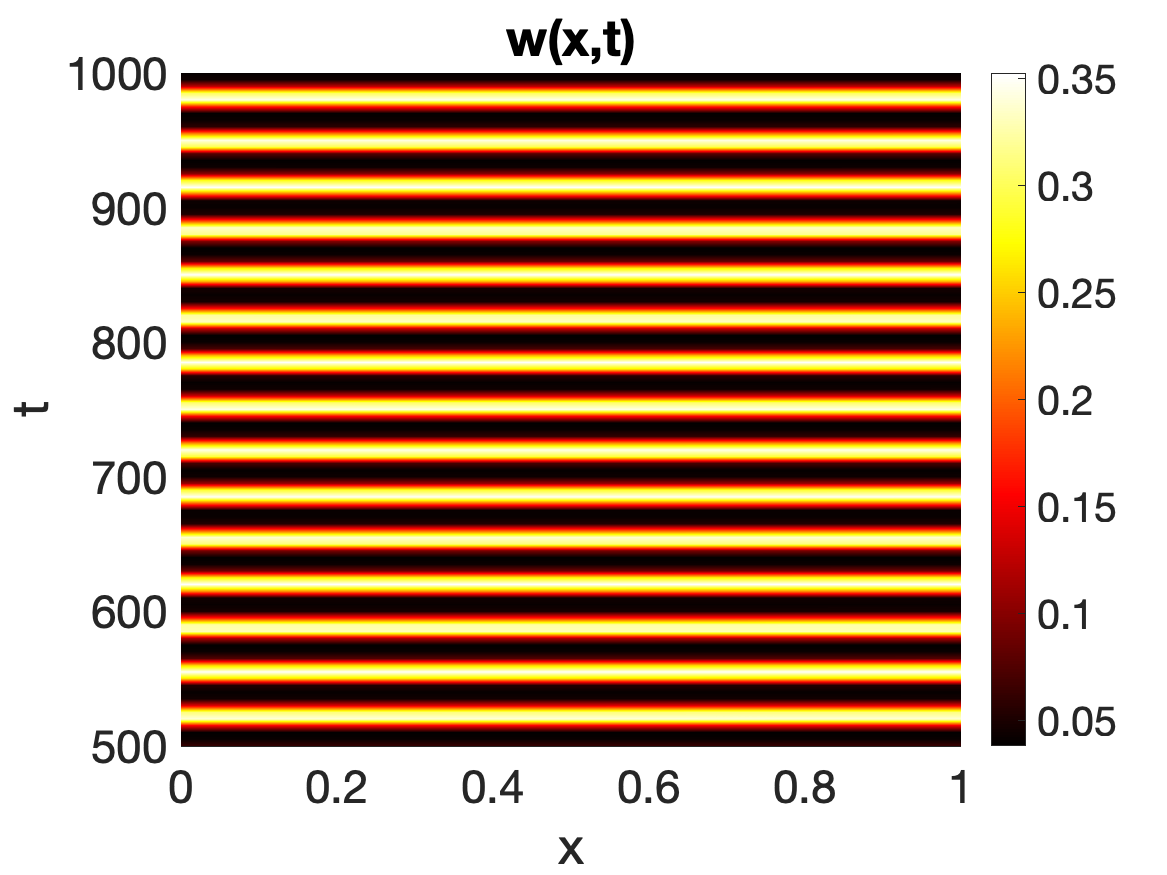}}
    		\subfigure[]{\includegraphics[scale=0.22]{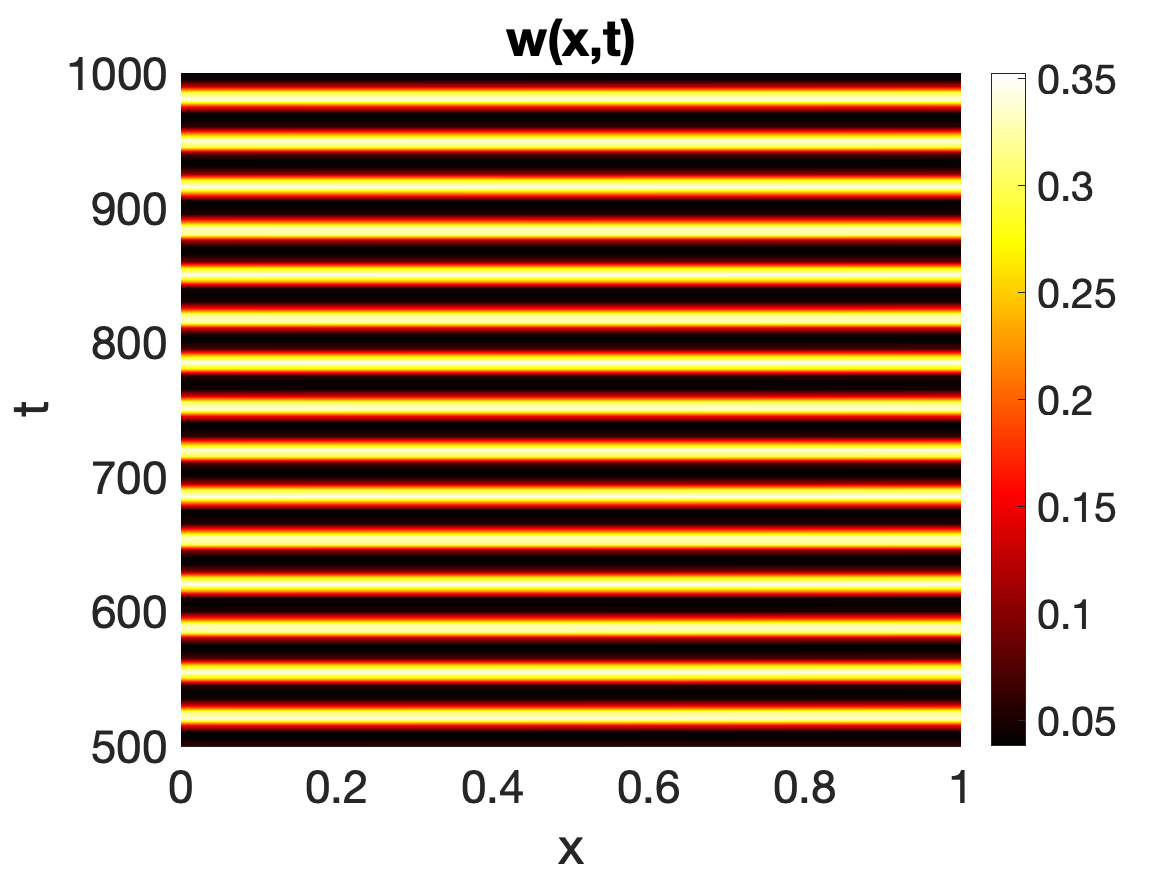}} 
    		\caption{Hopf bifurcation:  spatial distributions of $u,\,v,\, w$ for both untreated (first row) and treated (second row) patients. Parameters are set as follows: $c=0.25$, $p_2=0.55$, $d_{11}=0.001$, $d_{22}=0.00048$, $d_{33}=0.01$, $ d_{32}=-0.01$. The remaining  parameters are listed in Table \ref{tab:par}. }
      	 	\label{fig:hopf_space_time}
	\end{figure}

The reaction-diffusion system \eqref{eq:3} does not exhibit oscillatory Turing patterns, \emph{i.e.} the Turing-Hopf instability, for either $d_{22}=0.0000199$ or $d_{22}=0.00048$. However, in the latter case, $d_{22}=0.00048$, it is possible to identify certain frequencies $k$ for which the presence of a limit cycle is observed. However, this analysis can only be carried out numerically, as it is   analytically challenging to determine the critical parameter value at which the Hopf bifurcation occurs. This difficulty arises from the implicit form of the eigenvalue expressions  and the coexistence equilibrium point.  In this direction, some numerical simulations conducted in one spatial dimension (1D) highlight the emergence of time-periodic solutions (see the Figure \ref{fig:hopf_space_time}). These periodic solutions appear both in the presence and absence of therapeutic action.

	\begin{figure}[ht!]
		\centering
    		\subfigure[]{\includegraphics[scale=0.22]{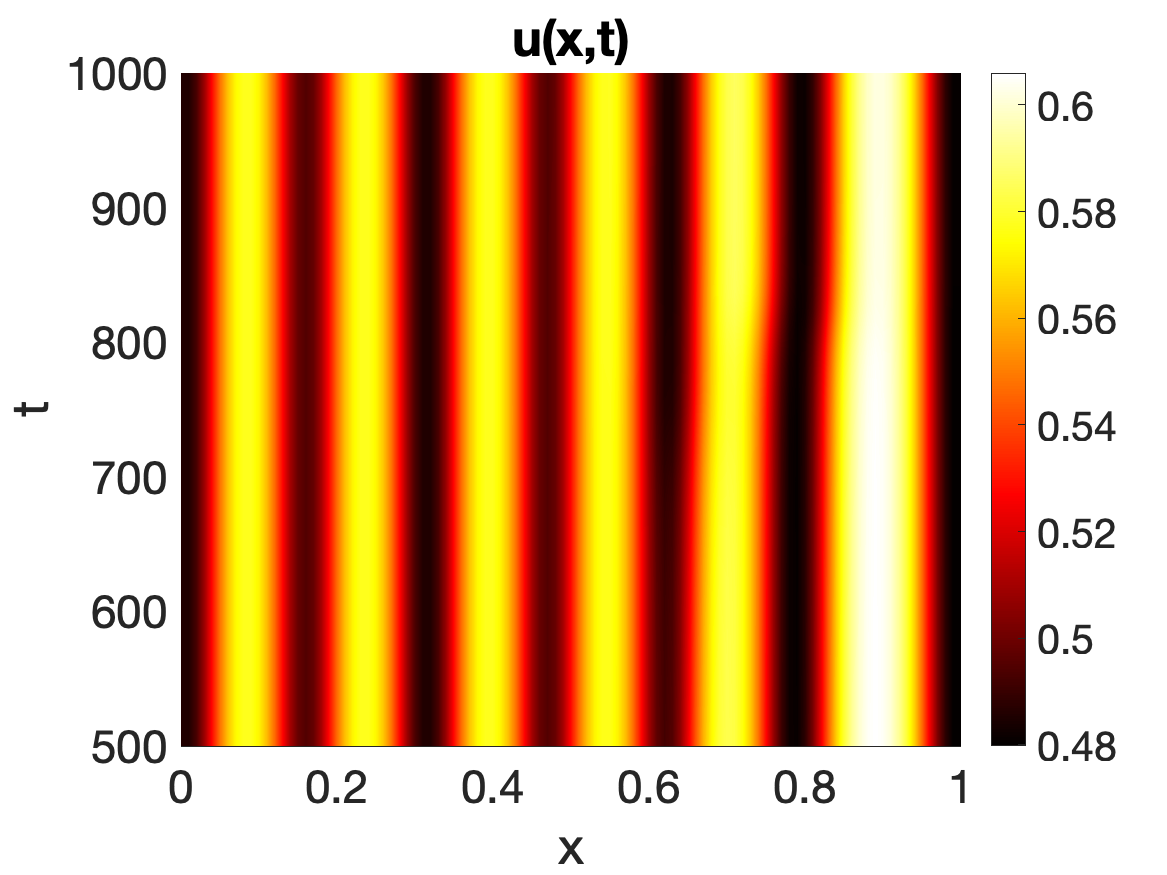}}
    		\subfigure[]{\includegraphics[scale=0.22]{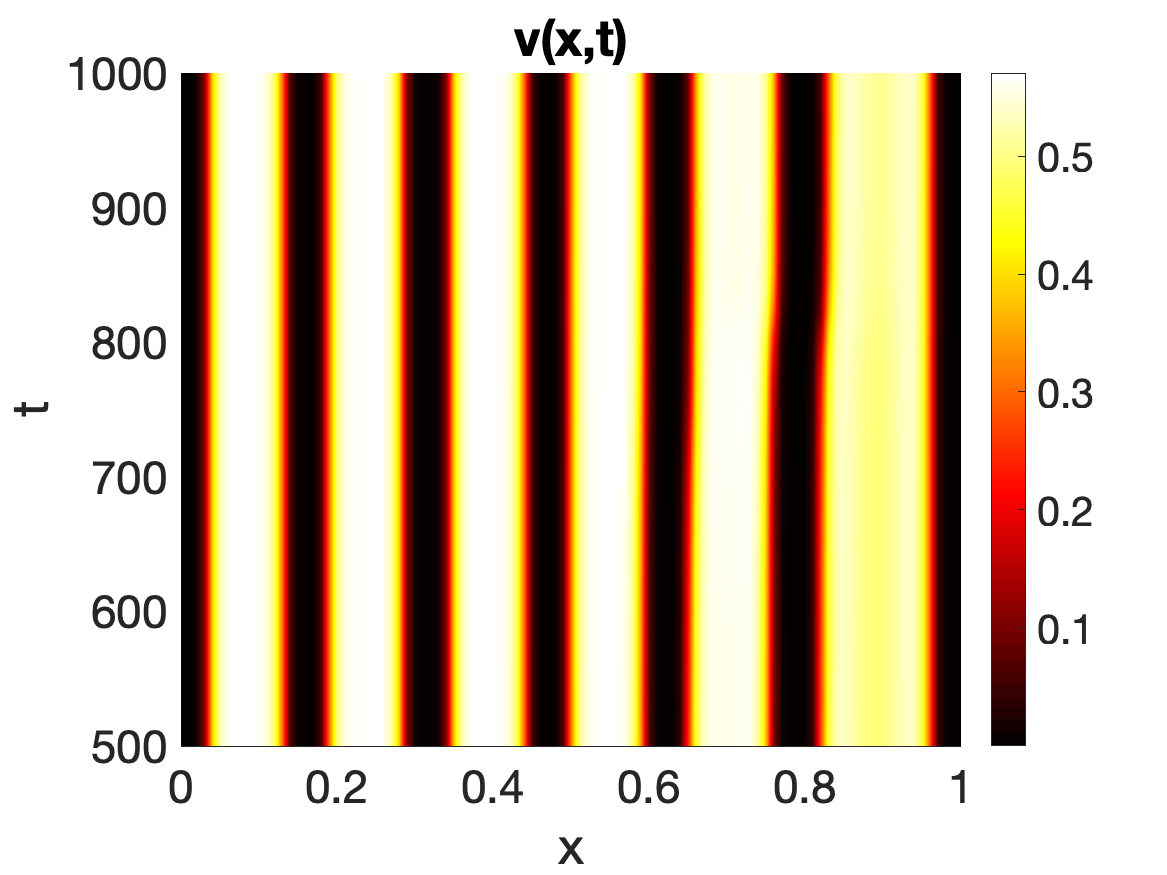}}
    		\subfigure[]{\includegraphics[scale=0.22]{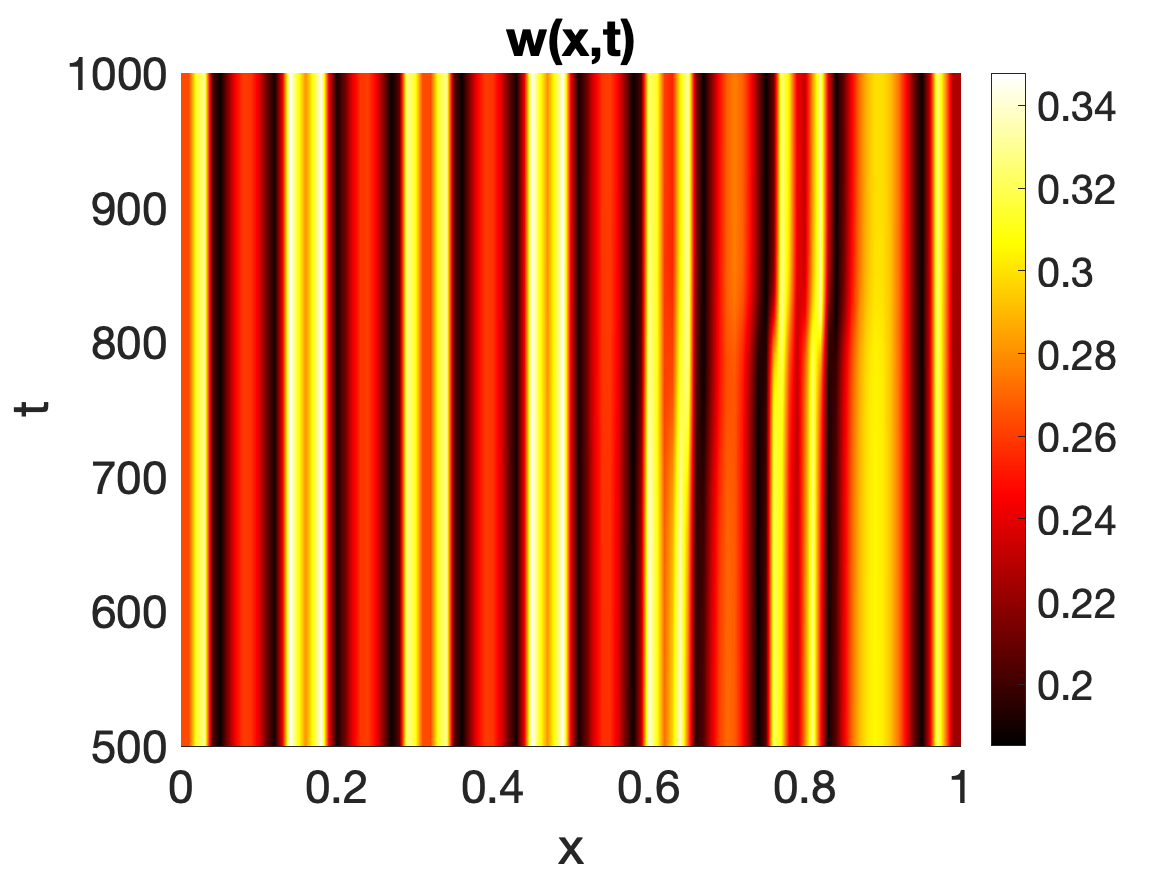}} \\
   		\subfigure[]{\includegraphics[scale=0.22]{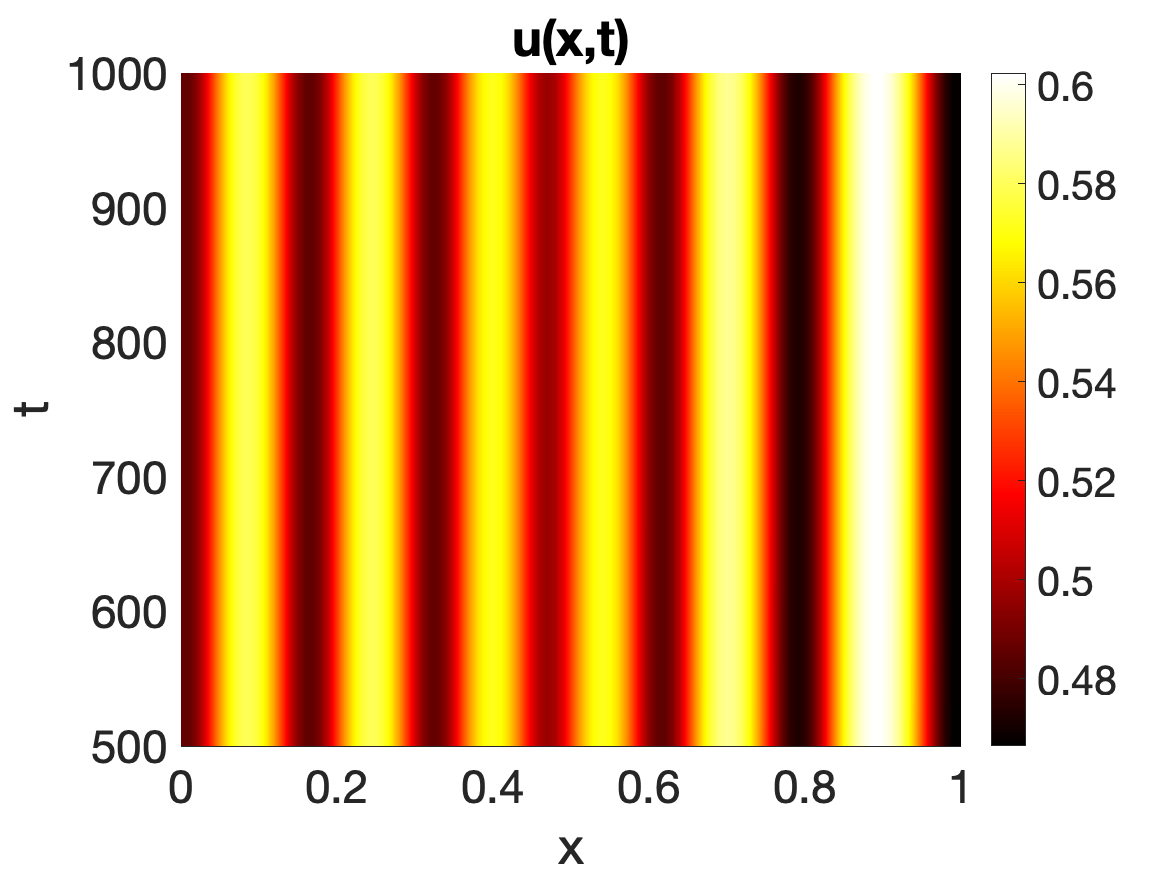}}
    		\subfigure[]{\includegraphics[scale=0.22]{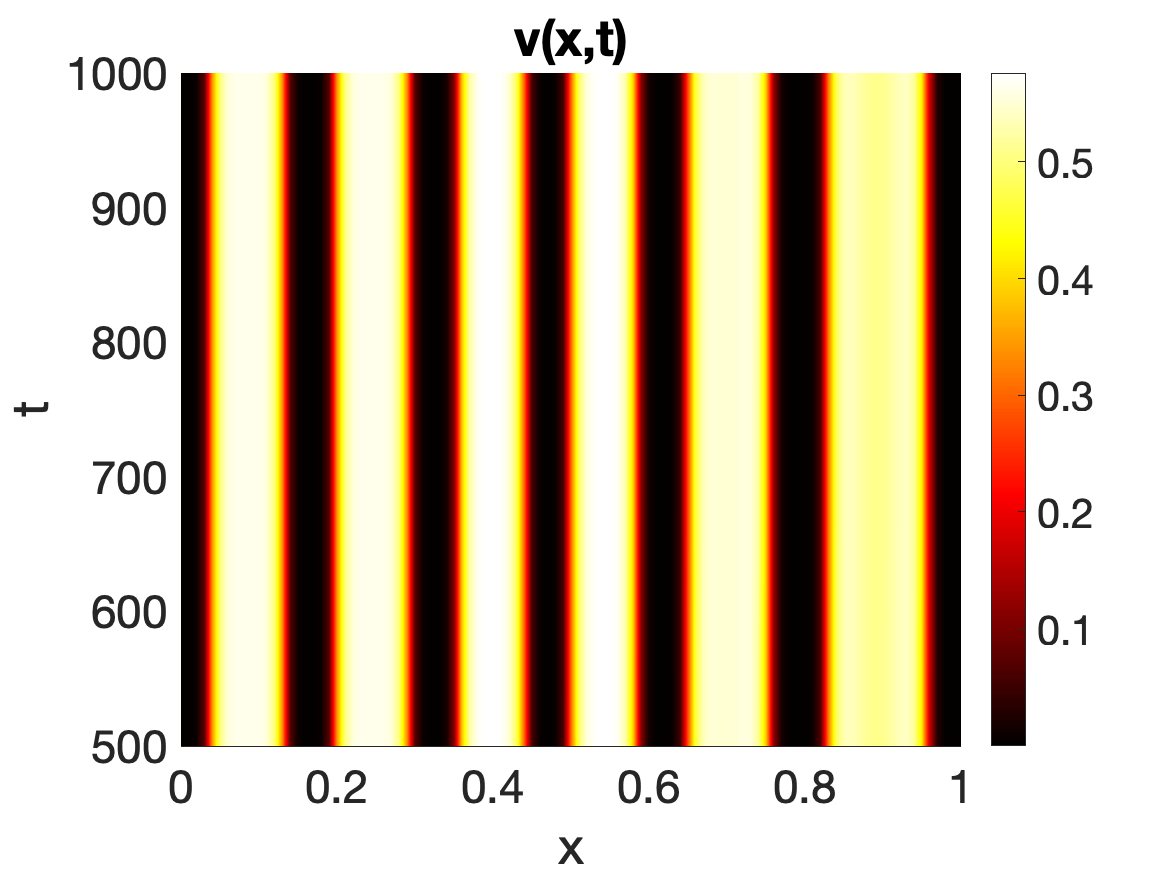}}
    		\subfigure[]{\includegraphics[scale=0.22]{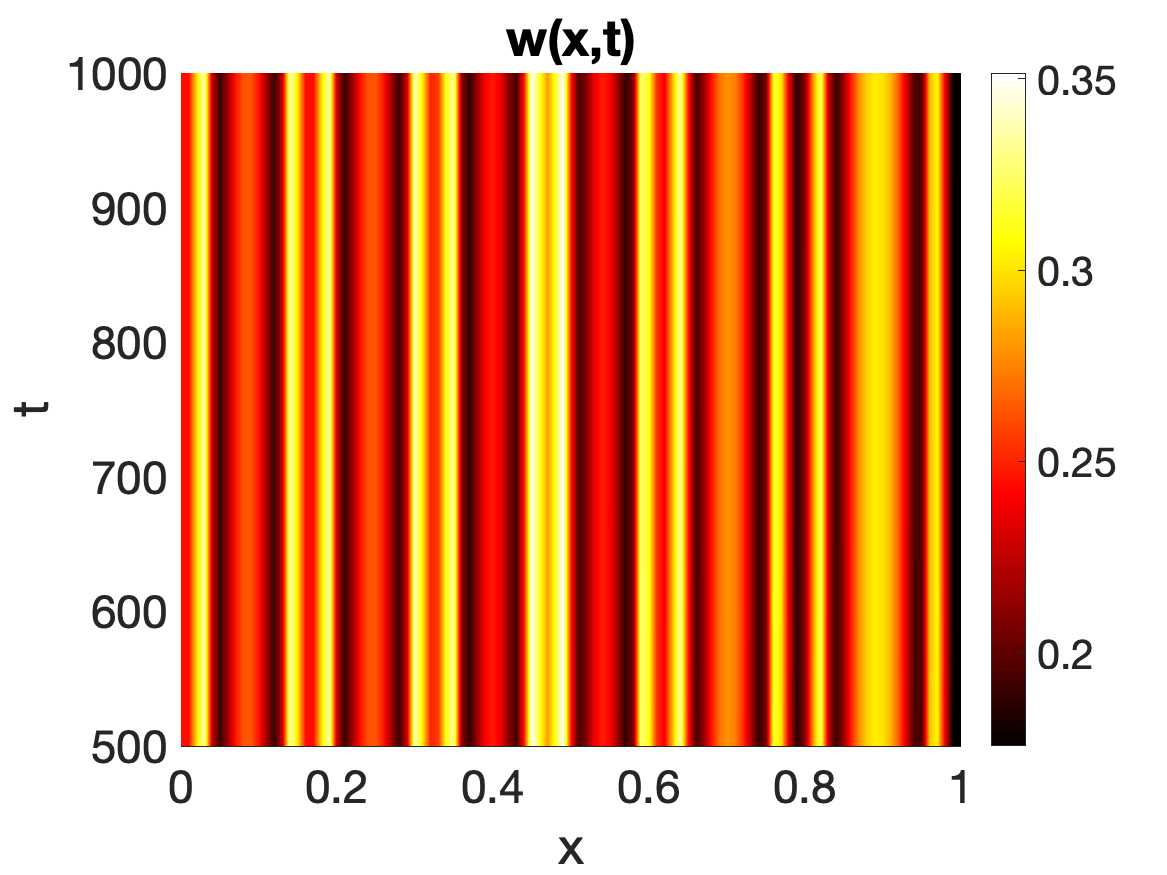}} 
    		\caption{Turing bifurcation for both untreated (first row) and treated (second row) patients. Spatial distributions of $u,\,v,\, w$ are shown with  $c=0.25$, $p_2=0.5$, $d_{11}=0.001$, $d_{22}=1.99\cdot 10^{-5}$, $d_{33}=0.01$, $ d_{32}=-0.01$. The remaining parameters are summarized in Table \ref{tab:par}. }
    	  	\label{fig:turing_1D}
	\end{figure}

\subsection{Results with COMSOL}
\begin{figure}[ht!]
        		\centering
            	\subfigure[]{\includegraphics[scale=0.3]{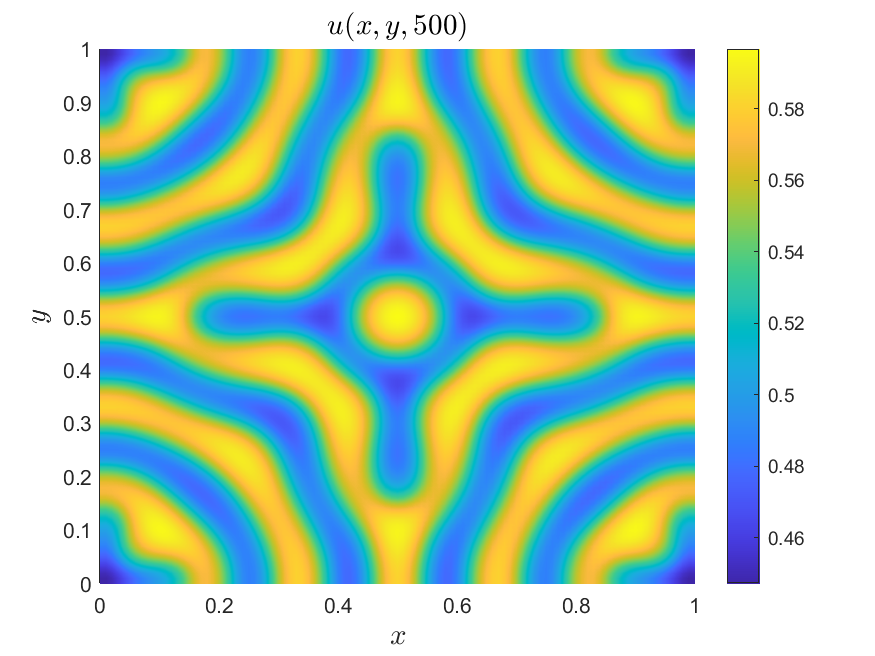}}
            	\subfigure[]{\includegraphics[scale=0.3]{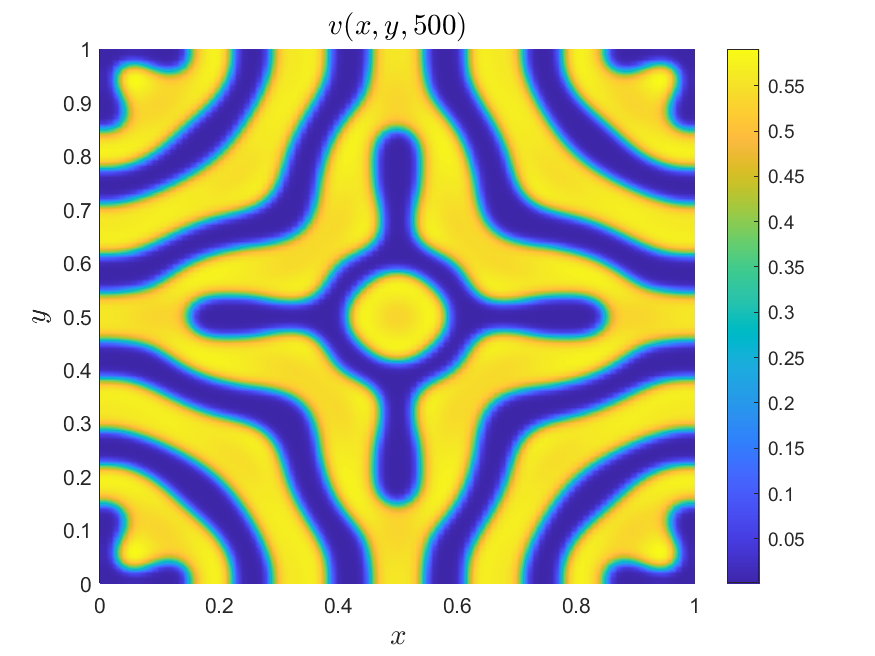}}
            	\subfigure[]{\includegraphics[scale=0.3]{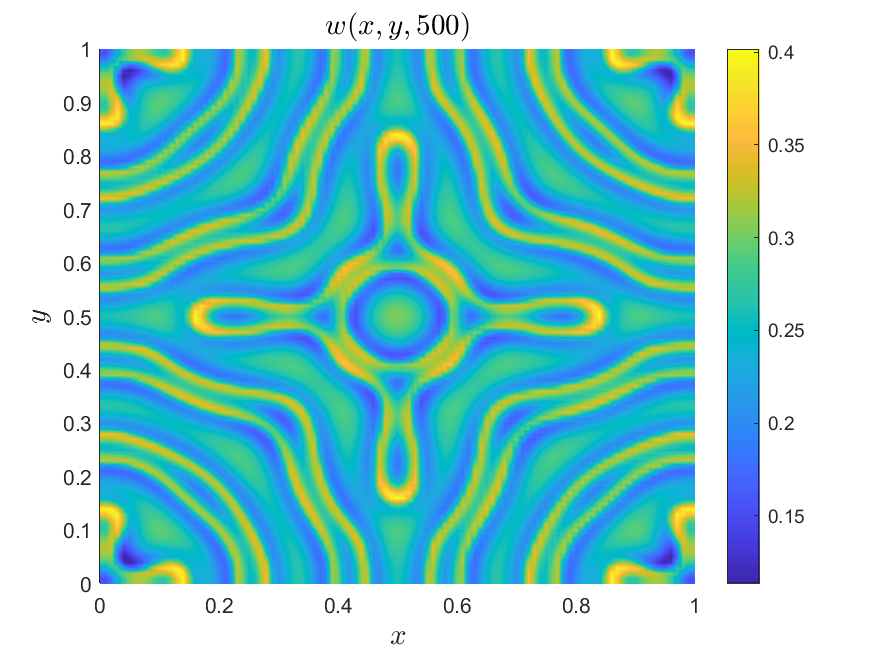}} \\
           	\subfigure[]{\includegraphics[scale=0.13]{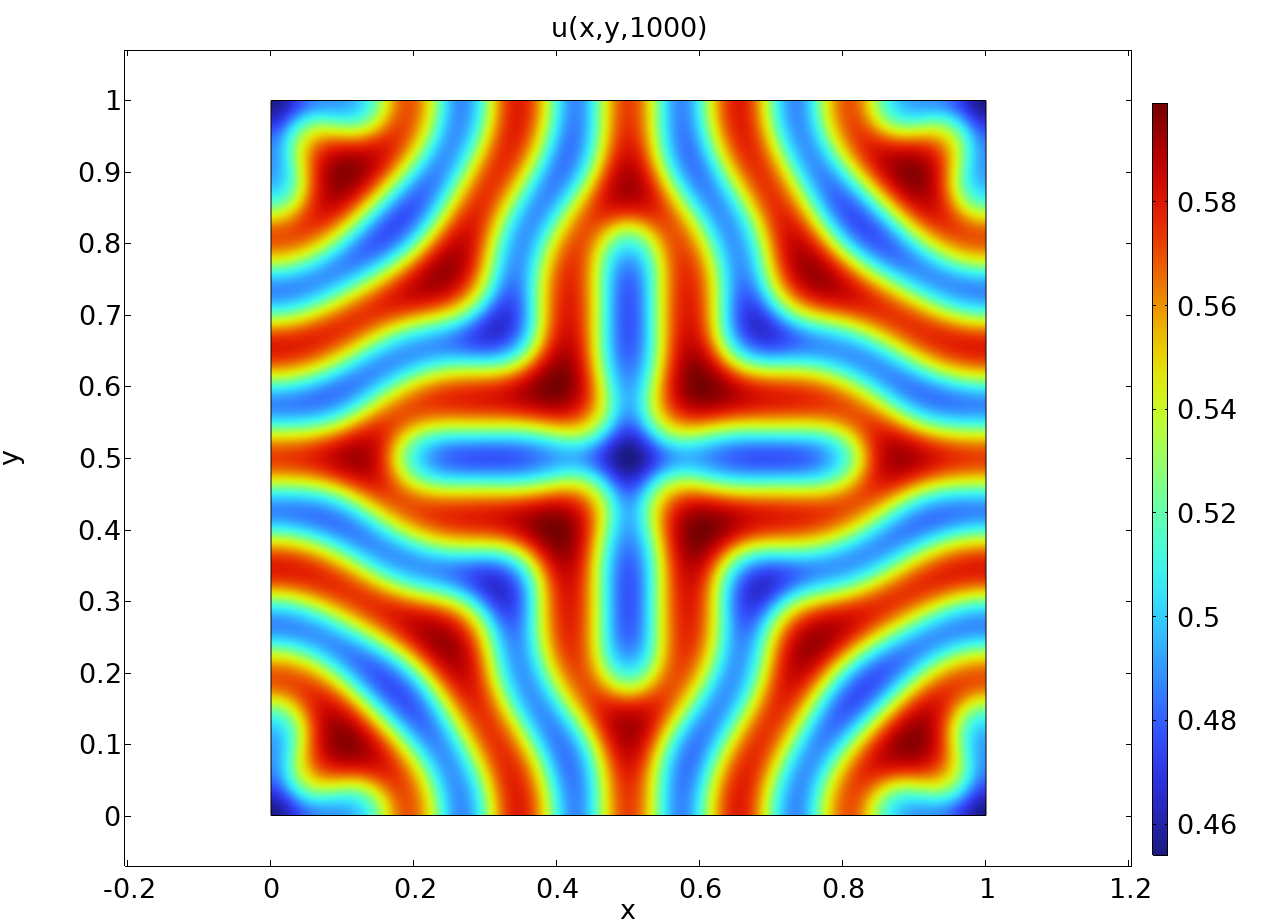}}
            	\subfigure[]{\includegraphics[scale=0.13]{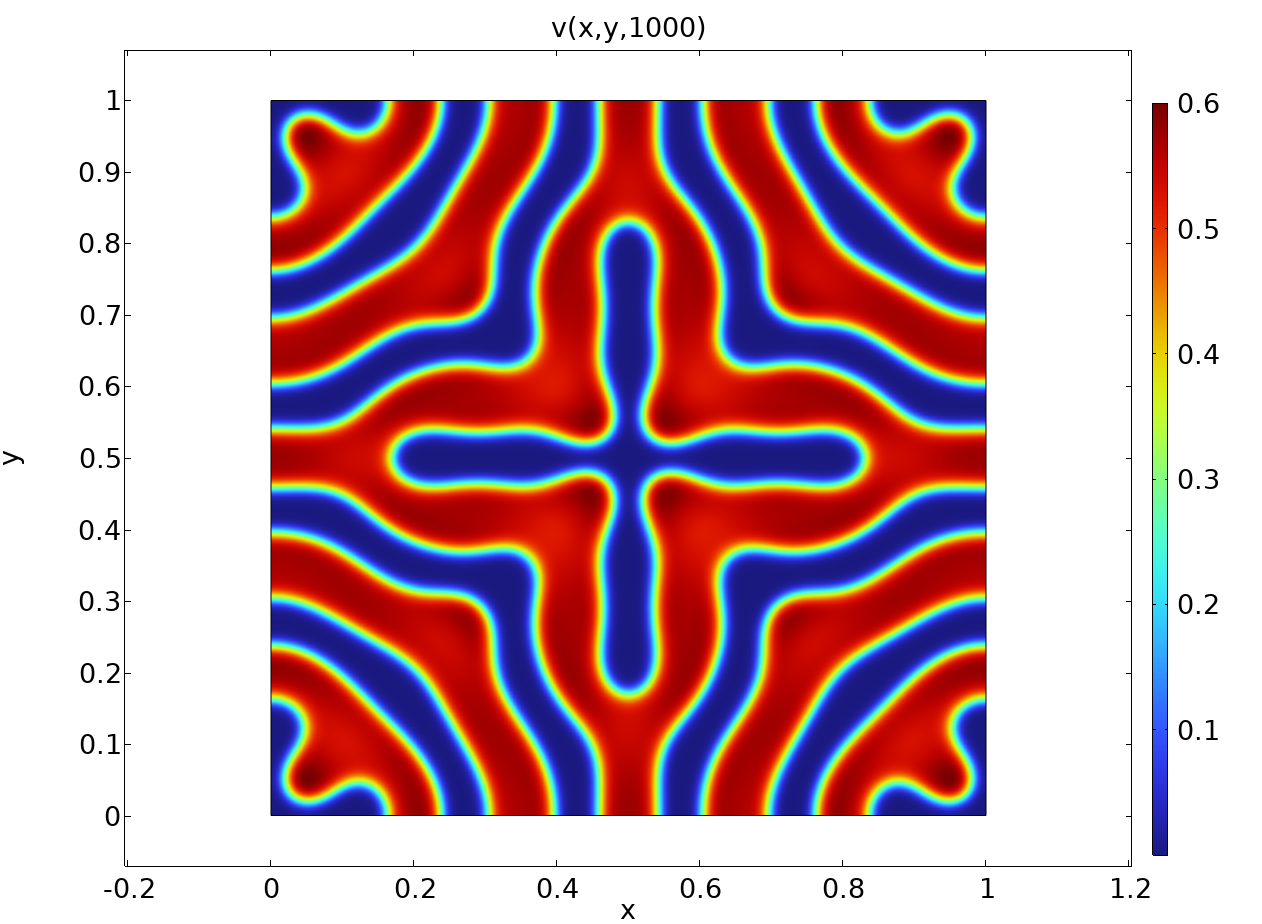}}
            	\subfigure[]{\includegraphics[scale=0.13]{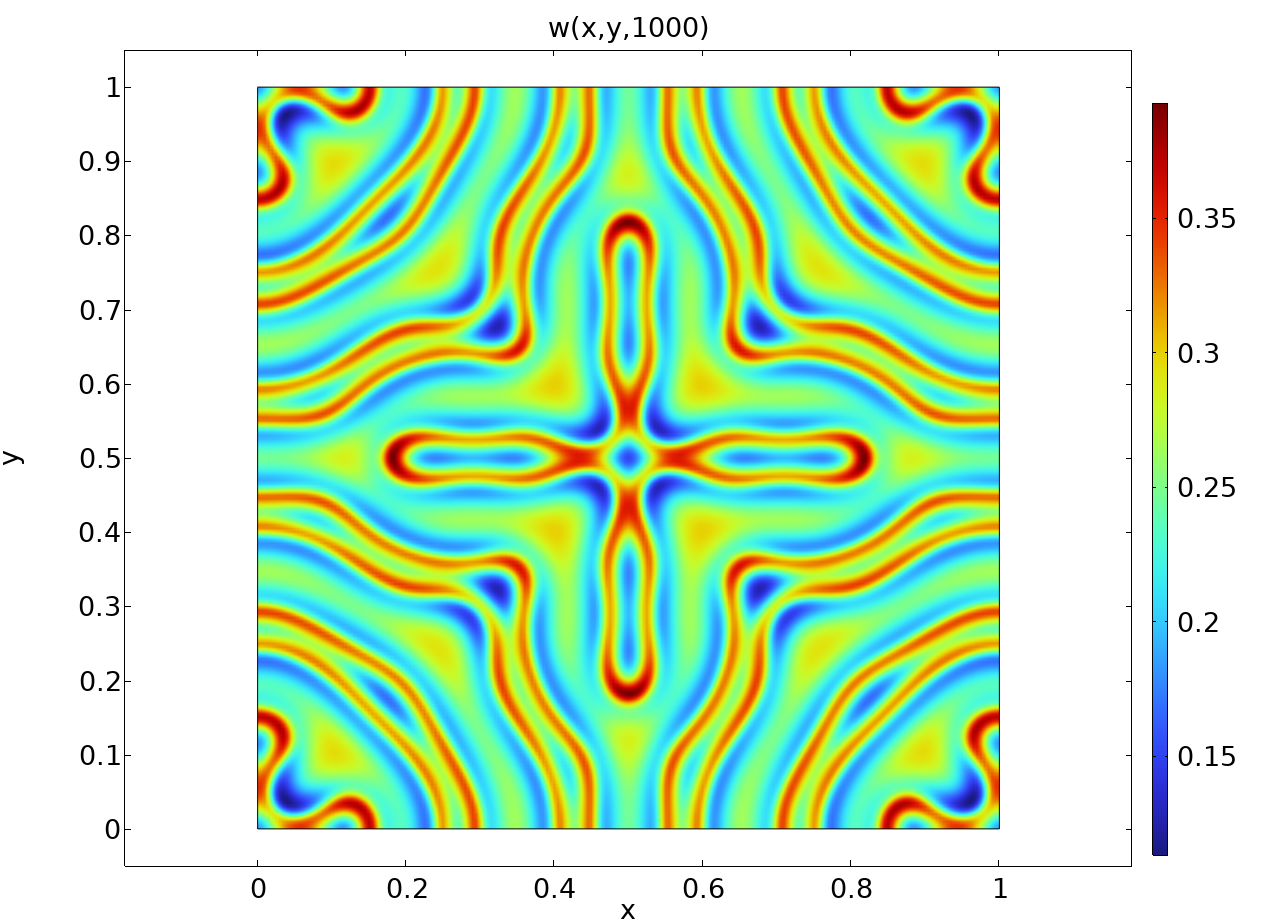}} \\
            	\subfigure[]{\includegraphics[scale=0.3]{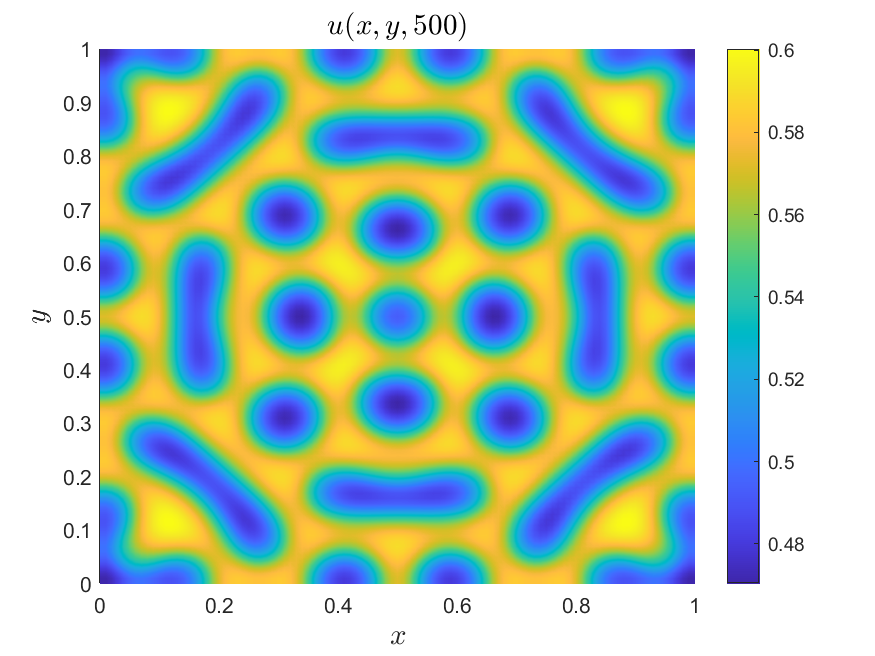}}
            	\subfigure[]{\includegraphics[scale=0.3]{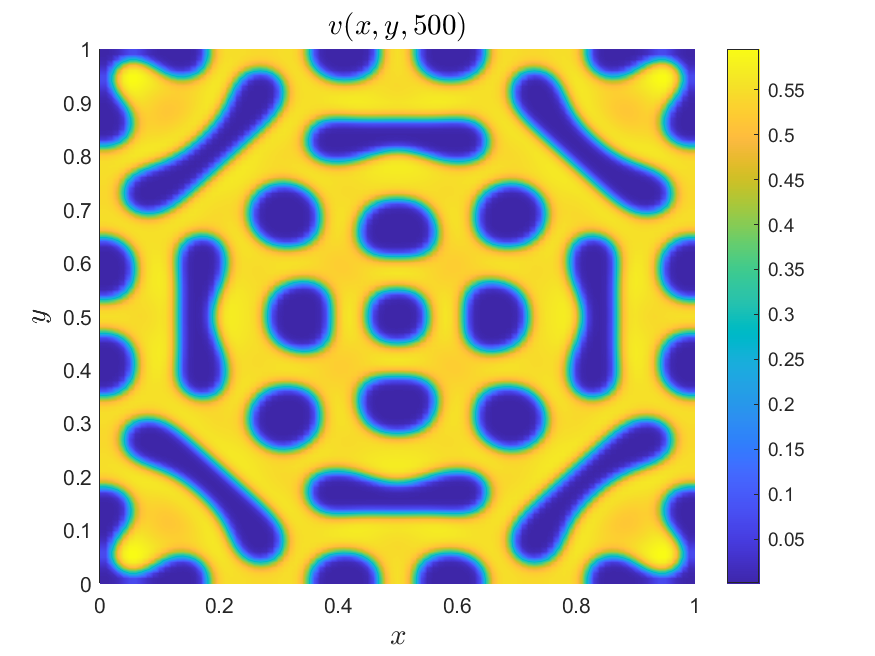}}
            	\subfigure[]{\includegraphics[scale=0.3]{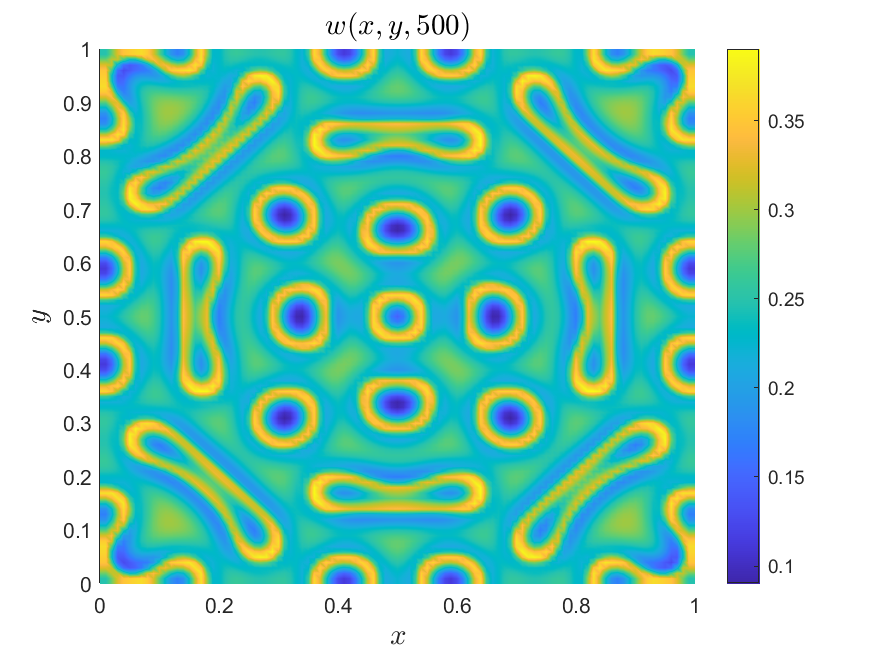}} \\
            	\subfigure[]{\includegraphics[scale=0.13]{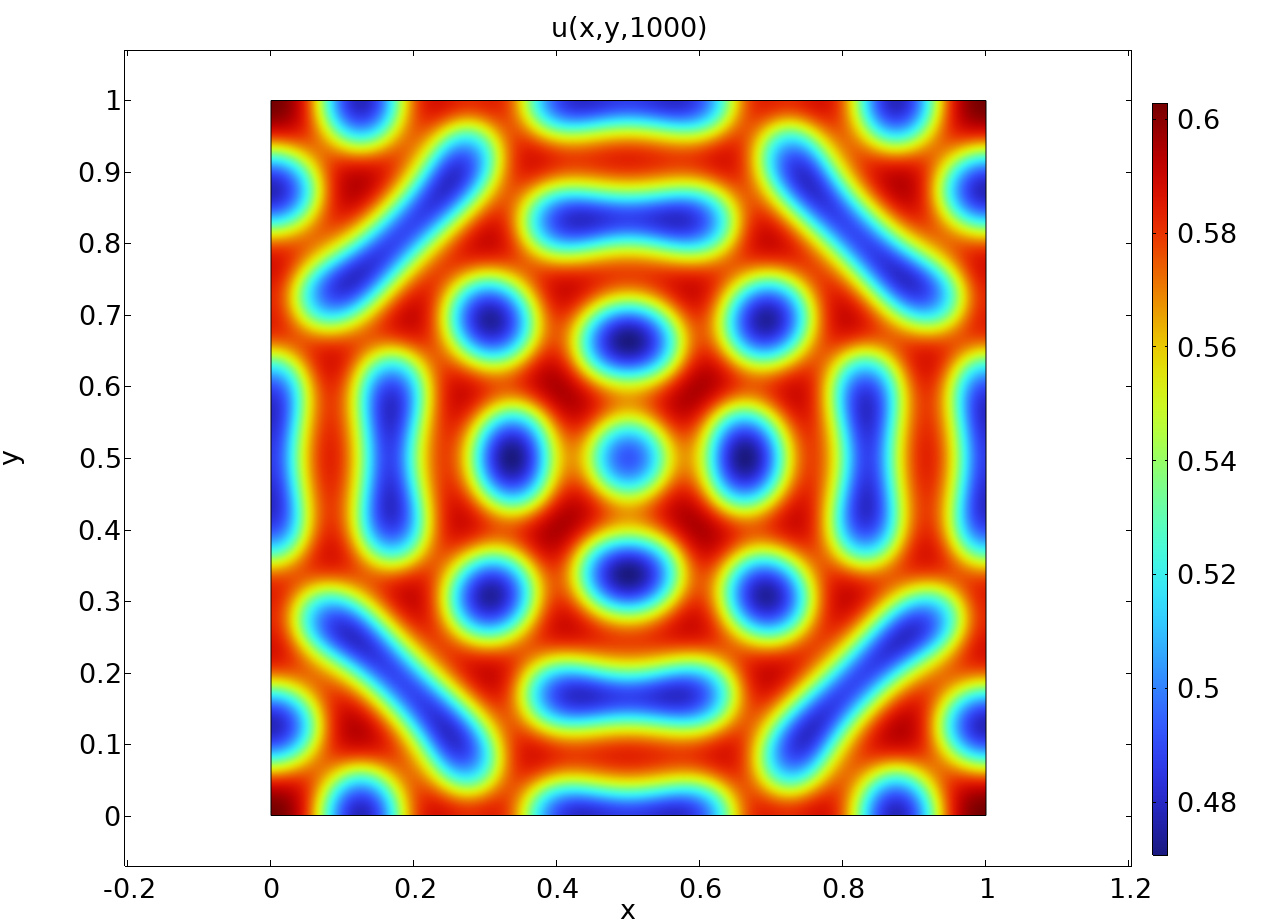}}
            	\subfigure[]{\includegraphics[scale=0.13]{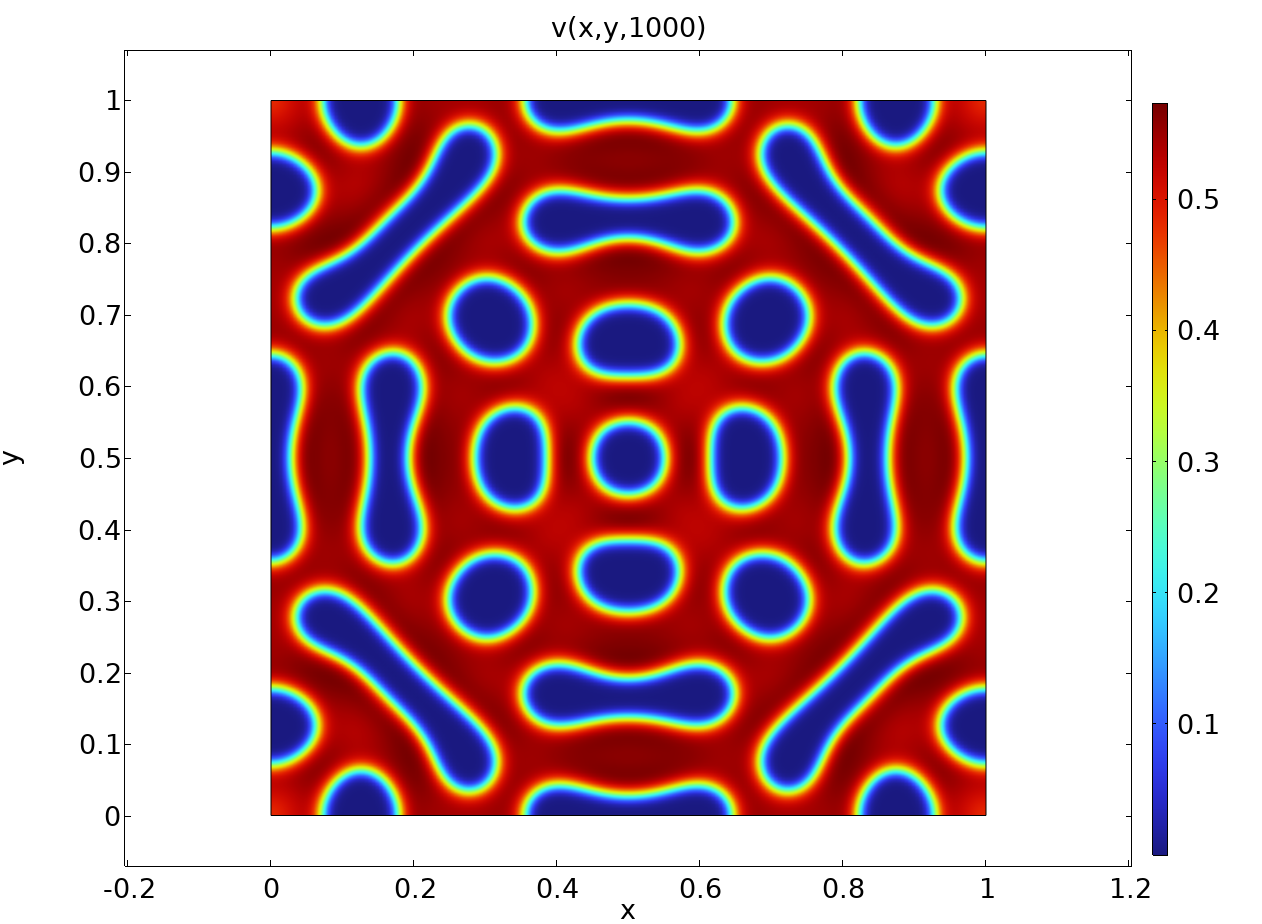}}
            	\subfigure[]{\includegraphics[scale=0.13]{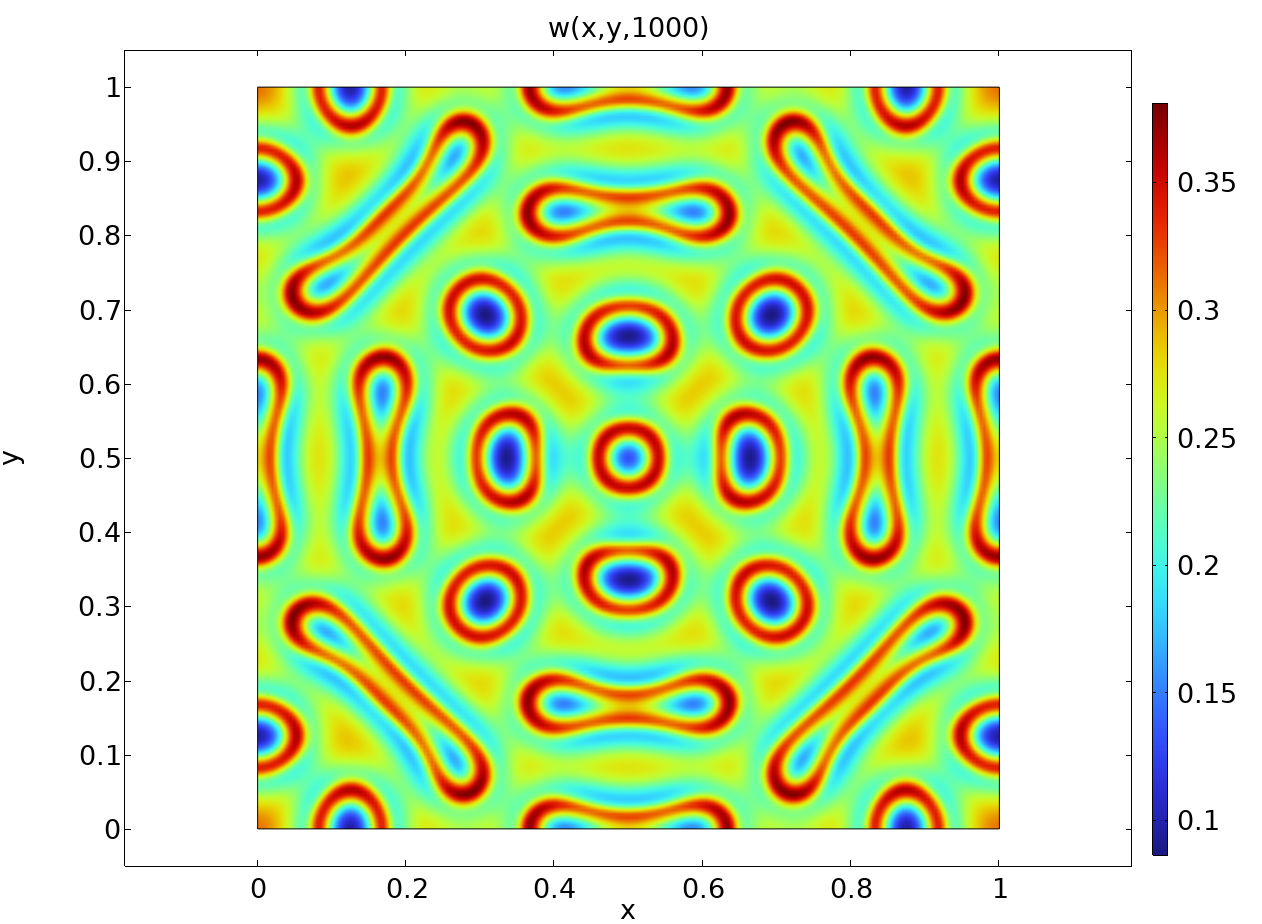}}
            	\caption{Spatial distributions of $u$, $v$, $w$ with $d_{32}=0.01$, both with and without therapy. Comparison between results obtained using the explicit finite difference scheme implemented in MATLAB and those obtained with COMSOL Multiphysics\textsuperscript{\textregistered}. The parameters are $c=0.25$, $p_2=0.5$, $d_{11}=0.001$, $d_{22}=1.99\cdot 10^{-5}$, $d_{33}=0.01$; the remaining ones are presented  in Table \ref{tab:par}}.
               \label{con_ter_2}
        \end{figure}

Numerical simulations obtained using COMSOL Multiphysics\textsuperscript{\textregistered} software are shown in Figure \ref{con_ter_2}, and the resulting patterns display the same spatial distribution as in the previous figures.

The 2D simulations of our immune–cancer model generated a rich variety of spot, stripe, and mixed patterns, whose stability properties are consistent with findings reported in previous studies.

\clearpage
\section{Conclusions}\label{sec_conclu}
In this paper, we implemented and analyzed a variant of a model originally proposed by Suddin et al., consisting of three coupled reaction-diffusion equations involving only self-diffusion terms. 

In particular,  we extended the model proposed in \cite{suddin2021reaction},  by incorporating the effects induced by the presence of cross-diffusion, specifically through the coefficient $d_{32}$.  From a biological perspective, it is particularly important to account for the migratory response of IL-2 cells, which are assumed to move (with $d_{32}\leq 0$) toward regions with a high concentration of tumor cells. This behavior is supported by biological reasoning and justifies the inclusion of cross-diffusion only for the IL-2 variable.

Cross-diffusion coefficients for the other components were not introduced, particularly  ($d_{23}$), as this would imply that tumor cells are either attracted to (with $d_{23}\leq 0$) or repelled by (with $d_{32}\geq 0$) 
 IL-2 cells. 
In our opinion, such an effect lacks experimental support and therefore does not appear to be biologically meaningful in this context. The inclusion of cross-diffusion
allows for the exploration of additional spatial interaction mechanisms between
species, leading to more biologically realistic pattern formation. From
a biomedical perspective, such mechanisms may provide useful insights into
the spatial organization of tumor cells and their microenvironment, which
can suggest more effective therapeutic strategies.

In our extended model, we derived numerical constraints on the parameters to guarantee the existence of an admissible, asymptotically stable, spatially homogeneous coexistence equilibrium, which may lose stability due to both self- and cross-diffusion terms. We analytically demonstrated that the model can exhibit Turing instability depending on the value of a control parameter. 

Numerical simulations reveal the emergence of characteristic patterns that remain stationary over time. These spatial structures are influenced by the therapy administered to patients. The resulting patterns show that cancer and  effector cells tend to concentrate in the same regions, while the IL-2 compound predominantly accumulates at the edges of the patterns formed by the cells.

In the presence of therapy, significant changes in  spatial organization are observed. The three cellular populations mainly form cold stripes localized  in the same region of the domain. Moreover, the diffusion of interleukins is directed toward areas of high tumor cell density when the cross-diffusion coefficient is negative, and toward areas with low density when the coefficient is positive.

Although we currently lack rigorous mathematical results establishing conditions for the positivity of solutions, the chosen parameters and the numerical solutions provided here describe biologically meaningful scenarios.

Based on our findings regarding cross-diffusion in cancer–cytokine cell populations, we anticipate that this research will inspire further empirical and theoretical studies in the field.

\paragraph*{Acknowledgements}
Work supported by the Italian National Group of Mathematical Physics (GNFM-IN$\delta$AM). 



\begin{thebibliography}{99}

	\bibitem{curti1996influence}
	B.D. Curti, A.C. Ochoa, W.J. Urba, W.G. Alvord, W.C. Kopp, G. Powers, C. Hawk, S.P. Creekmore, B.L. Gause, J.E. Janik, J.T. Holmlund, P. Kremers, R.G. Fenton, L. Miller,  Sznol, J. W. Smith II, W.H. Sharfman and D.L. Longo,
	\textit{Influence of interleukin-2 regimens on circulating populations of lymphocytes after adoptive transfer of anti-CD3-stimulated T cells: Results from a phase I trial in cancer patients},
	Journal of Immunotherapy, 
	pp. 296--308,
 	1996.
 
	\bibitem{gause1996phase}
	B.L. Gause, M. Sznol, W.C. Kopp, J.E. Janik, J.W. Smith II, R.G. Steis, W.J. Urba, W. Sharfman, R.G. Fenton, S.P. Creekmore, J. Holmlund, K.C. Conlon, L.A. VanderMolen, D.L. Longo,
	\textit{Phase I study of subcutaneously administered interleukin-2 in combination with interferon alfa-2a in patients with advanced cancer},
	 Journal of Clinical Oncology,
	14(8), 
	pp. 2234--2241,
	1996.
	
	\bibitem{hara1996rejection}
	I. Hara, H. Hotta, N. Sato, H. Eto, S. Arakawa, S. Kamidono,
	\textit{Rejection of mouse renal cell carcinoma elicited by local secretion of interleukin-2},
	Japanese journal of cancer research,
	87(7), 
	pp. 724--729,
	1996.
	
	\bibitem{KusSan}
	F. Adi-Kusumo and R. S. Winanda,
	\textit{Bifurcation Analysis of the cervical cancer cells, effector cells, and IL-2 compounds interaction model with immunotherapy},
	Far East Journal of Mathematical Sciences (FJMS),
	99(6), 
	pp. 869--883, 
	2016.

	\bibitem{kaempfer1996prediction}
	R. Kaempfer, L. Gerez, H. Farbstein, L. Madar, O. Hirschman, R. Nussinovich, A. Shapiro,
	\textit{Prediction of response to treatment in superficial bladder carcinoma through pattern of interleukin-2 gene expression},
	Journal of Clinical Oncology,
	14(6), 
	pp. 1778--1786, 
	1996.

	\bibitem{keilholz1994immunotherapy}
	U. Keilholz, C. Scheibenbogen, E. Stoelben,  H.D. Saeger, W. Hunstein,
	\textit{Immunotherapy of metastatic melanoma with interferon-$\alpha$ and interleukin-2: pattern of progression in responders and patients with stable disease with or without resection of residual lesions},
	European Journal of Cancer,
	30(7),
	pp. 955--958, 
	1994.

	\bibitem{valle2016pleiotropic}
	A. Valle-Mendiola, A. Guti{\'e}rrez-Hoya, M. Lagunas-Cruz, B. Weiss-Steider, I. Soto-Cruz,
	\textit{Pleiotropic Effects of IL-2 on Cancer: Its Role in Cervical Cancer},
	Mediators of inflammation,
	2016(1),
	pp. 2849523, 
	2016.
	
	\bibitem{kirschner1998modeling}
	D.~Kirschner and J.C.~Panetta, 
	\textit{Modeling immunotherapy of the tumor–immune interaction}, 
	Journal of Mathematical Biology,
	37(3),
	pp. 235-252, 
	1998.

	\bibitem{adam1997survey}
	J.A. Adam,and N. Bellomo,
	\textit{A survey of models for tumor-immune system dynamics},
	Springer Science and Business Media, 
	2012.

	\bibitem{bagheri2021}
	Y.~Bagheri, A.~Barati, A.~Aghebati‐Maleki, L.~Aghebati‐Maleki, and M.~Yousefi,
	\textit{Current progress in cancer immunotherapy based on natural killer cells},
	Cell Biology International, 
	45,
	pp. 2--17,
	2021.

	\bibitem{wang2021turing}
	J.~Wang and S.~Liu, 
	\textit{Turing and Hopf Bifurcation in a Diffusive Tumor-immune Model}, 
	Journal of Nonlinear Modeling and Analysis, 
	3(3),
	pp. 477--493,
	2021.

	\bibitem{wangberg2023}
	Y.~Wang, D.R.~Bergman, E.~Trujillo \textit{et al.}, 
	\textit{Mathematical model predicts tumor control patterns induced by fast and slow cytotoxic T lymphocyte killing mechanisms}, 
	Scientific Reports, 
	13,
	pp. 22541, 
	2023.

	\bibitem{turing1990chemical}
	A. M. Turing,
	\textit{The chemical basis of morphogenesis}, 
	Bulletin of mathematical biology, 
	52(1), 
	pp. 153--197,
	1990.
	
	
	\bibitem{Murray2}  
	J.~D. Murray
	\textit{Mathematical biology II: spatial models and biomedical applications},
 	Springer, New York,
	2003.
	
	
	\bibitem{LamLou} 
	K.-Y. Lam, Y.~Lou,
	\textit{Introduction to Reaction-Diffusion Equations: Theory and
  Applications to Spatial Ecology and Evolutionary Biology, Lecture Notes on
  Mathematical Modelling in the Life Sciences},
 	Springer, New York,
	2022.
		

	
	\bibitem{kan1998singular}  
	Y. Kan-On, M. Mimura,
	\textit{Singular perturbation approach to a 3-component reaction-diffusion system arising in population dynamics},
  	SIAM journal on mathematical analysis,
  	29(6),
  	pp. 1519--1536,
	1998.
		
	\bibitem{sherratt2013pattern}
	J.A. Sherratt,
	\textit{Pattern solutions of the Klausmeier model for banded vegetation in semiarid environments V: the transition from patterns to desert},
  	SIAM journal on mathematical analysis,
  	73(4),
  	pp. 1347--1367,
	2013.
		
	\bibitem{carfora2025turing}
	M.F.~Carfora, F.~Iovanna, and I.~Torcicollo,  
	\textit{Turing patterns in an intraguild predator--prey model}, 
	Mathematics and Computers in Simulation, 
	232,
	pp. 192--210, 
	2025.
	
	\bibitem{ali2025turing}
	G. Al{\`\i}, and I. Torcicollo,
	\textit{Turing pattern formation in a specialist predator--prey model with a herd-Holling-type II functional response},
	Mathematical Methods in the Applied Sciences,
  	48(1),
  	pp. 731--747,
  	2025,
	
	\bibitem{ali2026turing}
	G. Al\'i, C. Scuro, and I. Torcicollo,
	\textit{Pattern formation driven by cross-diffusion in the Klausmeier-Gray-Scott model}, 
	Mathematics and Computers in Simulation, 
	239,
	pp. 555-571, 
	2026.	
		
	\bibitem{tulumello2014cross}
	E. Tulumello, M.C. Lombardo, M. Sammartino,
	\textit{Cross-diffusion driven instability in a predator-prey system with cross-diffusion},
	Acta Applicandae Mathematicae,
  	132(1),
  	pp. 621--633,
	2014.
	
  

	\bibitem{giunta2021pattern}
	V. Giunta, M.C. Lombardo, M. Sammartino,
	\textit{Pattern formation and transition to chaos in a chemotaxis model of acute inflammation},
  	SIAM Journal on Applied Dynamical Systems,
 	20(4),
  	pp. 1844--1881,
	2021.
	
	
	\bibitem{lombardo2017demyelination} 
	M.C. Lombardo, R. Barresi, E. Bilotta, F. Gargano, P. Pantano, M. Sammartino,
	\textit{Demyelination patterns in a mathematical model of multiple sclerosis},
  	Journal of mathematical biology,
 	75(2),
  	pp. 373--417,
	2017.
	


	\bibitem{wang2012complex} 
	W. Wang, Y. Cai, M. Wu, K. Wang, Z. Li,
	\textit{Complex dynamics of a reaction--diffusion epidemic model},
  	Nonlinear Analysis: Real World Applications,
 	13(5),
  	pp. 2240--2258,
	2012.
	
	\bibitem{avila2022dynamics} 
	E. Avila-Vales, A.G.C. P{\'e}rez,
	\textit{Dynamics of a reaction--diffusion SIRS model with general incidence rate in a heterogeneous environment},
  	Zeitschrift f{\"u}r angewandte Mathematik und Physik,
 	73(1),
  	pp. 9,
	2022.
	
  	\bibitem{della2022mathematical}
	R.~Della Marca, R.~Machado Ramos, C.~Ribeiro, and A.~J.~Soares,  
	\textit{Mathematical modelling of oscillating patterns for chronic autoimmune diseases}, 
	Mathematical Methods in the Applied Sciences, 
	45(11), 
	pp. 7144--7161, 
	2022.
	
	\bibitem{della2025modelling} 
	R. Della Marca, M. Menale,
	\textit{Modelling the impact of opinion flexibility on the vaccination choices during epidemics},
  	Ricerche di Matematica,
 	74(2),
  	pp. 1003--1020,
	2025.
	
	\bibitem{inferrera2024reaction} 
	G.~Inferrera, C.F.~Munafò, F.~Oliveri, and P.~Rogolino, 
	\textit{Reaction-diffusion models of crimo--taxis in a street}, 
	Applied Mathematics and Computation, 
	67, pp.128504, 
	2024.
	
	\bibitem{petrovskii2020modelling} 
	S. Petrovskii, W. Alharbi, A. Alhomairi, A. Morozov,
	\textit{Modelling population dynamics of social protests in time and space: the reaction-diffusion approach},
  	Mathematics,
 	8(1),
  	pp. 78,
	2020.
	
	\bibitem{wen2009global} 
	Z. Wen, S. Fu, 
	\textit{Global solutions to a class of multi-species reaction-diffusion systems with cross-diffusions arising in population dynamics},
  	Journal of computational and applied mathematics,
 	230(1),
  	pp. 34--43,
	2009.
	
	\bibitem{ko2011stationary}
	W. Ko and I. Ahn,
	\textit{Stationary patterns and stability in a tumor-immune interaction model with immunotherapy},
	Journal of mathematical analysis and applications,
 	383(2),
 	pp. 307--329,
	2011.
	
	\bibitem{suddin2021reaction}
	S. Suddin,  F. Adi-Kusumo, L. Aryati and Gunardi,
	\textit{Reaction‐Diffusion on a Spatial Mathematical Model of Cancer Immunotherapy with Effector Cells and IL‐2 Compounds’ Interactions}, 
	International Journal of Differential Equations, 
	2021(1),
	pp. 5535447, 
	2021.
	
	\bibitem{brennan2025pattern}
  	M. Brennan and A. Krause and E. Villar-Sep{\'u}lveda and C. B. Prior,
	\textit{Pattern Formation as a Resilience Mechanism in Cancer Immunotherapy},
	Bulletin of Mathematical Biology,
 	87(8),
	pp. 1--31,
  	2025.

	\bibitem{oluwatosin2022investigating}
	O. Oluwatosin and A.K. Ibrahim and N. Hussaini,
  	\textit{Investigating Turing patterns in cancer-immune cells interaction model},
  	Bayero Journal of Pure and Applied Sciences,
  	13(1),
  	pp. 439--446,
  	2022.
	
	
	\bibitem{othmer1969interactions}
	H.G. Othmer, L.E. Scriven,
	\textit{Interactions of reaction and diffusion in open systems},
	Industrial $\&$ Engineering Chemistry Fundamentals,
 	8(2),
 	pp. 302--313,
  	1969.
		
	\bibitem{othmer1971instability}
	H.G. Othmer, L.E. Scriven,
	\textit{Instability and dynamic pattern in cellular networks},
	Journal of theoretical biology,
	32(3),
  	pp. 507--537,
	1971.
	
	\bibitem{othmer1974non}
	H.G. Othmer, L.E. Scriven,
	\textit{Non-linear aspects of dynamic pattern in cellular networks},
  	Journal of theoretical biology,
 	43(1),
 	pp. 83--112,
  	1974.
  	
	
	
	\bibitem{nakao2010turing}
	H. Nakao, A.S. Mikhailov,
	\textit{Turing patterns in network-organized activator--inhibitor systems},
  	Nature Physics,
  	6(7),
  	pp. 544--550,
  	2010.
		
	\bibitem{hata2014dispersal}
	S. Hata, H. Nakao, A. S. Mikhailov, 	
	\textit{Dispersal-induced destabilization of metapopulations and oscillatory Turing patterns in ecological networks},
  	Scientific reports,
  	4(1),
  	pp. 3585,
  	2014.
	
	\bibitem{sarker2025spatial}
	R.C. Sarker, S.K.Sahani,
	\textit{Spatial patterns in a cancer network mathematical model with subdiffusion},
  	Nonlinear Science,
  	3,
  	pp. 100026,
	2025.
		

	\bibitem{vanag2009cross}
	V.K.~Vanag and I.R.~Epstein, 
	\textit{Cross-diffusion and pattern formation in reaction--diffusion systems},
	Physical Chemistry Chemical Physics, 
	11(6),
	pp. 897--912, 
	2009.

	\bibitem{xie2012cross}
	Z.~Xie, 
	\textit{Cross-diffusion induced Turing instability for a three species food chain model},
	Journal of Mathematical Analysis and Applications, 
	388(1),
	pp. 539--547, 
	2012.

	\bibitem{hao2020spatial}
	W.~Hao and C.~Xue, 
	\textit{Spatial pattern formation in reaction--diffusion models: a computational approach},
	Journal of Mathematical Biology, 
	80(1), 
	pp. 521--543,
	2020.

	\bibitem{zincenko2021turing}
	A.~Zincenko, S.~Petrovskii, V.~Volpert, and M.~Banerjee,  
	\textit{Turing instability in an economic--demographic dynamical system may lead to pattern formation on a geographical scale},  
	Journal of the Royal Society Interface, 
	18(177), 
	pp. 20210034, 
	2021.

	\bibitem{chakraborty2021diffusion}
	B.~Chakraborty, H.~Baek, and N.~Bairagi,  
	\textit{Diffusion-induced regular and chaotic patterns in a ratio-dependent predator--prey model with fear factor and prey refuge}, 
	Chaos: An Interdisciplinary Journal of Nonlinear Science, 
	31(3),
	2021.

	\bibitem{aymard2022pattern}
	B.~Aymard,  
	\textit{On pattern formation in reaction--diffusion systems containing self-and cross-diffusion}, 
	Communications in Nonlinear Science and Numerical Simulation, 
	105, 
	pp. 106090, 
	2022.

	\bibitem{MitGri}
	A.~R.~Mitchell and D.~F.~Griffiths, 
	\textit{The finite difference method in partial differential equations}, 
	John Wiley $\&$ Sons, 
	New York, 
	1980.
	
\end{thebibliography}
\end{document}